\RequirePackage{fix-cm}
\documentclass[envcountsect,numbook]{svjour3_arxiv}   
\smartqed  
\usepackage{graphicx}
\graphicspath{{./figures/}}
\usepackage[utf8]{inputenc}
\usepackage{dsfont}
\usepackage{subfigure}
\usepackage{color}
\usepackage{url}
\usepackage{enumerate} 
\usepackage{booktabs}
\usepackage{tabularx}
\usepackage{ltablex}
\usepackage{mathptmx}
\usepackage{amssymb,amsmath,amsfonts,latexsym} 
\usepackage{longtable} 
\usepackage{textcomp}
\usepackage{lscape}
\usepackage[12pt]{extsizes}
\usepackage{geometry}
\usepackage{ifthen}
\usepackage{mathtools}
\mathtoolsset{showonlyrefs}
\usepackage{calrsfs}
\DeclareMathAlphabet{\pazocal}{OMS}{zplm}{m}{n}
\let\mathcal\pazocal
\allowdisplaybreaks

\spnewtheorem{theorem}{Theorem}{\bfseries}{\itshape}
\spnewtheorem{corollary}[theorem]{Corollary}{\bfseries}{\itshape}
\spnewtheorem{lemma}[theorem]{Lemma}{\bfseries}{\itshape}
\spnewtheorem{proposition}[theorem]{Proposition}{\bfseries}{\itshape}
\spnewtheorem{definition}[theorem]{Definition}{\bfseries}{\itshape}
\spnewtheorem{remark}[theorem]{Remark}{\bfseries}{\upshape}
\spnewtheorem{assumption}[theorem]{Assumption}{\bfseries}{\itshape}
\counterwithin{theorem}{section}
\smartqed

\renewcommand{\paragraph}[1]{{\bf #1.}}
\definecolor{myred}{rgb}{0.8,0,0}  

{\vskip\baselineskip\noindent\textbf{Proof of {#1}:}}%
{\hspace*{.1pt}\hspace*{\fill}\BOX\vskip\baselineskip}

\usepackage{tikz}

\def \R{\mathbb{R}}               
\def \N{\mathbb{N}}               
\def \1{{\bf 1}}                
\def \0{{\bf 0}}

\definecolor{myred}{rgb}{0.9,0,0}  
\definecolor{mygreen}{rgb}{0,0.7,0}  
\definecolor{myblue}{rgb}{0.2,0,0.8}  
\definecolor{orange}{rgb}{1,0.6,0}  
\definecolor{olive}{rgb}{0.5,0.5,0}  
\definecolor{mylila}{rgb}{0.8,0.5,0.2}  
\definecolor{mygrey}{rgb}{0.6,0.6,0.6}  
\definecolor{mybrown}{rgb}{0.65,0.16,0.16}  
\definecolor{mymaroon}{rgb}{0.11,0.0,0.0}

\def \Domainspace{\mathcal{D}}  

\def \omatrix{C}

\def \phx{\text{PHX} }
\def \phxk{\text{PHX}}
\def \phxs{\text{PHXs} }
\def \phxsk{\text{PHXs}}

\def \medium{M}
\def \fluid{F}
\def \outlet{O}
\def \inlet{I}
\def \bottom{B}
\def \interface{J}
\def \storage{S}
\def \Qav{{\overline{Q}}{}}
\def \Qm{\Qav^\medium}
\def \Qf{\Qav^\fluid}
\def \Qmf{ \Qav^{\storage}}
\def \Qout{\Qav^{\outlet}}
\def \Qin{Q^{\inlet}}
\def \QinC{\Qin_C}
\def \QinD{\Qin_D}
\def \Qbottom{\Qav^\bottom}

\def \Qg{Q^{G}}
\def \Qmm{Q^\medium}
\def \Qff{Q^\fluid}
\def \Dm{\Domainspace^\medium}
\def \Df{\Domainspace^\fluid}

\def \Din{\Domainspace^{I}}
\def \Dout{\Domainspace^{\outlet}}
\def \Dbottom{\Domainspace^\bottom}
\def \Dtop{\Domainspace^{T}}
\def \Dleft{\Domainspace^{L}}	
\def \Dright{\Domainspace^{R}}
\def \DInterface{\Domainspace^{\interface}}		
\def \DInterfaceL{\underline{\Domainspace}^{\interface}}
\def \DInterfaceU{{\overline{\Domainspace}}{}^{\interface}}

\def \Rp{R^P}
\def \Rb{R^B}
\def \Gm{G^\medium}
\def \Gf{G^\fluid}
\def \Gmf{G^{\storage}}
\def \Gp{G^P}
\def \Gb{G^\bottom}
\def \Cav{\omatrix}
\def \OutputM{\Cav^{\medium}}
\def \OutputF{\Cav^{\fluid}}
\def \OutputMF{ \Cav^{\storage}}
\def \OutputOut{\Cav^{\outlet}}
\def \OutputBottom{\Cav^{\bottom}}

\def \rhom{\rho^\medium}
\def \rhof{\rho^\fluid}
\def \kappam{\kappa^\medium}
\def \kappaf{\kappa^\fluid}
\def \cp{c_p}
\def \cpm{\cp^\medium}
\def \cpf{\cp^\fluid}
\def \am{a^\medium}
\def \af{a^\fluid}

\def \betam {\beta^\medium}

\newcommand{\mycaption}[1]{\caption{\footnotesize #1}}
\newcommand{\ncol}{q}
\newcommand{\vconst}{\overline{v}_0}

\newcommand{\heattransfer}{\lambda^{\!G}}
\newcommand{\Celsius}{{\text{\textdegree C}}}

\newcommand{\mat}[1]{{#1}}

\newcommand{\normalvec}{\mathfrak{n}}   
\newcommand{\iu}{{\underline{i}}}
\newcommand{\ju}{{\underline{j}}}
\newcommand{\io}{{\overline{i}}}
\newcommand{\jo}{{\overline{j}}}
\newcommand{\dom}{\dagger}

\begin{document}
	
	\title{Short-Term Behavior of a Geothermal Energy Storage: Numerical  Applications 
	%
}


\author{Paul Honore Takam  \and Ralf Wunderlich \and Olivier Menoukeu Pamen}

\authorrunning{P.H.~Takam, R.~Wunderlich, O.~Menoukeu Pamen} 

\institute{Paul Honore Takam \at
	Brandenburg University of Technology Cottbus-Senftenberg, Institute of Mathematics, P.O. Box 101344, 03013 Cottbus, Germany;  
	\email{\texttt{takam@b-tu.de}}           
	\and
	Ralf Wunderlich \at
	Brandenburg University of Technology Cottbus-Senftenberg, Institute of Mathematics, P.O. Box 101344, 03013 Cottbus, Germany;  
	\email{\texttt{ralf.wunderlich@b-tu.de}} 
	\and
	Olivier Menoukeu Pamen \at 
	University of Liverpool, Department of Mathematical Sciences, Liverpool L69 3BX, United Kingdom; 
	\email{\texttt{O.Menoukeu-Pamen@liverpool.ac.uk}} 
}

\date{Version of  \today}

\maketitle

\begin{abstract}
	This paper is devoted to numerical simulations of the short-term behavior  of the spatial temperature distribution in a geothermal energy storage.  Such simulations are needed for the optimal control and management  of residential heating systems equipped with an underground thermal storage.   We apply  numerical methods  derived in our companion paper \cite{Takam2021TheoResults} in which we study  the governing initial boundary value problem for a linear heat equation  with convection. Further, we perform extensive numerical experiments in order to investigate properties of  the spatio-temporal  temperature distribution and of its aggregated characteristics.

	%
\end{abstract}	
\section{Introduction}
\begin{figure}[!h]
	\centering
	\includegraphics[width=0.9\textwidth]{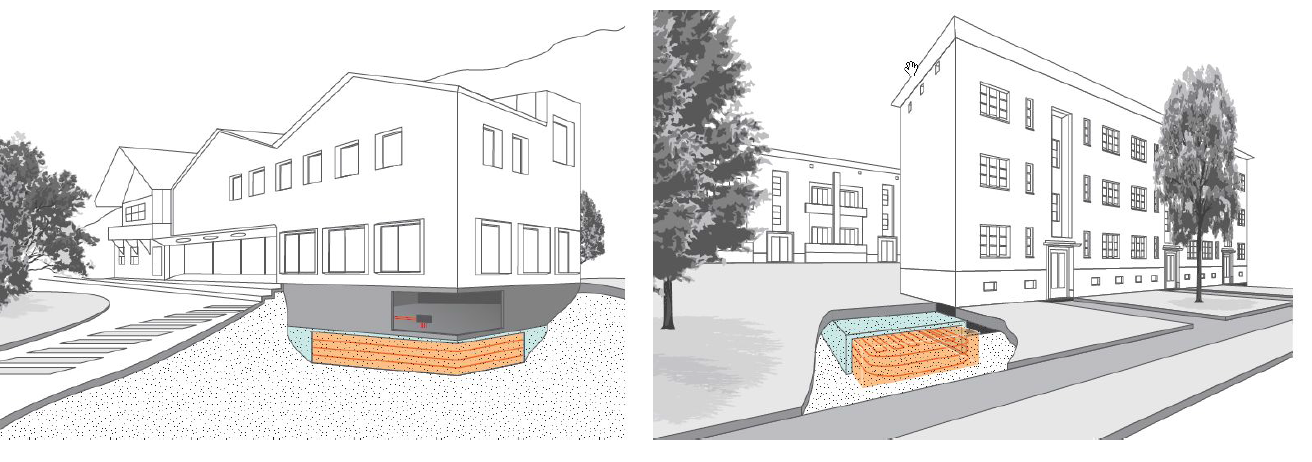}
	\caption[Geothermal storage]{Geothermal storage: in the new building, under a building (left) and in the renovation, aside of  the building (right), see
		\url{www.ezeit-ingenieure.eu}, \url{www.geo-ec.de}.}
	\label{fig:etank}
\end{figure}

This paper is devoted to the computation  of the spatial temperature distribution in a geothermal energy storage for short periods of time. We focus on underground
thermal storages as depicted in  Fig.~\ref{fig:etank} which can be found in heating systems of single buildings as well as of district heating systems.
Such storages have gained more and more importance and are quite attractive for residential heating systems since  construction and maintenance  are relatively inexpensive. Furthermore, they can be integrated both in new buildings and in renovations.   	
Such facilities are 	
used to  mitigate and to manage temporal fluctuations of   heat  supply and demand  and to move heat demand through time.  
It is well-known that thermal storages can significantly increase both the flexibility and the performance of district energy systems and enhancing the integration of intermittent  renewable energy 	sources into thermal networks (see Guelpa and  Verda \cite{guelpa2019thermal}, Kitapbayev et al.~\cite{KITAPBAYEV2015823}).     Since heat production is still  mainly based on burning fossil fuels (gas, oil, coal) these are important contributions for the reduction of carbon emissions and an increasing energy independence of societies.

The efficient operation of  geothermal storages requires a thorough design and planning because of the considerable  investment cost. For that purpose,   mathematical models and  numerical  simulations  are widely used. We refer to Dahash et al. \cite{dahash2020toward} and the references therein. In that paper the authors investigate  large-scale seasonal thermal energy storages allowing for buffering intermittent renewable heat production in district heating systems.  Numerical simulations are based on a multi-physics model of the thermal energy storage which was calibrated to measured data for a pit thermal energy storage in Dronninglund (Denmark). 
Another contribution is  Major et al. \cite{major2018numerical} which considers heat storage capabilities of deep sedimentary reservoirs. The governing heat and flow equations are solved using finite element methods.  Further, Regnier et al.~\cite{Regnier_et_al_2022} study the numerical simulation of aquifer thermal energy storages and focus on dynamic mesh optimisation for the finite element solution of the heat and flow equations. For an overview on thermal energy storages we refer to  Dincer and Rosen \cite{dincer2021thermal} and for further contributions on the numerical simulation of such storages to \cite{bazri2022thermal,dincer2021thermal,haq2016simulated,li2022modelling,soltani2019comprehensive,wu2022enhancing}.

\begin{figure}[h!]
	\centering
	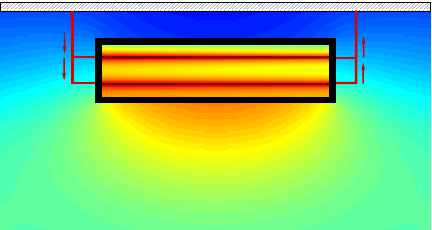
	\mycaption{\label{etank_longtermsimu} 
		2D-model of a geothermal storage  insulated to the top and the sides  while open at the bottom and spatial temperature distribution.}
\end{figure}

This paper is based on our  paper \cite{Takam2021TheoResults} where we give  a detailed description of the mathematical model of  an underground thermal energy storage and the derivation and theoretical justification of the numerical methods. The starting point is  a 2D-model, see Fig.~\ref{etank_longtermsimu}. A defined volume under or aside of a building is filled with soil and insulated to the surrounding ground.   The storage is charged and discharged via pipe heat exchangers (\phxsk) filled with some fluid (e.g.~water). 
Thermal energy is stored by raising the temperature of the soil inside the storage.
A special feature of the storage is its open architecture at the bottom. There is no insulation such that thermal energy can also flow into deeper layers   as it can be seen in Fig.~\ref{etank_longtermsimu}. This leads to a natural extension of the storage capacity since that heat can to some extent be retrieved if the storage is sufficiently discharged  (cooled) and a heat flux back to storage is induced. 

A similar model  has been already considered in B\"ahr et al.~\cite{bahr2017fast,bahr2022efficient} where the authors focus  on the numerical simulation of the long-term behavior over weeks and months of the spatial temperature distribution  and the interaction between a geothermal storage and its surrounding domain. For the sake of simplicity  the charging and discharging process using \phxs was not modeled in detail but described by a source term. In this work, we  focus on the short-term behavior of the spatial temperature distribution.  We believe that this is interesting for storages embedded into residential heating systems and the study of the storage's  response to charging and discharging operations on time scales from a few minutes to a few days. Contrary to \cite{bahr2017fast,bahr2022efficient} we include \phxs  for a more realistic model of the storage's  charging and discharging process.  However, we do not consider the surrounding medium but  reduce the computational domain to  the storage  depicted in Fig.~\ref{etank_longtermsimu} by a  black rectangle. Instead we 
set appropriate boundary conditions to mimic the interaction between storage and  environment.

The  temporal evolution  of  the spatial temperature distribution is governed by a linear heat equation with convection and appropriate boundary and interface conditions. 
A  numerical solution of that PDE  using finite difference schemes  is  sketched in Sec.~\ref{Discretization}. For more details we refer to our  paper \cite{Takam2021TheoResults}.
Management and control of a storage that is embedded into a residential heating system usually does not require the complete spatio-temporal temperature distribution   but	is  based only on certain  aggregated characteristics that can be computed in a post-processing step as  explained in Sec.~\ref{sec:Aggregate}. Examples are the average temperatures in the storage medium,  in the \phx fluid, at the outlet of the \phxs and at the storage's bottom boundary, respectively. From these quantities one can derive the  amount of available thermal energy that can be stored in or extracted from the storage in a given short period of time. 

In Sec.~\ref{Numerical} we present results of  extensive numerical experiments where we use simulations results for the temporal behavior of the spatial temperature distribution  to determine how much energy  can be stored in or taken from the storage within a given  short period of time. Special focus is laid on the dependence of these quantities on the arrangement of the $\phxs$  within the storage.

In another companion paper \cite{Takam2020Reduction}  we apply model reduction techniques known from control theory such as Lyapunov balanced truncation to derive low-dimensional approximations of the above mentioned aggregated characteristics. The latter is crucial if the cost-optimal management of residential heating systems equipped with a geothermal storage is studied mathematically in terms of optimal control problems.  It is well-known that most of the model reduction methods are developed for linear time-invariant (LTI) systems. However, the heat equation \eqref{heat_eq2} which we derive in Sec.~\ref{ExternalStorage} contains a convection term that is driven by the velocity of the fluid in the $\phxsk$. That velocity is time-dependent and typically piecewise constant during waiting, charging and discharging periods. Therefore, we are not in the framework of LTI systems and propose in  Sec.~\ref{sec:analog:system} an LTI analogous model that mimics the most important features of the original non-LTI model of the geothermal storage. 

The rest of the paper is organised as follows. In Sec.~\ref{ExternalStorage} we derive a linear heat equation with a convection term and appropriate boundary and interface conditions describing the dynamics of the spatial temperature distribution in the geothermal storage. Sec.~\ref{Discretization} is devoted to  the finite difference discretization of the heat equation. In Sec.~\ref{sec:Aggregate} we introduce aggregated characteristics of the spatio-temporal temperature distribution and explain their numerical approximation. Sec.~\ref{Numerical} presents numerical results of extensive case studies. 	We provide additional video material showing animations of the temporal evolution of the spatial temperature distribution in the geothermal storage.
The videos are available at {\small\url{www.b-tu.de/owncloud/s/D68fmqXRcgbesKj}}~.		
Finally, in Sec.~\ref{sec:analog:system} we derive  an LTI analogous model of the geothermal storage and present some numerical results. 	 
An appendix provides a list of frequently used notations and some    auxiliary results 	 removed from the main text.

\section{Dynamics of Spatial Temperature Distribution in a  Geothermal Storage}
\label{ExternalStorage}

The setting is based  on our  paper \cite[Sec.~2]{Takam2021TheoResults}. For self-containedness  and the convenience of the reader,  we recall in this section the description of the model.
The dynamics of the  spatial temperature distribution in a geothermal storage can be described mathematically by a linear heat equation with  convection term and appropriate boundary and interface conditions.  
We denote by $Q$ the temperature in the geothermal storage   depending  on time  as well as on the location in the storage.   

\subsection{2D-Model}
\label{model2D}

\begin{figure}[h!]		
	\begin{center}
		\includegraphics[width=0.8\linewidth,height=.6\linewidth]{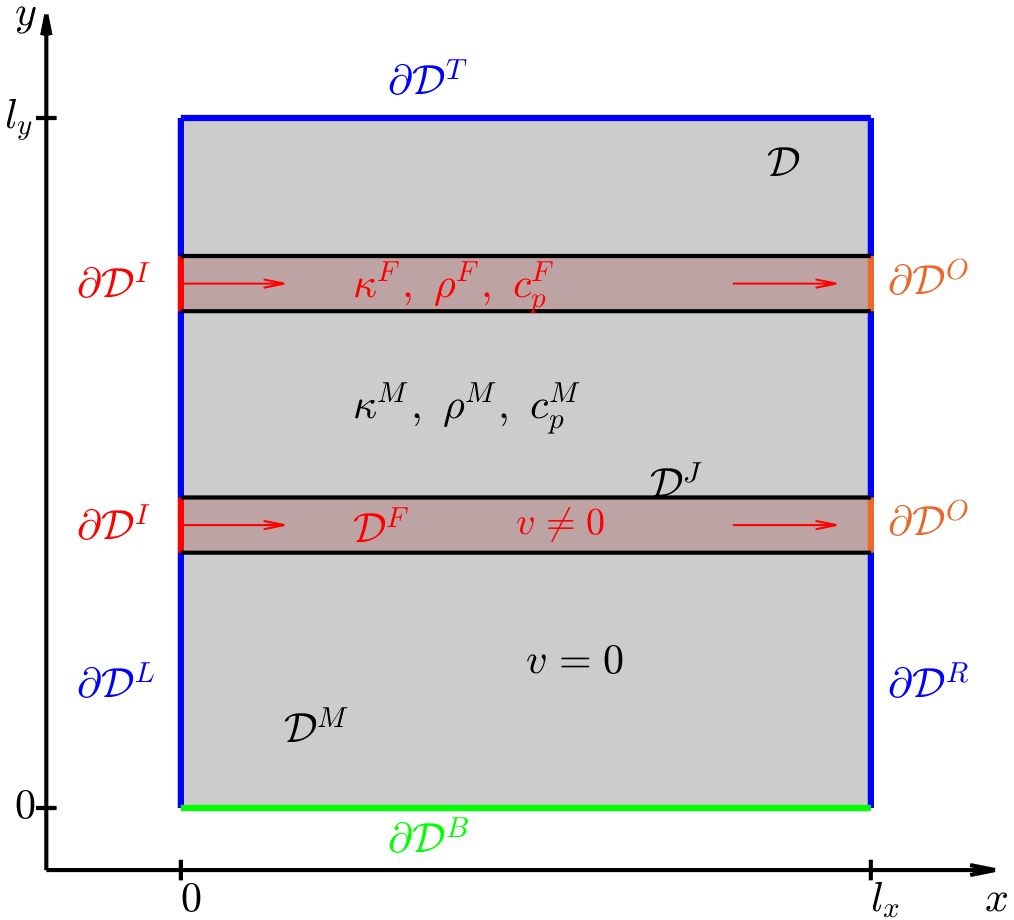}
	\end{center}
	\mycaption{\label{bound_cond} 2D-model of the geothermal storage:  decomposition of the domain $\Domainspace$ and the boundary $\partial \Domainspace $.}			
\end{figure}

We assume that the domain of the geothermal storage is a cuboid and consider a two-dimensional rectangular cross-section.
We denote by  $Q=Q(t,x,y)$ the temperature at time $t \in [0,T]$ in the point $(x,y)\in \Domainspace=(0,l_x) \times (0,l_y)$ with $l_x,l_y$ denoting the width and height of the storage. 
The domain $\Domainspace$ and its boundary $\partial \Domainspace$  are depicted  in Fig.~\ref{bound_cond}.  $\Domainspace$ is divided into three parts. The first is $\Dm$ and is  filled with a homogeneous medium (soil) characterized by constant material parameters $\rhom, \kappam$ and $\cpm$ denoting  mass density,    thermal conductivity and   specific heat capacity, respectively. The second is $\Df$,  it represents the $\phxs$  filled with a fluid (water) with constant material parameters $\rhof, \kappaf$ and $\cpf$. The fluid moves with time-dependent velocity $v_0(t)$ along the $\phxk$. For the sake of simplicity we  restrict ourselves to the case, often observed in applications, where the pumps moving the fluid are either on or off. Thus the velocity $v_0(t)$ is piecewise constant taking  values $\vconst>0$ and zero, only.  Finally, the third part is the interface $\DInterface$ between $\Dm$ and $\Df$. That interface is split into  upper and lower interfaces $\DInterfaceU$ and $\DInterfaceL$, respectively. Observe that we neglect modeling the wall of the $\phx$ and suppose perfect contact between the $\phx$ and the soil. Details are given  in \eqref{Interface} and \eqref{eq: 13f} below. The above can be summarized in the following
\begin{assumption}
	\label{assum1}~
	\begin{itemize}
		\item [1.] Material parameters of the medium   $\rhom, \kappam, \cpm$ in the domain $\Dm$  and of the fluid  $\rhof, \kappaf, \cpf$ in the domain  $\Df$ are constants.
		\item [2.] Fluid velocity is piecewise constant, i.e. $v_0(t)=\begin{cases}
			\vconst>0,  &\text{pump~on},\\
			0, & \text{pump off}.
		\end{cases}$
		\item [3.] Perfect contact at the interface between fluid and  medium.
	\end{itemize}
\end{assumption}

\begin{remark}\label{rem:3D} Results obtained for our 2D-model, where $\Domainspace$ represents the rectangular cross-section of a box-shaped storage can be extended to the 3D-case if we assume that the 3D storage domain is a cuboid of depth $l_z$ with a homogeneous temperature distribution in $z$-direction. 
	A $\phx$ in  the 2D-model then represents a horizontal snake-shaped $\phx$ densely filling a small layer of the storage. 
\end{remark}

\paragraph{Heat equation}
The temperature $Q=Q(t,x,y)$ in the external storage is governed by the linear heat equation with convection term
\begin{align}
	\rho c_p \frac{\partial Q}{\partial t}=\nabla \cdot (\kappa \nabla Q)-
	\rho v \cdot \nabla (c_p Q),\quad  (t,x,y) \in (0,T]\times \Domainspace \backslash \DInterface,   \label{heat_eq}
\end{align}
where  $\nabla=\big(\frac{\partial}{\partial x},\frac{\partial}{\partial y}\big)$ denotes the gradient operator. The  first term  on the right hand side describes diffusion while the second represents convection of the moving fluid in the \phxsk.
Further,  
$v=v(t,x,y)$ $=v_0(t)(v^x(x,y),v^y(x,y))^{\top}$  denotes  the velocity vector with $(v^x,v^y)^\top$ being the normalized directional vector of the flow. 
According to Assumption \ref{assum1} the material parameters $\rho,\kappa, c_p$ 
depend on the position $(x,y)$ and take the values $\rhom,\kappam, \cpm$ for points in $\Dm$ (medium) and  $\rhof, \kappaf, \cpf$ in  $\Df$ (fluid).

Note that there are no sources or sinks inside the storage and therefore the above heat equation appears without forcing term.
Based on this assumption, the heat equation (\ref{heat_eq}) can be written as
\begin{align}
	\frac{\partial Q}{\partial t}=a\Delta Q-
	v \cdot \nabla Q,\quad  (t,x,y) \in (0,T]\times \Domainspace \backslash \DInterface,   \label{heat_eq2}
\end{align}
where $\Delta=\frac{\partial^2}{\partial x^2}+\frac{\partial^2}{\partial y^2}$ is the Laplace operator  and  $ a=a(x,y)$ is the thermal diffusivity which is piecewise constant with values  $a^\dom=\frac{\kappa^\dom}{\rho^\dom \cp^\dom}$  with $\dom=\medium$ for  $(x,y)\in \Dm $ and $\dom=\fluid$   for  $(x,y)\in \Df$, respectively.   
The initial condition $Q(0,x,y)=Q_0(x,y)$ is given by the initial temperature distribution $Q_0$ of the storage.

\begin{remark}\label{rem_heat_exchanger}
	In real-world geothermal storages \phxs  are often designed in a snake form located in the storage domain at multiple horizontal layers. There may be  only a single inlet and a single outlet. We will mimic that design by a computationally more tractable design characterized by multiple horizontal straight \phxs as it is sketched in Fig.~\ref{etank_longtermsimu}.
	This allows to control   the \phxs in different layers separately. For a  topology with single inlet and  outlet snake-shaped \phxs  the outlet of a straight \phx in one layer  can be connected with the inlet of the straight \phx in the next layer. 
\end{remark}

\subsection{Boundary and Interface Conditions}
For the description of the boundary conditions we decompose the boundary $\partial\Domainspace$ into several subsets as depicted in Fig.~\ref{bound_cond} representing the insulation on the top and the side, the open bottom, the inlet and outlet of the $\phxsk$.  Further, we have to specify conditions at the interface between  $\phxs$  and  soil. The inlet, outlet and the interface  conditions model the heating and cooling of the storage via $\phxsk$. We distinguish between the two regimes 'pump on' and 'pump off'. For simplicity we assume perfect insulation at inlet and outlet if the pump is off.
This leads to the following boundary conditions.

\begin{itemize}				
	\item \textit{Homogeneous Neumann condition} describing perfect insulation on the top  and the side 
	\begin{align}\frac{\partial Q}{\partial \normalvec }=0,\qquad (x,y)\in 
		\partial \Dtop \cup \partial \Dleft  \cup \partial \Dright , 
		\label{Neumann}
	\end{align}
	where $\partial \Dleft =\{0\} \times  [0,l_y] \backslash \partial \Din $, ~
	$\partial \Dright =\{l_x\} \times  [0,l_y] \backslash \partial \Dout$, $\partial \Dtop = [0,l_x] \times \{l_y\}$   and $\normalvec$ denotes the outer-pointing normal vector.
	\item \textit{Robin condition} describing heat transfer at the bottom 
	\begin{align}
		-\kappam\,\frac{\partial Q}{\partial \normalvec }=\heattransfer(Q-\Qg(t)), \qquad  (x,y)\in 
		\partial \Dbottom ,
		\label{Robin}
	\end{align}
	with $\partial \Dbottom = [0,l_x] \times \{0\}$, where $\heattransfer>0$ denotes the  heat transfer coefficient  and $\Qg(t)$ the underground temperature.
	\item  \textit{Dirichlet condition} at the inlet if the pump is on ($v_0(t)>0$), i.e.~the fluid arrives at the storage with a given temperature $\Qin(t)$. If pump is off ($v_0(t)=0$), we set a homogeneous Neumann condition describing perfect insulation.
	\begin{align}
		\begin{cases}
			\begin{array}{rll}	
				Q &=\Qin(t), &\text{ ~pump on,} \\
				\frac{\partial Q}{\partial \normalvec }&=0, &\text{ ~pump off,} 
			\end{array}
		\end{cases} 
		\qquad  (x,y)\in 	\partial \Din .
		\label{input}
	\end{align}
	
	\item \textit{``Do Nothing'' condition} at the outlet in the following sense. If  the  pump is on ($v_0(t)>0$) then the total heat flux directed outwards can be decomposed into a diffusive  heat flux given by $\kappaf\frac{\partial Q}{\partial \normalvec }$ and a convective  heat flux given by $v_0(t) \rhof \cpf Q$. Since in real-world applications the latter is much larger than the first we neglect the diffusive heat flux. This leads to a homogeneous Neumann condition
	\begin{align}\frac{\partial Q}{\partial \normalvec }=0,\qquad (x,y)\in 
		\partial \Dout . 
		\label{output}
	\end{align}
	If the pump is off then we  assume (as already for the inlet) perfect insulation which is also described by the above condition.
	
	\item \textit{Smooth heat flux} at interface $\DInterface$ between fluid and soil leading to a coupling condition
	\begin{align}
		\kappaf\,\frac{\partial \Qff }{\partial \normalvec  }=\kappam\,\frac{\partial \Qmm }{\partial \normalvec  },
		\qquad  (x,y)\in \DInterface.
		\label{Interface}
	\end{align}	
	Here, $\Qff , \Qmm $ denote the temperature of the fluid inside the $\phx$ and of the soil outside the $\phxk$, respectively.
	Moreover, we assume that the contact between the $\phx$ and the medium is perfect which leads to a smooth transition of a temperature, i.e., we have 
	\begin{align}
		\Qff =\Qmm ,\qquad  (x,y)\in \DInterface. \label{eq: 13f}
	\end{align} 
\end{itemize}

\section{Discretization of the Heat Equation}
\label{Discretization}
We now sketch the discretization of the  heat equation \eqref{heat_eq2} together with the boundary and interface conditions given in \eqref{Neumann} through \eqref{eq: 13f}. For details we refer to our  paper \cite[Sec.~3 and 4]{Takam2021TheoResults}.
We proceed in two steps. In the first step we apply semi-discretization in space and approximate only spatial derivatives by their respective finite differences. This approach is also known as 'method of lines' and leads  to a high-dimensional system of ODEs  for the temperatures at the grid points.  In the second step also time is discretized  resulting in  an implicit finite difference scheme.

\subsection{Semi-Discretization of the Heat Equation}
\label{Semi-Discret}
\begin{figure}[h!]
	\centering
	\begin{center}
		\resizebox{0.7\textwidth}{0.5\textwidth}{%
			\begin{tikzpicture}[thick,scale=1, every node/.style={scale=1}]
				\draw[thick,->] (-1.5,-1.5) -- (8.7,-1.5) node[anchor=north west] {x};
				\draw[thick,->] (-1.5,-1.5) -- (-1.5,8.3) node[anchor=south east] {y};
				\draw[step=1.5cm,black,very thin] (-1.5,-1.5) grid (7.5,7.5);					
				\fill (3,3) circle (3pt);
				\node at (3.5,2.7) {$(i,j)$};
				\fill (3,4.5) circle (3pt);
				\node at (3.8,4.8) {$(i,j+1)$};
				\fill (3,1.5) circle (3pt);
				\node at (3.8,1.2) {$(i,j-1)$};
				\fill (4.5,3) circle (3pt);
				\node at (5.3,2.7) {$(i+1,j)$};
				\fill (1.5,3) circle (3pt);
				\node at (2.2,2.7) {$(i-1,j)$};
				\node at (7.7,-2) {$l_x=N_xh_x$};
				\node at (-2.6,7.5) {$l_y=N_yh_y$};
				\fill (7.5,7.5) circle (3pt);
				\node at (8.3,7.75) {$(N_x,N_y)$};
				\fill (-1.5,-1.5) circle (3pt);
				\node at (-2.1,-1.8) {$(0,0)$};	
				\fill (-1.5,7.5) circle (3pt);
				\node at (-0.8, 7.75) {$(0,N_y)$};
				\fill (7.5,-1.5) circle (3pt);
				\node at (8.2,-1.1) {$(N_x,0)$};
			\end{tikzpicture}
		} 
	\end{center}
	\mycaption{\label{grid}Computational grid.}	
\end{figure}
The spatial domain depicted in Fig.~\ref{bound_cond} is discretized by the means of a mesh   with grid points  $(x_i,y_j)$ as shown in Fig.~\ref{grid} where
$	x_i =ih_x, ~y_j =jh_y,~ i ={ 0},...,N_x, ~j ={ 0},...,N_y.$
Here,  $N_x$ and $N_y$ denote the number of grid points  while   $h_x={l_x}/{N_x}$ and $h_y={l_y}/{N_y}$  are the step sizes in $x$ and $y$-direction, respectively. 
We denote by $Q_{ij}(t)\simeq Q(t,x_i,y_j)$ the semi-discrete approximation  of the temperature   and by $v_0(t)(v^x_{ij},v^y_{ij})^{\top} =v_0(t)(v^x(x_i,y_j),v^y(x_i,y_j))^{\top} =v(t,x_i,y_j)$  the velocity vector  at the grid point $(x_i,y_j)$ at time $t$.

For the sake of simplification and tractability of our analysis we restrict ourselves to the following assumption on the arrangement of $\phxs$  and impose  conditions on the location of grid points  along the \phxsk.
\begin{assumption}\label{assum2}~%
	\begin{enumerate}
		\item  There are $n_P \in \N$  straight horizontal \phxsk, the fluid moves in positive $x$-direction.
		\item The interior of $\phxs$  contains grid points.
		\item Each interface between medium and fluid contains grid points.
	\end{enumerate}		
\end{assumption}	

We approximate the spatial derivatives in the heat equation \eqref{heat_eq2}, the boundary and interface conditions by finite differences as in \cite[Subsec.~3.1--3.3]{Takam2021TheoResults} where we apply upwind techniques for the convection terms. The result is the system of ODEs \eqref{Matrix_form1} (given below) for a vector function $Y:[0,T]\to \R^n$ collecting  the semi-discrete approximations $Q_{ij}(t)$ of the temperature $Q(t,x_i,y_j)$ in the ``inner'' grid points, i.e., all grid points except those on the boundary $\partial \Domainspace$ and the interface $\DInterface $. 
For a model with $n_P$ $\phxs$  the dimension of $Y$ is   $n=(N_x-1)(N_y-2n_P-1)$, see \cite{Takam2021TheoResults}.

Using the above notation the semi-discretized  heat equation  together with the given initial, boundary and interface conditions reads as
\begin{align}
	\frac{d Y(t)}{dt}= \mat{A}(t)Y(t)+\mat{B}(t)g(t), ~~t \in (0,T],
	\label{Matrix_form1}
\end{align}
with the initial condition $Y(0)=y_0$ where the vector $y_0\in \R^n$ contains the initial temperatures $Q(0,\cdot,\cdot)$ at the corresponding grid points.  
The system matrix  $\mat{A}$  results from the spatial discretization of the convection and diffusion term in the heat equation (\ref{heat_eq2}) together with the Robin and linear heat flux boundary conditions.  It has the tridiagonal structure  
\begin{align}
	\label{matrix_A}
	\mat{A}=\begin{pmatrix}
		\mat{A}_{L} ~&~ \mat{D}^{+} ~& &&&\text{\LARGE0} \\
		\mat{D}^{-} ~&~ \mat{A}_{M} ~&~ \mat{D}^{+} \\
		& \mat{D}^- ~&~ \mat{A}_{M} ~&~ \mat{D}^{+} \\
		&& \ddots &\ddots & \ddots &\\
		& && \mat{D}^- ~&~ \mat{A}_{M}~ &~ \mat{D}^+ \\
		\text{\LARGE0}& & &&~ \mat{D}^- ~&~ \mat{A}_{R}
	\end{pmatrix}
\end{align}
and consists of $(N_x-1)\times( N_x-1)$ block matrices of dimension $\ncol=N_y-2n_P-1$.
The block matrices $\mat{A}_{L},\mat{A}_{M},\mat{A}_{R}$ on the diagonal have a tridiagonal structure and are given
in \cite[Tables 3.1 and B.1]{Takam2021TheoResults}. 		
The block matrices on the subdiagonals $\mat{D}^{\pm} \in \R^{\ncol \times \ncol}$, $i=1,\ldots, N_x-1$, are  diagonal matrices and given in \cite[Eq.~(3.12)]{Takam2021TheoResults}. 

As a result of the discretization of the  Dirichlet  condition at the inlet boundary and the Robin condition at the bottom boundary, we get the  function  $g:~[0,T] \to \R^2$ called input function and the  $n\times 2$ input matrix $\mat{B}$  called input matrix. The entries of the input matrix  $B_{lr}, ~l=1,\ldots,n,~~ r=1,2,$ are derived in \cite[Subsec 3.4]{Takam2021TheoResults} and are given by 
\begin{align}\label{eq:input_matrix}
	\begin{array}{rl@{\hspace*{2em}}l}
		B_{l1}&=B_{l1}(t)=\begin{cases}
			\frac{\af }{h^2_x}+\frac{\vconst}{h_x},  & \text{pump on,}\\
			0, & \text{pump off,}
		\end{cases}
		& l=\mathcal{K}(1,j), (x_0,y_j)\in\Din,\\[3ex]	
		B_{l2}& = \frac{\heattransfer h_y}{\kappam+\heattransfer h_y}\betam ,
		&  l=\mathcal{K}(i,1), (x_i,y_0)\in \Dbottom ,
	\end{array}		
\end{align}   
with $\beta^\medium=\am/{h^2_y}$.
The entries for other $l$ are zero. Here,  $\mathcal{K}$ denotes the mapping   $(i,j)\mapsto  l=\mathcal{K}(i,j)$  of pairs of indices   of  grid point $(x_i,y_j) \in\Domainspace$ to the single index $l \in \{1,\ldots,n \}$ of the corresponding entry in the vector $Y$. The input function  reads as
\begin{align}
	g(t)=\begin{cases}
		(\Qin(t),~\Qg(t))^{\top}, & \quad \text{pump on},\\
		~~~~(0,~~~~\Qg(t))^{\top}, & \quad \text{pump off}.
	\end{cases}
	\label{eq:input}
\end{align}
Recall that  $\Qin$ is the inlet temperature  of the $\phx$ during pumping  and $\Qg$ is the underground temperature.

\subsection{Full Discretization}
\label{sec:full_discrtization}
After discretizing the heat equation \eqref{heat_eq2} w.r.t.~spatial variables  we will now also discretize the temporal derivative and derive a family of implicit finite difference schemes.

We introduce the notation  $N_{\tau}$ for the number of grid points in  $t$-direction,   $\tau={T}/{N_{\tau}}$ the time step and $t_k=k \tau, ~ k = 
0,...,N_{\tau}$. Let $\mat{A}^k,\mat{B}^k, g^k, v_0^k$ be the values of $\mat{A}(t),\mat{B}(t), g(t), v_0(t)$ at time  $t=t_k$. Further, we denote by 
$ Y^k=(Y^k_1,\ldots,Y^k_n)^\top$ the discrete-time approximation of the vector function $Y(t)$ at time  $t=t_k$. 
Discretizing the temporal derivative in \eqref{Matrix_form1} with the forward difference gives 
\begin{align}
	\frac{d Y(t_k)}{d t} = \frac{Y^{k+1}-Y^k}{\tau}+ \mathcal{O}(\tau).
	\label{time_d}
\end{align} 

Substituting \eqref{time_d} into \eqref{Matrix_form1} and replacing the r.h.s.~of \eqref{Matrix_form1} by a convex combination of the values at time $t_k$ and $t_{k+1}$ with the weight  $\theta \in [0,1]$ gives the following general $\theta$-implicit finite difference scheme 
\begin{align}
	\frac{Y^{k+1}-Y^k}{\tau} = \theta[\mat{A}(t_{k+1})Y^{k+1} +\mat{B}(t_{k+1}) g^{k+1}] +(1-\theta)[\mat{A}(t_{k})Y^{k} +\mat{B}(t_{k}) g^{k}]
	\nonumber
\end{align}
for which we provide in our  paper \cite[Sec.~4]{Takam2021TheoResults} a detailed stability analysis. For our numerical experiments in Sec.~\ref{Numerical} we use an explicit scheme which is obtained for  $\theta=0$ and given by the recursion as
\begin{align}
	Y^{k+1}&= (\mathds{I}_{n}+\tau\mat{A}^k)Y^k +\tau \mat{B}^kg^k, \quad k=0,\ldots,N_\tau-1, 
	\label{Matrix_explicit} 
\end{align}
with the initial value $Y^0=Y(0)$ and the notation $\mathds{I}_{n}$ is the  $n \times n$ identity matrix.
The advantage of an explicit scheme  is that it avoids the time-consuming solution of systems of linear equations but one has to satisfy stronger conditions on the time step $\tau$ to ensure stability of the scheme. In \cite[Theorem 4.2]{Takam2021TheoResults}, we show that the above explicit scheme is stable  if the time step $\tau$ satisfies the condition 
\begin{align}
	\nonumber 
	\tau \leq 	\Big(2 \max\{\af ,\am \}\Big(\frac{1}{h_x^2}+\frac{1}{h_y^2}\Big)+\frac{\vconst }{h_x}\Big)^{-1}.
\end{align}

\section{Aggregated Characteristics}	
\label{sec:Aggregate}
The numerical methods introduced in Sec.~\ref{Discretization}	 allow the approximate  computation of  the  spatio-temporal temperature distribution in the geothermal storage. In many applications it is  not necessary to know  the complete information about that distribution. An example is the management and control of a storage which is embedded into a residential heating system. Here it is sufficient to know  only  the response of  a few aggregated characteristics of the temperature distribution  to charging and discharging operations. These quantities  can be computed via  a post-processing procedure.   In this section we introduce some of these aggregated characteristics and describe their approximate computation  based on the solution vector $Y$ of the finite difference scheme.

\subsection{Aggregated Characteristics Related to the Amount of Stored Energy}

We start with aggregated characteristics given by the average temperature in some  subdomain of the storage which are related to the amount of stored energy in that domain. 

Let  $\mathcal{B}\subset \Domainspace$ be a  generic  subset of the 2D computational domain. We denote by $|\mathcal{B}|=\iint_{\mathcal{B}} dxdy$ the area of $\mathcal{B}$. 
Then
$W_{\mathcal{B}}(t)=l_z \iint_{\mathcal{B}} \rho c_p Q(t,x,y) dxdy$ represents the thermal energy contained in the 3D spatial domain $\mathcal{B}\times[0,l_z]$ at time $t\in[0,T]$. Then for $0\le t_0<t_1\le T$ the difference $G_{\mathcal{B}}(t_0,t_1)=W_{\mathcal{B}}(t_1)-W_{\mathcal{B}}(t_0)$ is the gain of thermal energy during the period $[t_0,t_1]$. While positive values correspond to warming of $\mathcal{B}$,  negative values indicate cooling and  $-G_{\mathcal{B}}(t_0,t_1)$ represents the magnitude of the  loss of thermal energy. 

For $\mathcal{B}=\Domainspace^\dom, \dom=\medium,\fluid$, we can use that the material parameters on $\Domainspace^\dom$ equal  the constants  $\rho=\rho^\dom,c_p=c_p^\dom$. Thus,  for the corresponding gain of thermal energy we obtain
\begin{align*}
	G^\dom=G^\dom(t_0,t_1)& := G_{\Domainspace^\dom}(t_0,t_1) = \rho^\dom c_p^\dom|\Domainspace^\dom| l_z~(\Qav^\dom(t_1)-\Qav^\dom(t_0)),\\
	\text{where}\quad 
	\Qav^\dom(t) &= \frac{1}{|\Domainspace^\dom|} \iint_{\Domainspace^\dom} Q(t,x,y) dxdy, \quad \dom=\medium,\fluid,
\end{align*}
denotes the  average temperature in the medium ($\dom=\medium$) and the fluid  ($\dom=\fluid$), respectively.
We denote by $\Qmf$ the average temperature  in the whole storage. It  can be obtained from $\Qm $ and $\Qf$ by
\begin{align}
	\label{av_temp_storage}
	\Qmf(t)= \frac{1}{|\Domainspace|} \big(\Qm (t)\,|\Dm | +\Qf(t)\,|\Df|\big).
\end{align}
Further, the total gain in the storage denoted by $\Gmf$ is obtained by 
\begin{align*}
	\Gmf= \Gmf(t_0,t_1)=\Gm(t_0,t_1)+\Gf(t_0,t_1).
\end{align*}

\subsection{Aggregated Characteristics Related to the Heat Flux at the Boundary}
Now we consider the convective heat flux at the inlet and outlet boundary and the heat transfer at the bottom boundary. 
Let  $\mathcal{C}\subset \partial \Domainspace $ be a generic curve on the boundary,  then  we denote by $|\mathcal{C}|=\int_{\mathcal{C}} ds$ the curve length. 

The rate  at which the energy is injected or withdrawn via the $\phx$ is given by 
\begin{align}
	\nonumber	
	\Rp(t)&=\rho^\dom c_p^\dom v_0(t) \Big[\int_{\Din}  Q(t,x,y)\, ds  -\int_{\Dout} Q(t,x,y)\, ds\Big]\\
	\label{Rout}	
	&=\rho^\dom c_p^\dom v_0(t) |\partial  \Dout|[\Qin(t) -\Qout(t)],\\
	\nonumber	
	\text{where}\quad 
	\Qout(t) &=	 \frac{1}{|\partial \Dout|} \int_{\partial \Dout} Q(t,x,y)ds
\end{align}
is the average temperature at the outlet boundary. Here, we have used that in our model we have  horizontal $\phxs$  such that  $|\partial \Din|=|\partial \Dout|$ and a uniformly distributed  inlet temperature at the inlet boundary $\partial \Din$.   Note that the fluid moves at time $t$ with velocity $v_0(t)$ and arrives at the inlet with temperature $\Qin(t)$ while it leaves at the outlet with  the average temperature $\Qout(t)$. For a given interval of time $[t_0,t_1]$ the quantity 
$$\Gp=\Gp(t_0,t_1) =l_z\int_{t_0}^{t_1} \Rp(t)\, dt$$ 
describes the amount of heat injected  ($\Gp>0$)  to or withdrawn ($\Gp<0$) from the storage due to convection of the   fluid.

Next we look at the  diffusive heat transfer via the bottom boundary and define the rate 
\begin{align}
	\nonumber
	\Rb(t)&=  \int_{ \Dbottom }\kappam\frac{\partial Q}{\partial \normalvec }\, ds = \int_{ \Dbottom }\heattransfer(\Qg(t)-Q(t,x,y))\, ds	\\[0.5ex]
	\label{RB}		
	& = \heattransfer |\partial   \Dbottom | (\Qg(t)- \Qbottom(t) ),\\
	\nonumber	
	\text{where}\quad 
	\Qbottom(t) &=	 \frac{1}{|\partial  \Dbottom |} \int_{\partial  \Dbottom } Q(t,x,y)ds
\end{align}
is the average temperature at the bottom boundary.
Note that the second equation in the first line  follows from the Robin boundary condition.
The quantity $$\Gb=\Gb(t_0,t_1) =l_z\int_{t_0}^{t_1} \Rb(t)\, dt$$ describes the amount of  heat  transferred via the bottom boundary of the storage. 

\subsection{Energy Balance}
\label{subsec:energy_balance}
In our model we assume perfect thermal insulation at all boundaries except the inlet, outlet  and the bottom boundary. At the outlet we impose a  homogeneous Neumann condition describing zero diffusive heat transfer.  At the inlet we also have a zero diffusive heat transfer under the reasonable assumption that  the temperature in the supply pipe is constant and equals $\Qin(t)$, thus the normal derivative $\frac{\partial Q}{\partial \normalvec }$ is zero. This implies that gains and losses of thermal energy in the storage are caused either by injections or withdrawals via the $\phxs$   or by heat transfer via the open bottom boundary. Thus, we  can decompose the total  gain $ \Gmf$ to obtain  the following  energy balance
\begin{align}
	\label{energy_balance}
	\Gmf&= \Gm+\Gf=\Gp+\Gb. 
\end{align}

\subsection{Numerical  Computation of Aggregated Characteristics}
\label{Aggregate:Num}
In this subsection we consider the approximate computation of aggregated characteristics introduced in the previous subsections by using finite difference  approximations of the temperature $Q=Q(t,x,y)$. The approximations are given in terms of  the entries of the vector function $Y(t)$ satisfying  the system of ODEs \eqref{Matrix_form1} and containing the semi-discrete finite difference approximations of the temperature in the inner grid points of the computational domain $\Domainspace$. Recall that  the temperatures  at boundary and interface grid points can be determined by  linear combinations from the entries of $Y(t)$. The extension to approximations based of the solution of the fully discretized PDE \eqref{Matrix_explicit} is straightforward using the relation $Y(t_k)=Y(k\tau)\approx Y^k, k=0,\ldots,N_\tau$.

Let us start with the average temperatures $\Qm $ and $\Qf$, where the temperature $Q(t,x,y)$ is averaged over unions of disjoint rectangular subsets  of the computational domain $\Domainspace$.
Assume that $\mathcal{B}\subset \Domainspace$ is a generic rectangular subset with corners defined by the grid points $(x_i,y_j)$ with indices $(\iu,\ju), (\io,\ju),(\io,\jo),(\io,\ju)$,  where $0\le \iu<\io\le N_x$ and $0\le \ju<\jo\le N_y$. We assume further that  the  domain $\mathcal{B}$ contains at least one layer of horizontal and vertical inner grid points, respectively. Thus we require   $ \iu+2\le\io$ and    $\ju+2\le \jo$. We denote by $\Qav^{\mathcal{B}}=\Qav^{\mathcal{B}}(t)=\frac{1}{|\mathcal{B}|}\iint_{\mathcal{B}} Q(t,x,y)dxdy$ the average temperature in ${\mathcal{B}}$. Rewriting  the double integral as two iterated single integrals and applying trapezoidal rule to the single integrals    
the average temperature $\Qav^{\mathcal{B}}$  can be approximated by (for details see Appendix \ref{append:quadformula})
\begin{align}
	\label{quad_2D_0}
	\Qav^{\mathcal{B}} & =\frac{1}{|\mathcal{B}|}\iint_{\mathcal{B}} Q(t,x,y)dxdy \approx \sum_{(i,j)\in \mathcal{N}_\mathcal{B}} \mu_{ij}\,Q_{ij}, 
\end{align}
where $\mathcal{N}_\mathcal{B}=\{ (i,j): i=\iu,\ldots, \io, j=\ju,\ldots,\jo\}$  and the coefficients $d_{ij}$ of the above quadrature formula are given by
\begin{align}\label{mu_coeff}
	\mu_{ij} &= \frac{1}{(\io-\iu)(\jo-\ju)}  \left\{
	\begin{array}{cl@{\hspace*{1em}}l}
		1, &  \text{for }~~\iu<i<\io, ~~\ju<j<\jo, & \text{(inner points)}\\[1ex]
		\frac{1}{2}, & \text{for }~~
		\begin{array}[t]{ll}
			i\;=\iu,\io, &\ju<j<\jo,\\
			j=\ju,\jo,& ~\iu<i<\io, 
		\end{array}
		&  \text{(boundary points, except corners)}\\[1ex]
		\frac{1}{4}, & \text{for }~~i=\iu,\io, ~~j=\ju,\jo& \text{(corner points)}.
	\end{array} \right.
\end{align}

Next we want to rewrite approximation \eqref{quad_2D_0} in terms of the vector $Y=Y(t)$. Recall that $Y$ contains the finite difference approximations of the temperature in the inner grid points of the computational domain $\Domainspace$. Let us introduce the vector  $\overline Y$  of dimension $\overline n=(N_x+1)(N_y+1)-n$  containing the  temperature approximations at the remaining grid points located on the boundary $\partial \Domainspace$ and the interface $\DInterface $.  These values can be determined by  the discretized boundary and interface conditions and expressed  as linear combinations of the entries of $Y$. This allows for a  representation $\overline Y=\overline C Y$ with some $\overline n\times n-$matrix $\overline C$.

Now, let $\mathcal{N}_\mathcal{B}^0\subset \mathcal{N}_\mathcal{B}$ and $\overline{\mathcal{N}_\mathcal{B}^0}= \mathcal{N}_\mathcal{B}\setminus \mathcal{N}_\mathcal{B}^0 $ be the subsets (of index pairs $(i,j)\in \mathcal{N}_\mathcal{B}$ of grid points) for which the finite difference approximation  $Q_{ij}$ is contained in the vector $Y$ and the  vector $\overline Y$, respectively.   Further, let $\mathcal{K}:\mathcal{N}_\mathcal{B}^0\to \{1,\ldots,n \}$ and $\overline{\mathcal{K}}:\overline{\mathcal{N}_\mathcal{B}^0}\to \{1,\ldots, \overline n \}$ denote the mappings $(i,j)\mapsto  l=\mathcal{K}(i,j)$ and $(i,j)\mapsto  \overline l=\overline{\mathcal{K}}(i,j)$
of pairs of indices $(i,j)$ to the single indices  $l$ and $\overline l$ of the corresponding entries in the vectors $Y$ and $\overline Y$, respectively.
Then it holds
\[ Q_{ij}= \begin{cases}
	Y_{\mathcal{K}(i,j)}, & (i,j) \in \mathcal{N}_\mathcal{B}^0,\\
	\overline{Y}_{\overline{\mathcal{K}}(i,j)}, & (i,j) \in \overline{\mathcal{N}_\mathcal{B}^0},
\end{cases}\]
and we can rewrite approximation \eqref{quad_2D_0} as 
\begin{align}
	\label{quad_2D_1}
	\begin{array}{rcccc}	\Qav^{\mathcal{B}} & \approx &\sum\limits_{(i,j)\in \mathcal{N}_\mathcal{B}^0} \mu_{ij}\,Q_{ij} &+&\sum\limits_{(i,j)\in \overline{\mathcal{N}_\mathcal{B}^0}} \mu_{ij}\,Q_{ij} \\[0.5ex]
		& =& \sum\limits_{l=\mathcal{K}(i,j):  (i,j)\in \mathcal{N}_\mathcal{B}^0} d_{l}\,Y_{l} &+&\sum\limits_{\overline{l}=\overline{\mathcal{K}}(i,j): (i,j)\in \overline{\mathcal{N}_\mathcal{B}^0}} \overline{d}_{\overline{l}}\;\overline{Y}_{\overline{l}}\\[0.5ex]
		&=& D\,Y&+&\overline{D}\,\overline{Y},
	\end{array}
\end{align}
with an  $1\times n-$matrix $D$ and an  $1\times \overline{n}-$matrix $\overline{D}$, whose entries are given for $l=1,\ldots,n,~\overline{l}=1,\ldots, \overline{n}$ by 
\begin{align}
	\label{matrix_D}
	d_l=\begin{cases}
		\mu_{ij}, & l=\mathcal{K}(i,j),~~ (i,j)\in \mathcal{N}_\mathcal{B}^0,\\
		0 & \text{else},
	\end{cases} 
	\quad\text{and}\quad 
	\overline{d}_{\overline{l}}=\begin{cases}
		\mu_{ij}, & \overline{l}=\overline{\mathcal{K}}(i,j),~~ (i,j)\in \overline{\mathcal{N}_\mathcal{B}^0},\\
		0 & \text{else},
	\end{cases}
\end{align}
respectively.
Finally, substituting $\overline Y=\overline C Y$ into \eqref{quad_2D_1} yields a representation of the average temperature $\Qav^{\mathcal{B}}$ as a linear combination of entries of the vector $Y$ which reads as
\begin{align}
	\label{quad_2D}
	\Qav^{\mathcal{B}} & \approx    \Cav^{\mathcal{B}} \,Y \quad \text{with}\quad \Cav^{\mathcal{B}}=D+\overline{D}\,\overline{C}.
\end{align}
Based on the above representation we can derive similar approximations for the average temperatures $\Qm $ and $\Qf$ in the medium and the fluid, respectively.
Our model assumptions imply that for a storage with $n_P$ $\phxs$  the domain  $\Df$ splits into $n_P$ disjoint rectangular subsets  $\Df_j, j=1,\ldots,n_P$ (\phxsk), whereas   $\Dm $  consists of $n_P+1$  of such  subsets between the $\phxs$  and the top and bottom boundary of $\Domainspace$ which we denote by $\Dm _j,j=0,\ldots,n_P$. Then we can apply \eqref{quad_2D_0} to derive the approximation 
\begin{align}
	\label{eq:av:fluid}
	\Qf  \approx\frac{1}{|\Df|} \sum_{j=1}^{n_P} |\Df_j| \Qav^{\Df_j} 
	= \OutputF \,Y \quad\text{where}\quad  \OutputF= \frac{1}{|\Df|} \sum_{j=1}^{n_P} |\Df_j| \Cav^{\Df_j}.
\end{align}
An approximation of the form  $\Qm \approx \OutputM \,Y$ can be obtained analogously.  Further, from Eq.~\eqref{av_temp_storage} the approximation for the average temperature in the whole  storage can be derived as
\begin{align}
	\label{eq:av:storage}
	\Qmf \approx \OutputMF \,Y \quad \text{with}\quad  \OutputMF= \frac{|\Dm |}{|\Domainspace|} \, \OutputM + \frac{|\Df |}{|\Domainspace|}\,\OutputF .
\end{align}

In Appendix \ref{append:out:bottom} we derive  approximations $\Qout\approx \OutputOut \,Y$  and $\Qbottom\approx \OutputBottom \,Y$ for the average temperatures at the outlet and the bottom boundary, respectively. Here, the line integrals in the definitions  \eqref{Rout} and \eqref{RB} of these two characteristics are approximated by trapezoidal rule.

\section{Numerical Results}
\label{Numerical}
In this section we present results of numerical experiments based on the finite difference discretization \eqref{Matrix_explicit} of the heat equation \eqref{heat_eq2}.
We determine the spatio-temporal  temperature distribution in the storage. Further, we study the impact of the  $\phx$ topology  and vary the number and arrangement of the $\phxsk$. In Subsecs.~\ref{sub: num1}, \ref{sub: num2} and \ref{sub: num3} we present results for a storage with one, two and three \phxsk, respectively.
For these experiments we also compute and compare  certain aggregated characteristics which are introduced in Sec.~\ref{sec:Aggregate} and   computed via  post-processing of the temperature distribution. 

Note that we provide additional video material showing animations of the temporal evolution of the spatial temperature distribution for which in the following we can present snapshots only.
The videos are available at {\small\url{www.b-tu.de/owncloud/s/D68fmqXRcgbesKj}}~.

\subsection{Experimental Settings}	
\label{Num_Settings}
The  model and discretization parameters are given in Table \ref{tab:cap1}.
The storage is charged and discharged via $\phxs$  filled with a moving fluid and thermal energy is stored by raising the temperature of the storage medium.
We recall the open architecture of the storage which is only insulated at the top and the side but not at the bottom. This leads to an additional heat transfer to the underground for which we assume a constant temperature of  $\Qg(t)=15 \Celsius$.
In the simulations the fluid is assumed to be water while the storage medium is dry soil. During charging  a pump moves the fluid with constant velocity $\vconst$ arriving with constant temperature $\Qin(t)=\QinC=40\Celsius$  at the inlet.  If this temperature   is higher than in the vicinity of the \phx, then a heat flux into the storage medium is induced. During discharging the inlet  temperature is $\Qin(t)=\QinD=5\Celsius$ leading to a cooling of the storage. At the outlet we impose a vanishing diffusive heat flux, i.e.~during pumping there is only a convective heat flux. 
We also consider waiting periods where the pump is off. This helps to mitigate saturation effects in the vicinity of the $\phxs$    which reduce the injection and extraction efficiency. During that waiting periods the injected heat (cold) can propagate to other regions of the storage. Since pumps are off  we have only diffusive propagation of heat in the storage and the transfer over the bottom boundary. 

\begin{table}[h]
	\centering	
	\begin{tabular}[T]{|lc|rrr|} \hline
		Parameters&&&Values& Units\\
		\hline		
		\textbf{Geometry} & & &&\\
		width  &$l_x$ &&$10$~&$m$ \\
		height  &$l_y$ &&$1$&$m$ \\
		depth  &$l_z$ &&$10$&$m$ \\
		diameter of $\phx$ &$d_P$ && $0.02$&$ ~m$ \\
		number of $\phxs$  &$n_P$&& $1,2,3$&\\\hline 
		\textbf{Material} & && &\\
		\textit{medium (dry soil)} & && &\\
		\hspace*{1em} mass density &$\rhom $ && $2000$ &$~kg/m^3$\\
		\hspace*{1em}  specific heat capacity &$c_p^m$ && $800$& $~J/kg\, K$\\
		\hspace*{1em}  thermal conductivity  &$\kappam $ &&$1.59$ &~ $W/m\,  K$\\	
		\hspace*{1em}  thermal diffusivity ~~$\kappam (\rhom  c_p^m)^{-1}$& $\am $ && $9.9375 \times 10^{-7}$&$m^2/s$\\		
		\textit{fluid (water)} & && &\\
		\hspace*{1em} mass density &$\rhof $ &&$998$ ~&$kg/m^3$\\
		\hspace*{1em} specific heat capacity&$\cpf$&&$4182$ & $~J/kg\, K$\\
		\hspace*{1em}  thermal conductivity  &$\kappaf$ &&$0.60$ &~ $W/m\,  K$\\
		\hspace*{1em}  thermal diffusivity ~~$\kappaf(\rhof \cpf)^{-1}$ &$\af $&& $1.4376\times 10^{-7}$&$m^2/s$\\			
		velocity during pumping & $~\vconst$ && $ 10^{-2}$&$~m/s$\\
		heat transfer coeff.~ to underground & $\heattransfer$ && $10$&$ W/(m^2~ K)$\\	
		initial temperature 	&$Q_0$ &&  $10 $  and $35$ &$~\Celsius$\\
		inlet temperature: charging  & $\QinC$ && $40 $&$~\Celsius$\\
		\phantom{inlet temperature:} discharging  &	$\QinD$ && $5 $&$\Celsius$\\
		underground temperature &	$\Qg$ && $15$&$ ~\Celsius$\\
		\hline 
		\textbf{Discretization} &&& & \\
		step size&$h_x$&&$0.1$&$m$\\
		step size& $h_y$&&$0.01$&$~m$\\
		time step&  $\tau$ && $1$&$ ~s$\\		
		time  horizon &  $~T$ &&  $36$ and $72$&$ ~h$\\
		\hline
	\end{tabular}

	\medskip	
	\mycaption{ Model and discretization parameters.}	
	\label{tab:cap1}	
\end{table}

\subsection{Storage With  One Horizontal Straight \phx}
\label{sub: num1}

\begin{figure}[h]
		\centering
		\hspace*{-0.01\textwidth}
		\includegraphics[width=0.49\textwidth,height=0.2\textwidth]{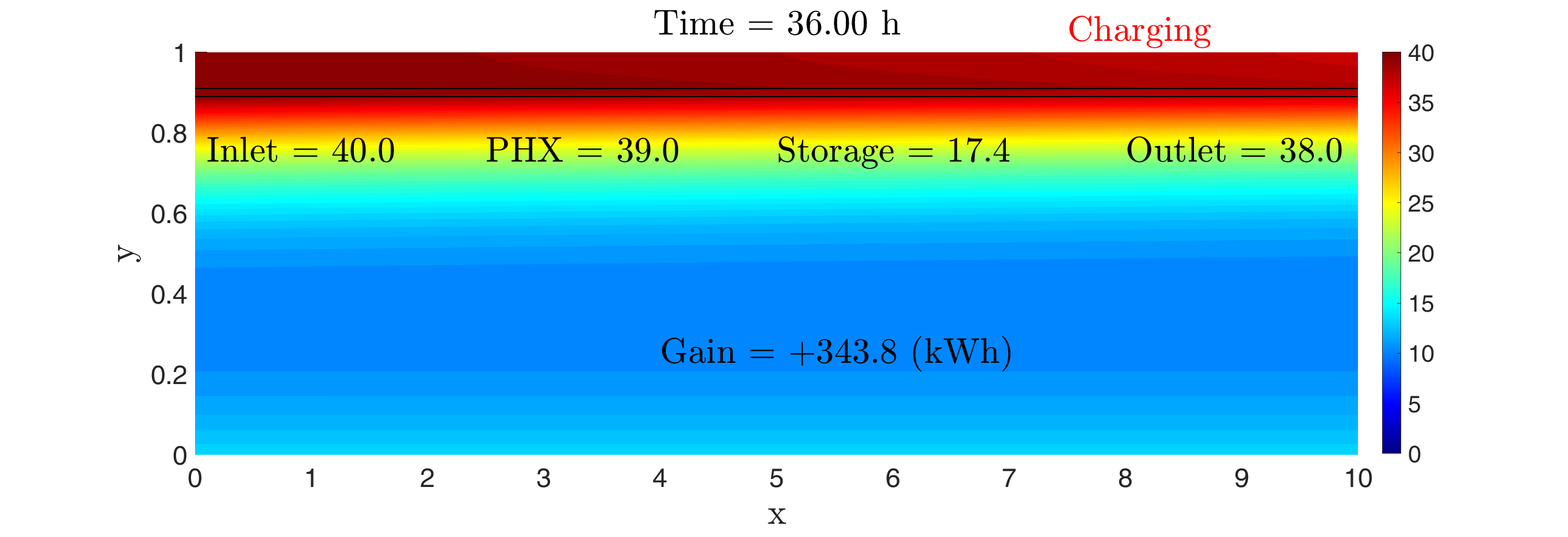}
		\hspace*{-0.001\textwidth}
		\includegraphics[width=0.49\textwidth,height=0.2\textwidth]{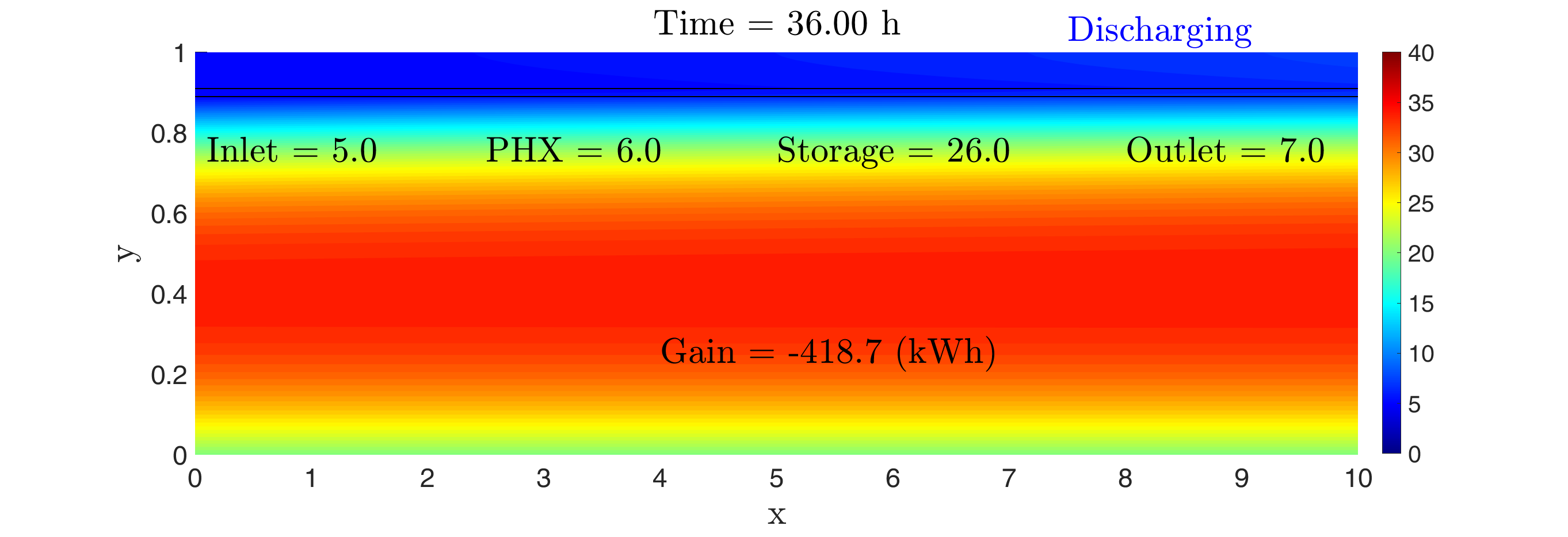} 
		
		\centering
		\hspace*{-0.01\textwidth}
		\includegraphics[width=0.49\textwidth,height=0.2\textwidth]{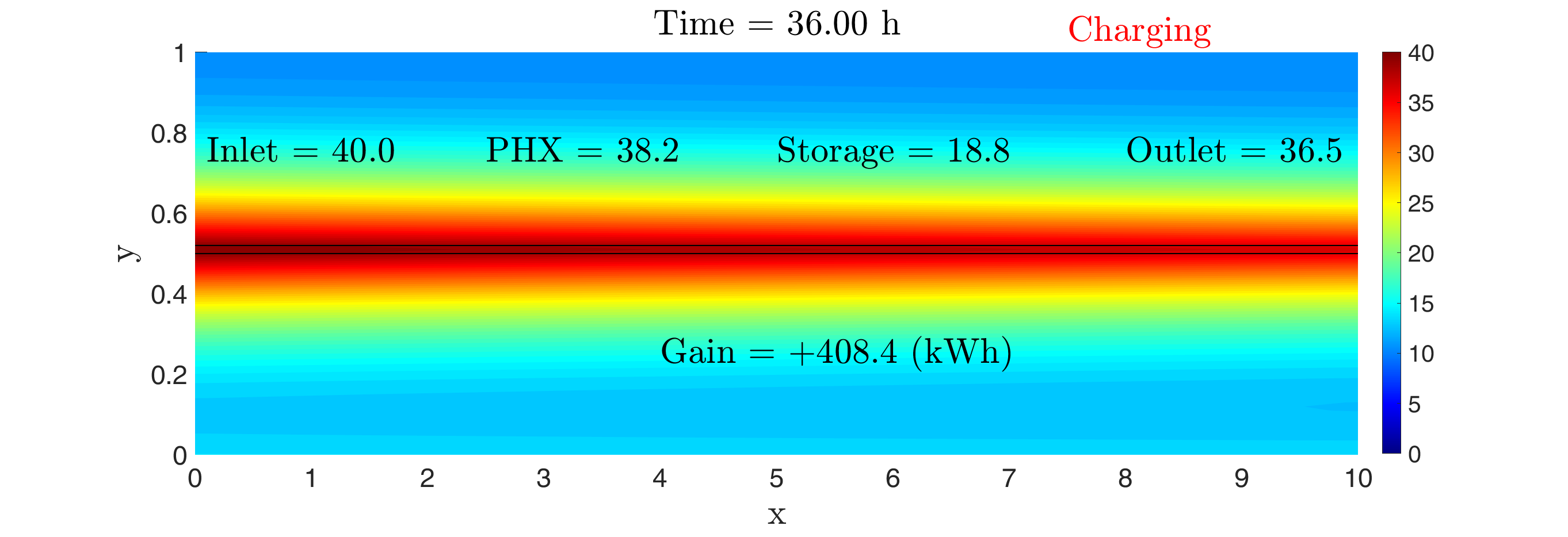}
		\hspace*{-0.001\textwidth}
		\includegraphics[width=0.49\textwidth,height=0.2\textwidth]{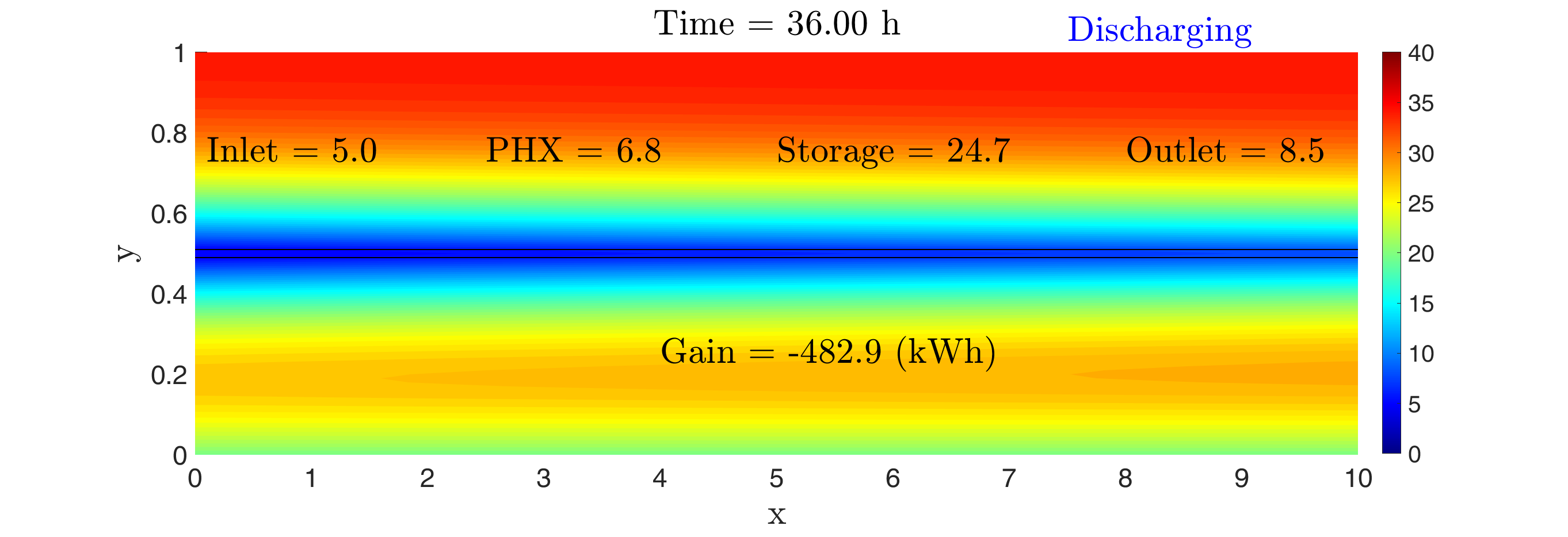}

		\centering
		\hspace*{-0.01\textwidth}
		\includegraphics[width=0.49\textwidth,height=0.2\textwidth]{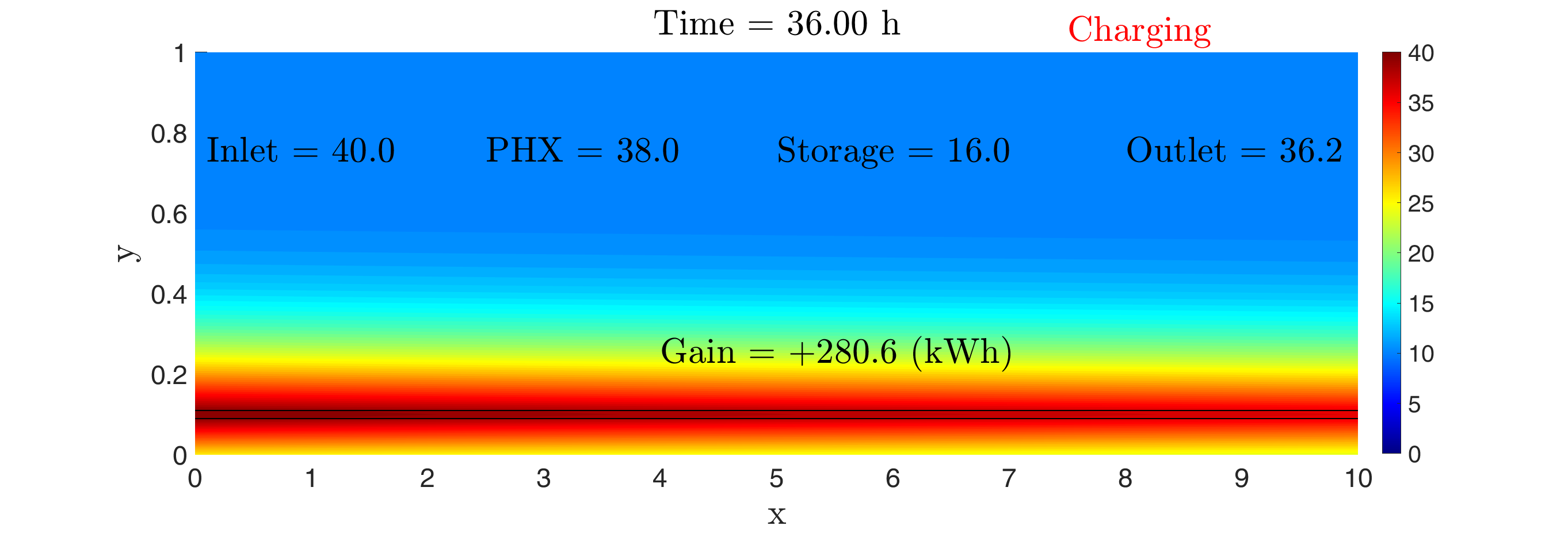}
		\hspace*{-0.001\textwidth}
		\includegraphics[width=0.49\textwidth,height=0.2\textwidth]{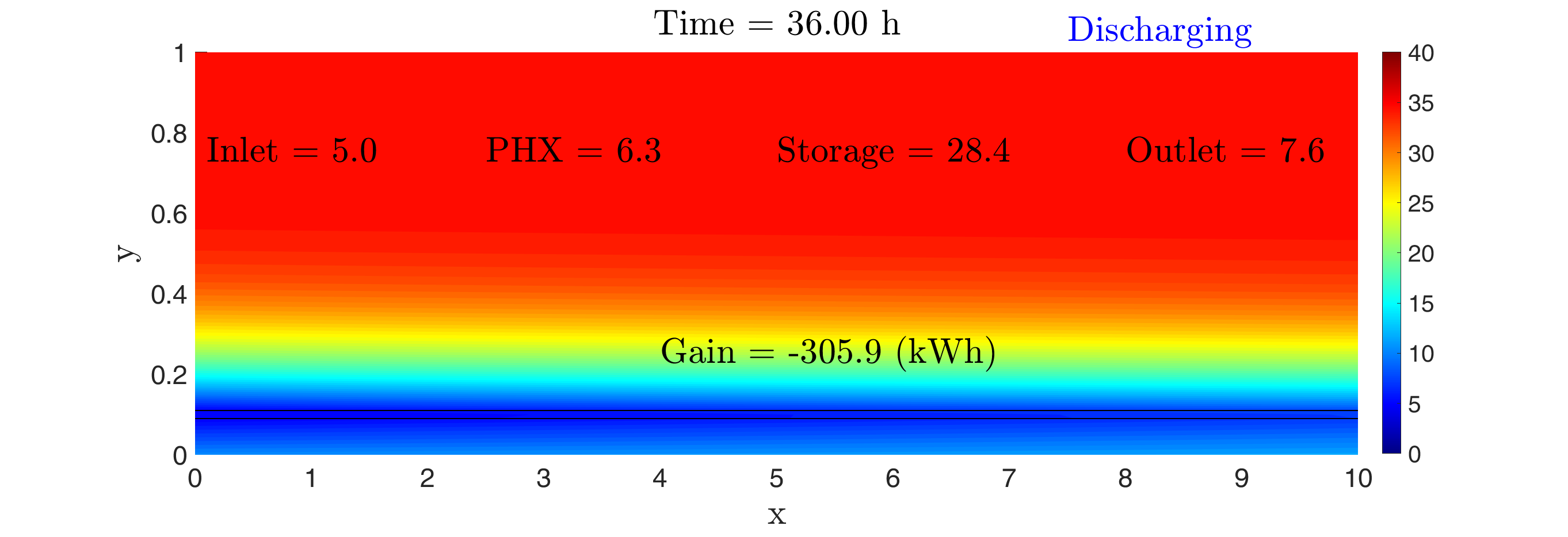} 
		
		\mycaption{Spatial distribution of the temperature in the storage with  one horizontal $\phx$ at vertical position $p$ after  of $36$ hours of charging (left) and discharging (right). \newline 
			Top: $p=90~cm$. Middle:   $ p=50~cm$. Bottom:  $ p=10~cm$.  }
		\label{fig:fig1c}
	\end{figure}  
	
	\begin{figure}[h]
		
		\centering
		\hspace*{-0.05\linewidth}
		\includegraphics[width=0.49\textwidth,height=0.35\textwidth]{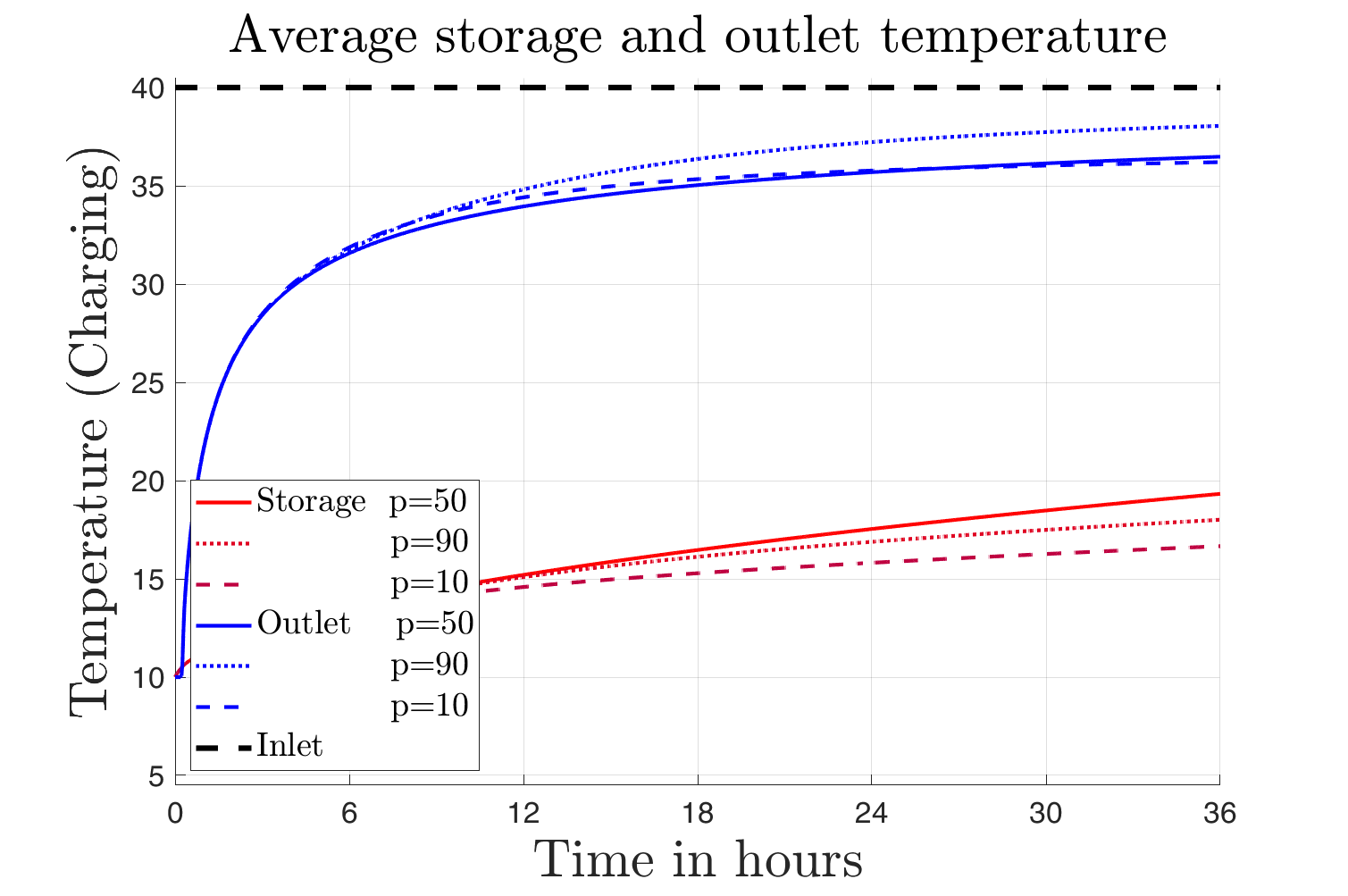}
		\includegraphics[width=0.49\textwidth,height=0.35\textwidth]{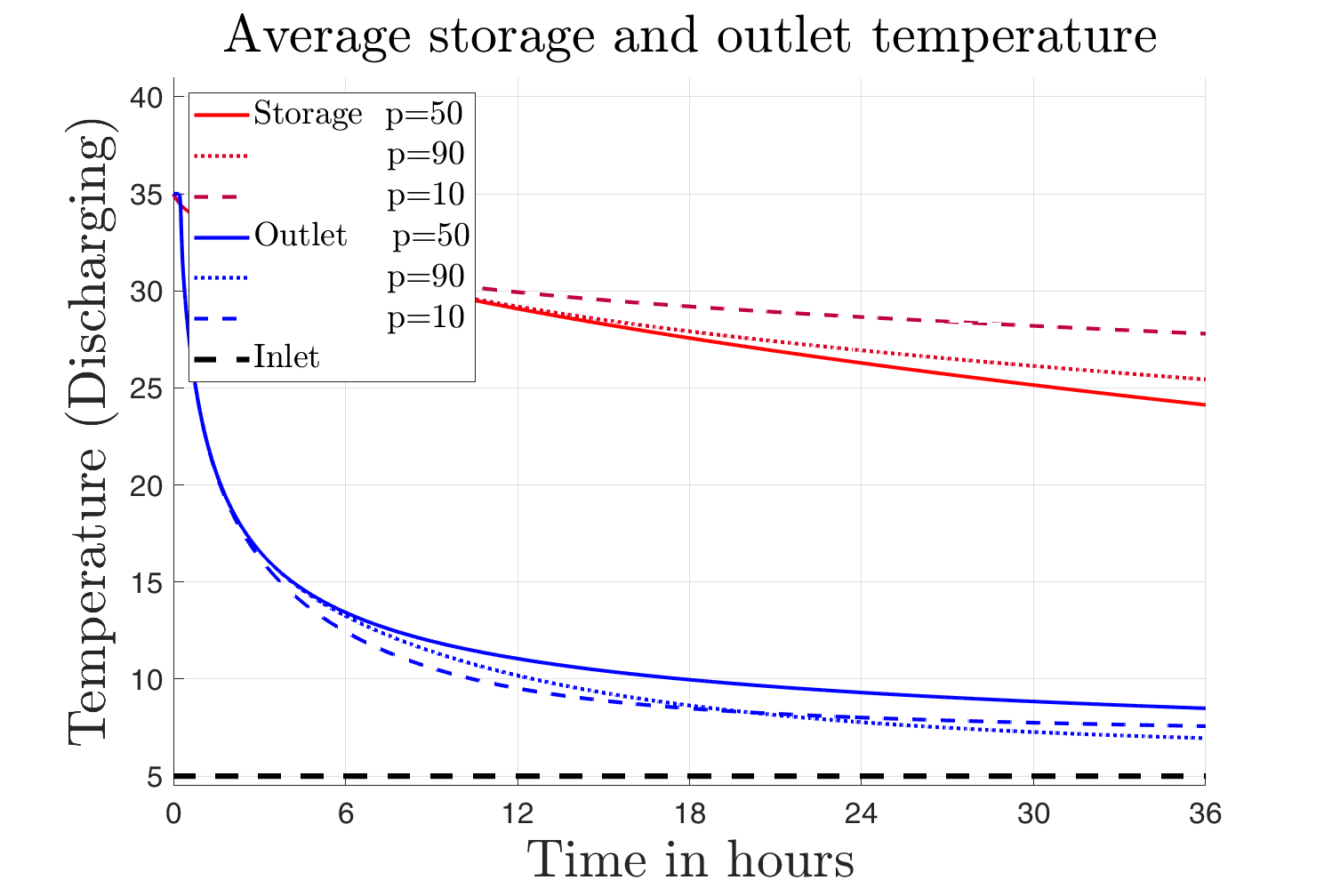}
		
		\mycaption{ Average temperature in the storage $\Qmf$  and average outlet temperature $\Qout$ after $36$ hours for  a storage with one horizontal $\phx$ at different vertical positions.~
			Left: Charging.~
			Right: Discharging.}
		\label{fig:fig1b}
	\end{figure}
	In this experiment we run simulations with one horizontal $\phx$ located at different vertical positions $p$ between the bottom ($p=0~cm$) and the top  ($p=l_y=100~cm$) of the storage. 
	We compare the spatial temperature distributions as well as aggregated characteristics such as the average temperature {in the storage} $\Qmf(t)$,   the average outlet temperature $\Qout(t)$, and the gain or loss of energy $\Gmf(0,T)$ in the storage during a period of $T=36$ hours. Charging is realized by sending fluid through the $\phx$ for $36$ hours. It arrives at the inlet  with constant temperature $\QinC(t)=40\,\Celsius$. We start with  an initial temperature  $Q(0,x,y)=10\,\Celsius$, uniformly distributed in the storage. In the experiment with discharging  we start with  an uniformly distributed  initial temperature $35\,\Celsius$. For $36$ hours the storage is cooled by the moving  fluid arriving at the storage inlet with  constant temperature $\QinD(t)=5\,\Celsius$.
	
	Fig.~\ref{fig:fig1c} shows the spatial distribution of the temperature in the storage after  $36$ hours of charging (left) and discharging (right) where we used three different vertical positions $p$ of the $\phxk$. In the top panels  the $\phx$ is located close to the insulated top boundary ($p=90~cm$). The panels in the middle show the results for a $\phx$ in the center ($p=50~cm$) while in the bottom panels the $\phx$ is close to the bottom boundary ($p=10~cm$). Recall that the bottom is open and allows for heat transfer to the underground with constant temperature $\Qg(t)=15~\Celsius$. 
	Fig.~\ref{fig:fig1b} plots the corresponding average temperatures in the storage and at the  outlet  against time.
	In Fig.~\ref{fig:fig1c} it can be seen that warming  and cooling mainly takes places in a vicinity of the $\phx$ and after $36$ hours the temperature in more distant storage domains is only slightly changed. Due to the direction of the moving fluid from left to right, warming and cooling in the left part  of the storage is  slightly stronger than in the right part. 
	A closer inspection of the results shows that  except in the experiment with the \phx  close to the bottom boundary ($p=10~cm$), after $36$ hours of charging  the temperatures in the  vicinity of that boundary are below the underground temperature $\Qg=15\,\Celsius$. Thus in addition to the injection of heat via the $\phx$  we also have an inflow of thermal energy from the warmer underground into the storage. This results in a  ``boundary layer''  which is slightly warmer than in the inner storage region.
	The reverse effect can be observed during  discharging where close to the bottom boundary the temperature is always above  $\Qg=15\,\Celsius$. This  induces a heat flux from the storage to the colder underground which contributes together with the extraction of heat via the $\phx$ to the total loss of thermal energy in the storage.

	\begin{figure}[h]
			\centering
			\hspace*{-0.05\linewidth}
			\includegraphics[width=0.49\textwidth,height=0.35\textwidth]{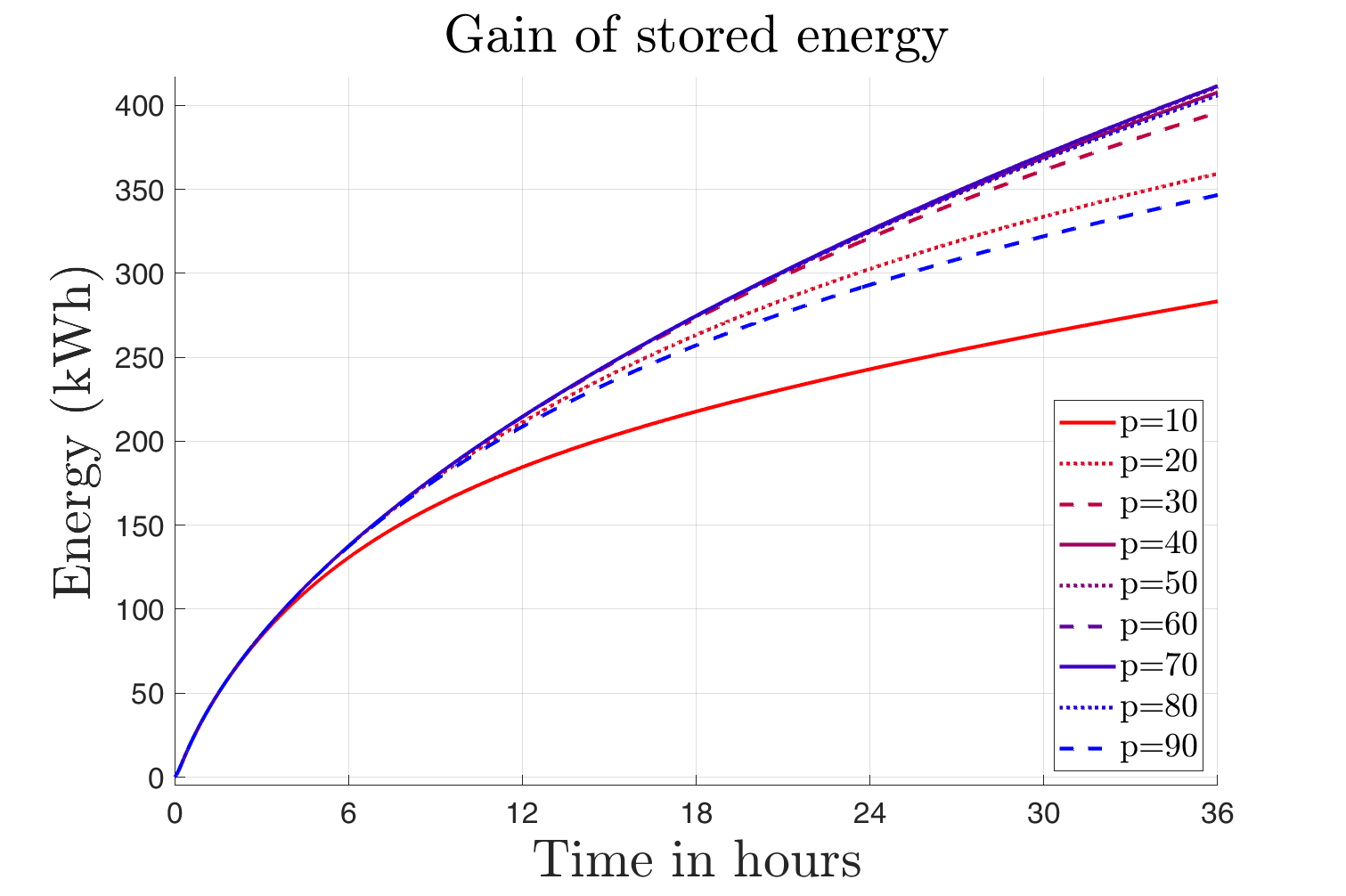}
			\includegraphics[width=0.49\textwidth,height=0.35\textwidth]{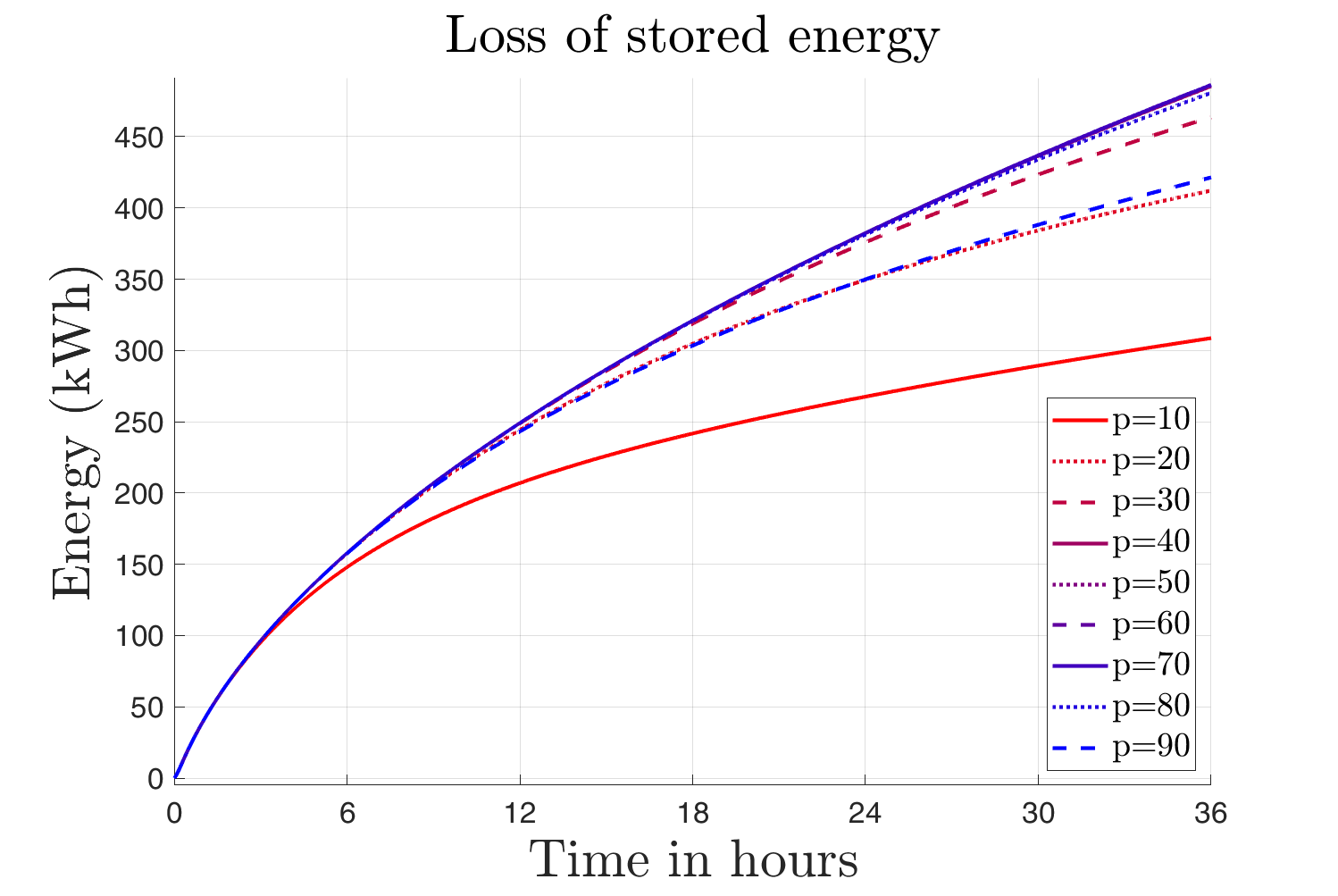}
			\\[1ex]
			\includegraphics[width=0.49\textwidth,height=0.35\textwidth]{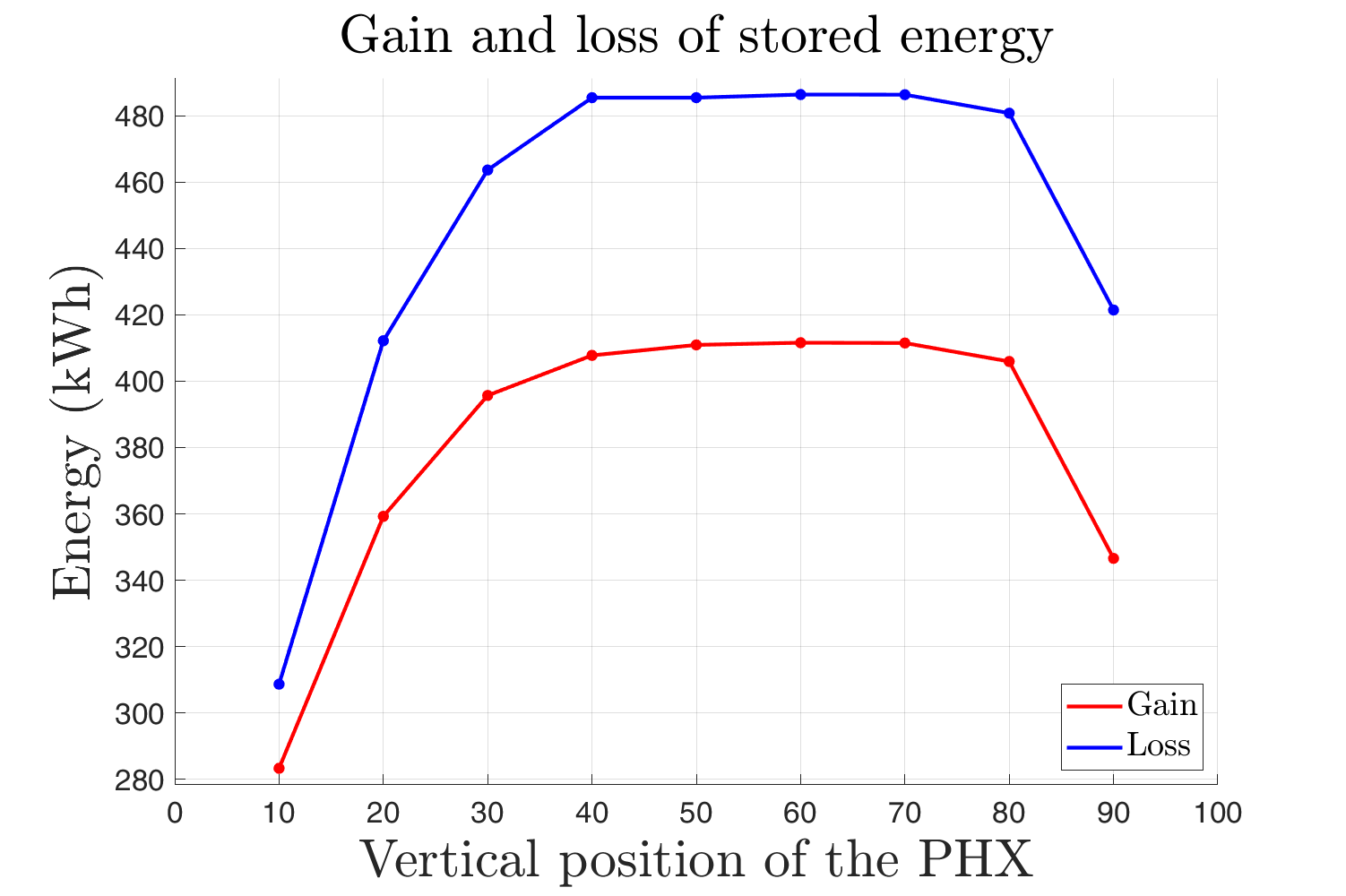}
		\mycaption{Gain and loss of stored energy  
			for a storage with one horizontal $\phx$ at different vertical positions. \newline
			Top left:  Gain  of stored energy  $\Gmf$ during charging. ~
			Top right: Loss of stored energy  $-\Gmf$ during discharging.\newline
			Bottom:  Gain  $\Gmf$ and loss  $-\Gmf$ of stored energy after $36$ hours of charging and discharging, respectively, depending on vertical $\phx$ position $p$.}
		\label{fig:fig1e}
	\end{figure}
	In Fig.~\ref{fig:fig1e} we plot in the upper panels the gain $\Gmf$ (respectively loss $-\Gmf$) of thermal energy during $36$ hours of charging (respectively discharging) against time for vertical positions $p=10,~20,\ldots,90~cm$. The lower panel shows these  quantities at the end of the 36 hour charging and discharging period, depending on the vertical $\phx$ position $p$. 
	In the first 4 hours of charging there are almost no visible deviations in the gains and losses,  but after $36$ hours  we can see a clear dependence  of the \phx's vertical position $p$. Further, for all $p$ we observe a decaying slope  of the curves in the upper plots. This can be explained by the ``thermal saturation'' in the vicinity of the $\phx$ and the slow diffusive propagation of the heat to the more distant regions of the storage. It shows that (dis)charging  the storage becomes less efficient after longer periods of operation. Injecting (extracting)  a certain amount of energy takes longer and needs more electricity consumed by the pumps. This effect suggests to interrupt (dis)charging and include waiting periods in which the heat (cold) in the vicinity of the $\phxs$  can propagate to other regions of the storage. The impact of such waiting periods will be studied in more detail in  Subsec.~\ref{sub: num2}.

	The results for $p=40,\ldots,70~cm$ are quite similar. However,  for  $\phx$ locations close to the open bottom boundary ($p=10,~20~cm$) and the insulated top boundary ($p=90~cm$) we observe remarkable deviations. Here charging and discharging is considerably slower and gains and losses of thermal energy are smaller. For a $\phx$ close to the top this can be explained by the saturation of the storage domain in the vicinity of the $\phxk$. During charging (discharging) the boundary and its insulation prevent  the propagation of heat into (from) the  inner storage regions. On the bottom boundary that effect is combined with  heat transfer to the underground. During charging a part of the injected heat is lost to the underground while during discharging the vicinity of the $\phx$ is also cooled by the colder underground. Thus as expected,  for an efficient operation of the storage the $\phx$ should be located in the central region of the storage.

	\subsection{Storage With  Two Horizontal Straight \phxs}
	\label{sub: num2}
	In this experiment we run the simulations with two horizontal $\phxs$  located symmetrically to the vertical mid level of $p=50 ~cm$ and  separated by a distance $d$ varying between $10~ cm$ and  $90~ cm$. Recall that placing a single $\phx$ at  $p=50 ~cm$ showed quite  good performance in the last subsection. First we study the spatial temperature distribution and some aggregated characteristics during (dis)charging for $T=36$ hours. Then we introduce waiting periods allowing the injected heat (cold) to spread within the storage. 
	
	\begin{figure}[h]
		\centering
		\hspace*{-0.01\textwidth}
		\includegraphics[width=0.49\textwidth,height=0.2\textwidth]{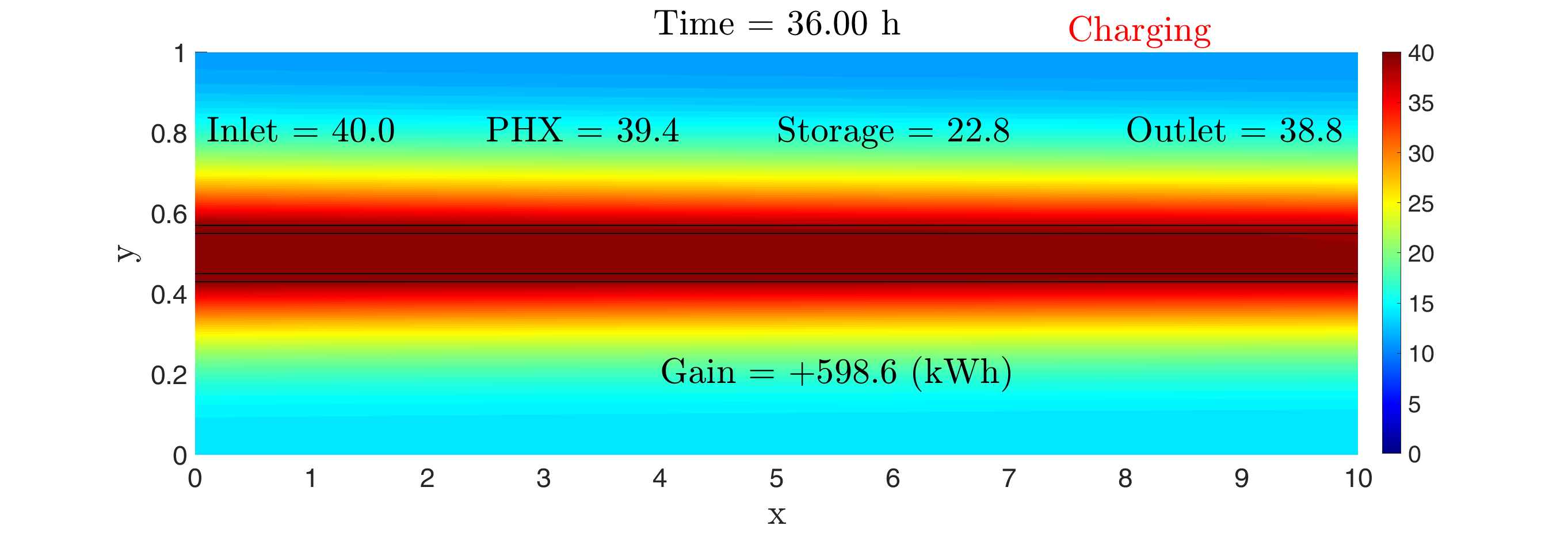}
		\hspace*{-0.0001\textwidth}
		\includegraphics[width=0.49\textwidth,height=0.2\textwidth]{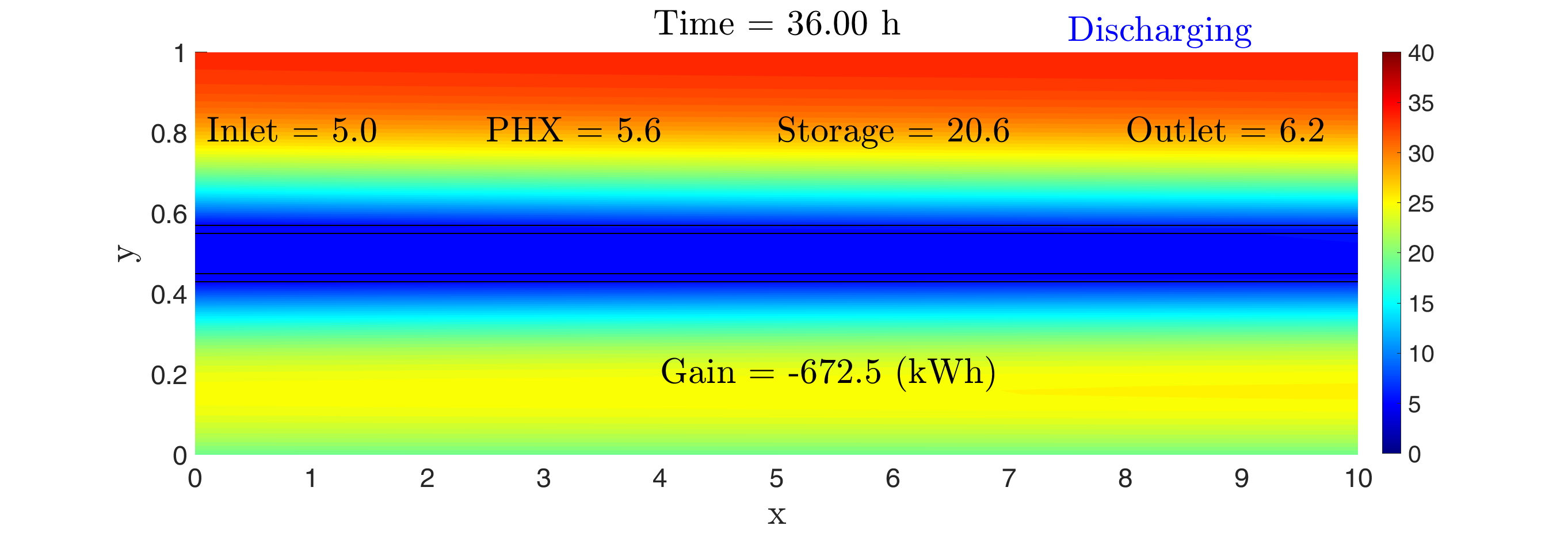} 
		
		\centering
		\hspace*{-0.01\textwidth}
		\includegraphics[width=0.49\textwidth,height=0.2\textwidth]{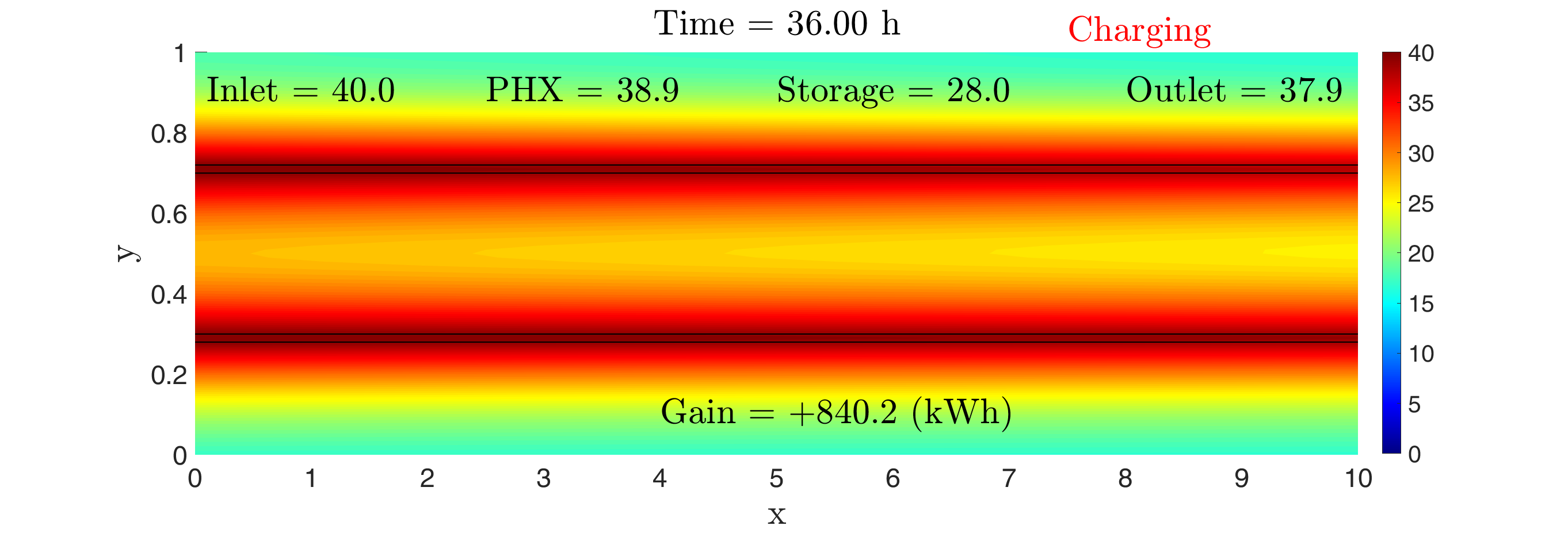}
		\hspace*{-0.0001\textwidth}
		\includegraphics[width=0.49\textwidth,height=0.2\textwidth]{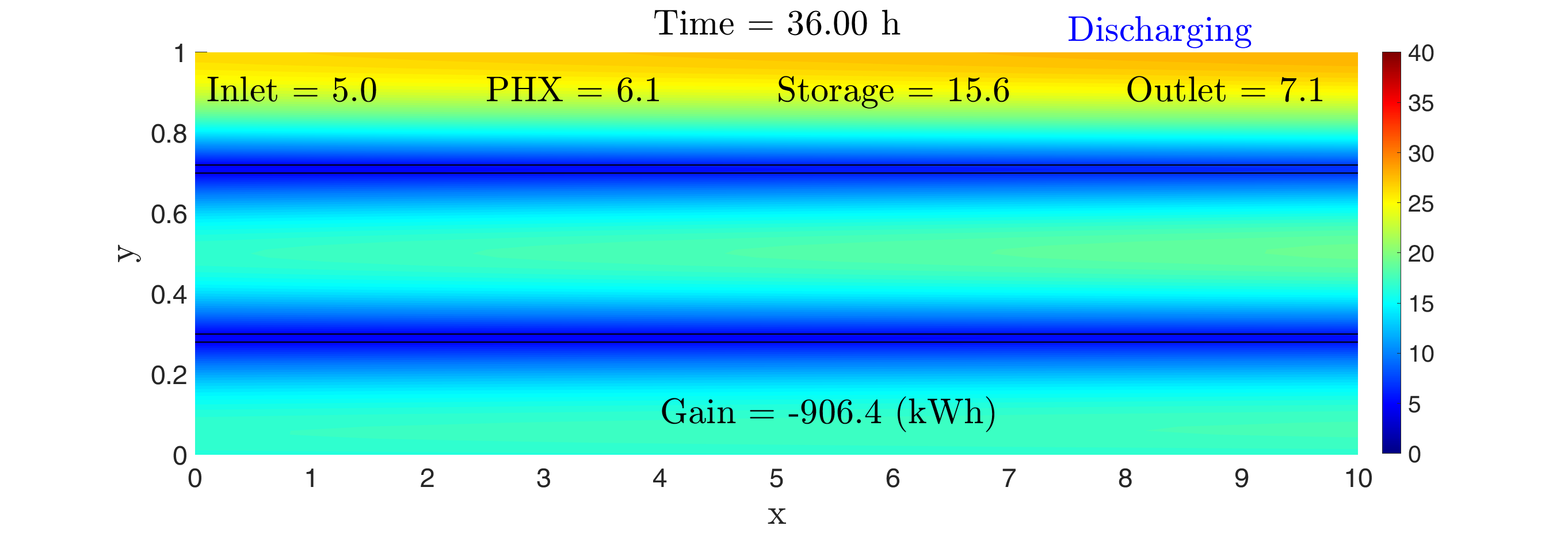}

		\centering
		\hspace*{-0.01\textwidth}
		\includegraphics[width=0.49\textwidth,height=0.2\textwidth]{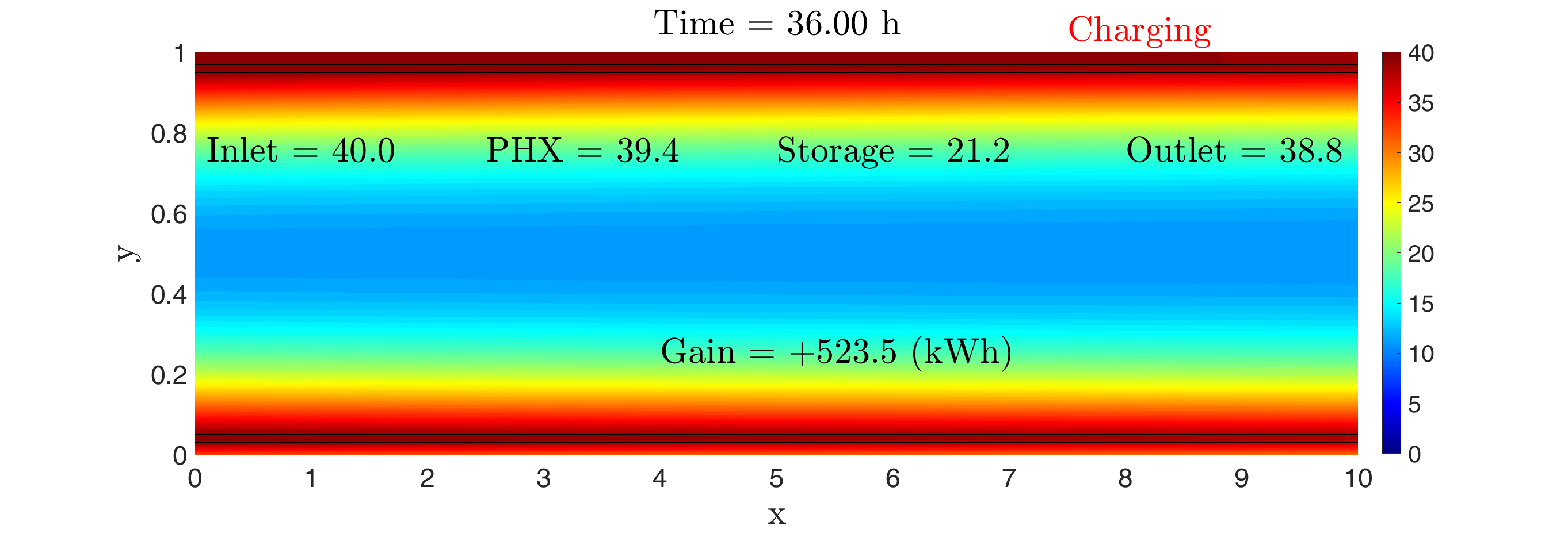}
		\hspace*{-0.0001\textwidth}
		\includegraphics[width=0.49\textwidth,height=0.2\textwidth]{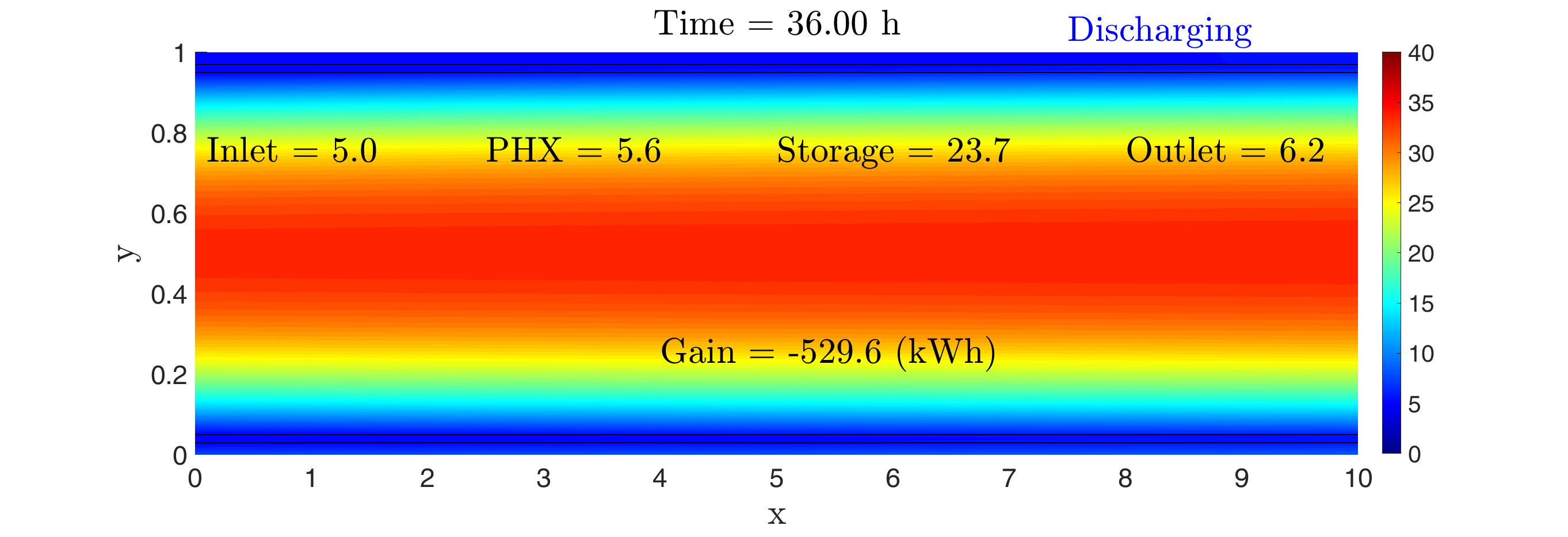} 
		\mycaption{Spatial distribution of the temperature in the storage with two horizontal $\phxs$  of vertical distance $d$ after  $36$ hours of charging (left) and discharging (right). \newline 
			Top: $d= 10~cm$. Middle: $d= 40~cm$. Bottom: $d= 90~cm$. }
		\label{fig:fig2c}
	\end{figure}  
	
	\begin{figure}[h]
			\centering
			\hspace*{-0.05\linewidth}
			\includegraphics[width=0.49\textwidth,height=0.35\textwidth]{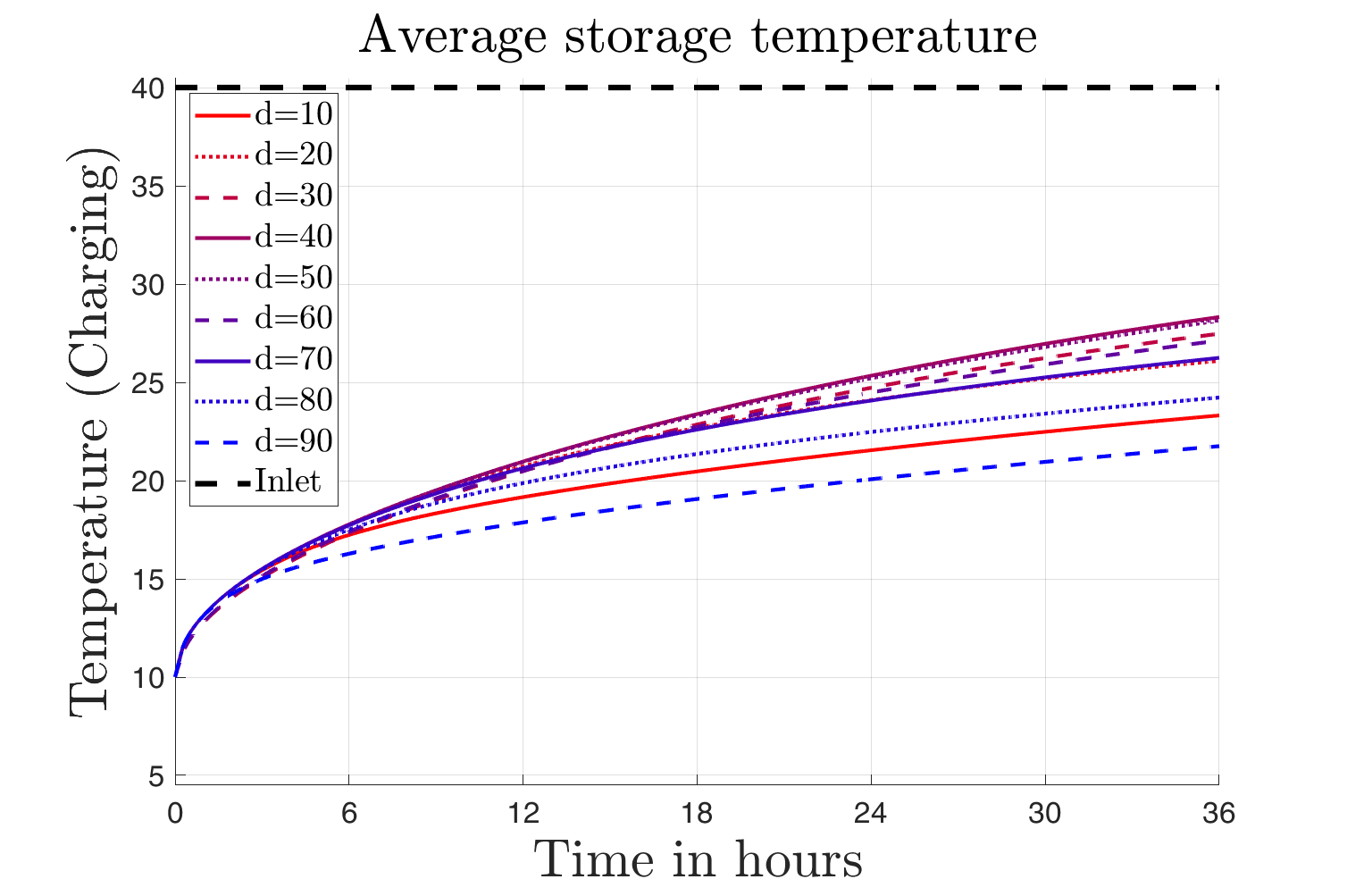}
			\includegraphics[width=0.49\textwidth,height=0.35\textwidth]{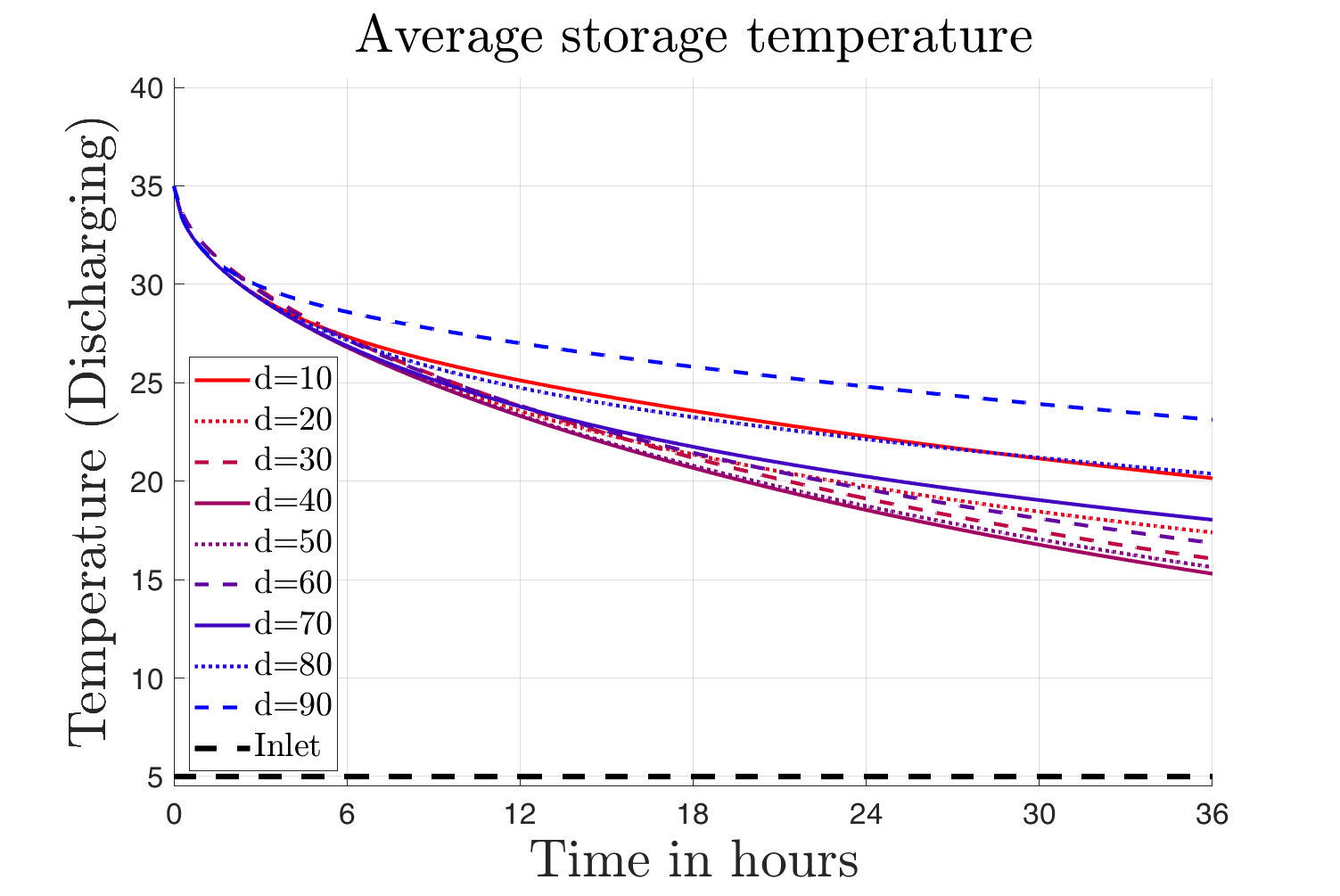}
		\mycaption{Average temperature in the storage $\Qmf$ during $36$ hours for  a storage with two horizontal $\phxs$  of different vertical distances.~
			Left: Charging.~
			Right: Discharging.}
		\label{fig:fig2b}
	\end{figure}

	\begin{figure}[h]
		\centering
		\hspace*{-0.05\linewidth}
		\includegraphics[width=0.49\textwidth,height=0.35\textwidth]{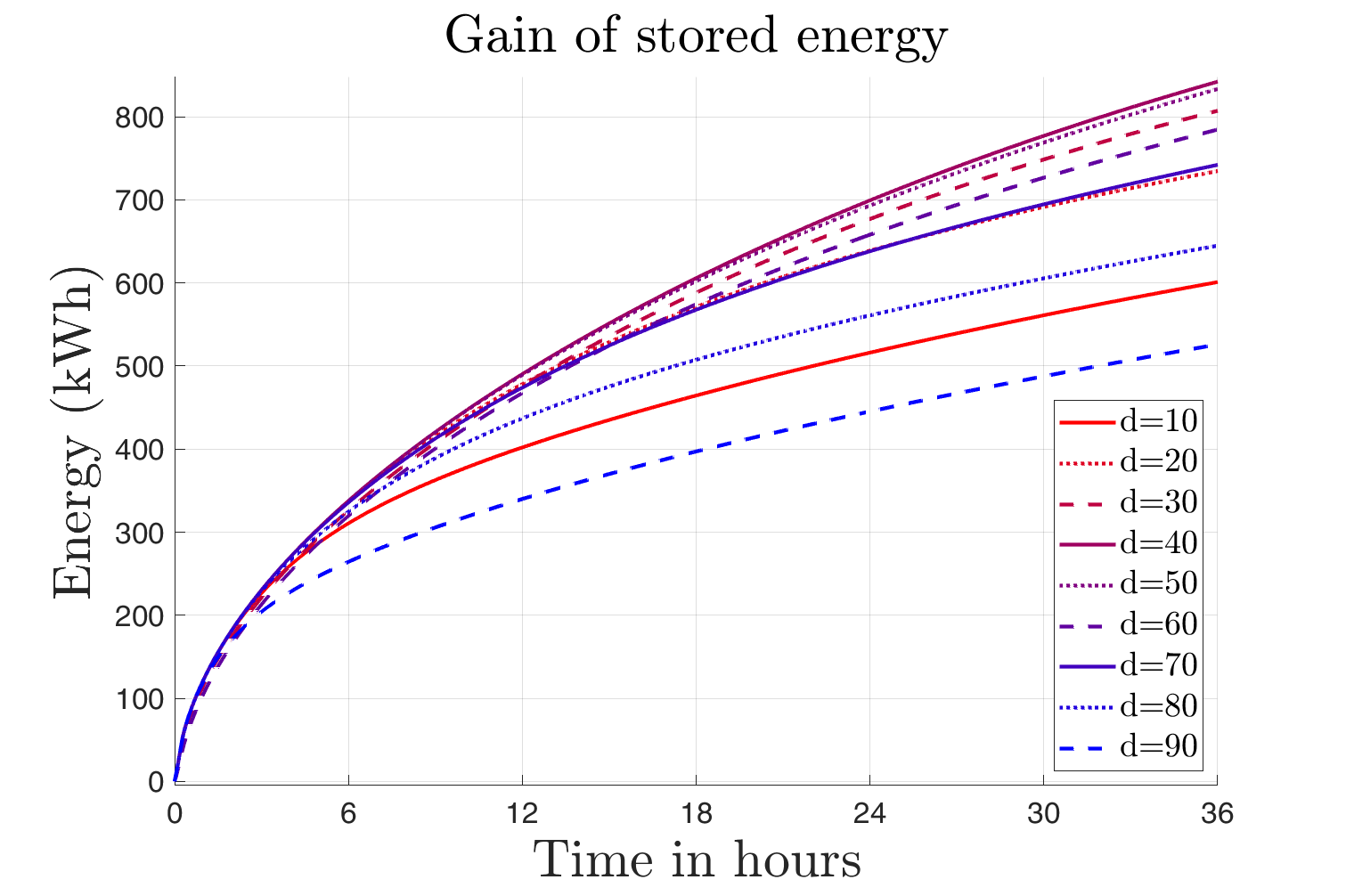}
		\includegraphics[width=0.49\textwidth,height=0.35\textwidth]{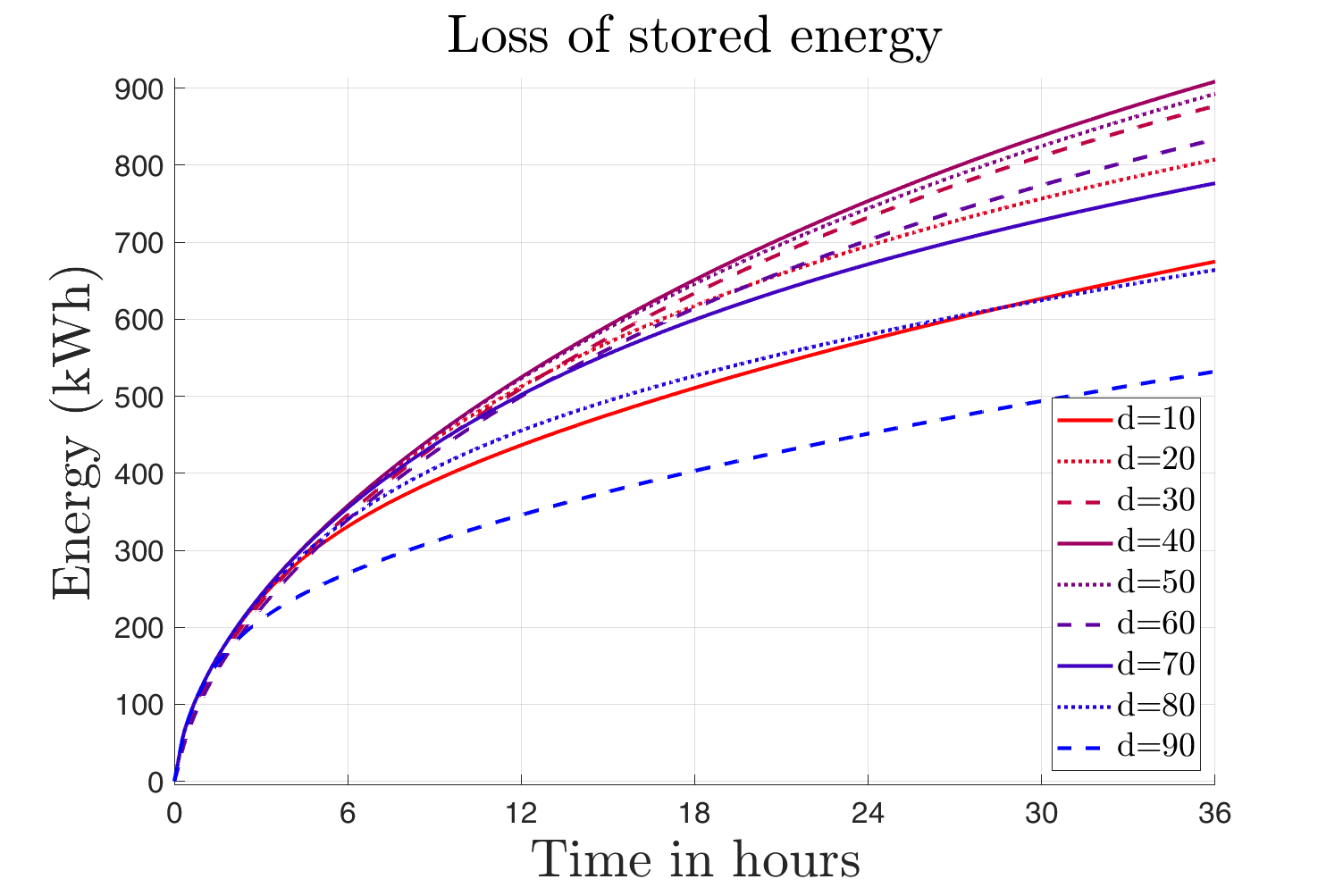}
		\\[1ex]
		\includegraphics[width=0.49\textwidth,height=0.35\textwidth]{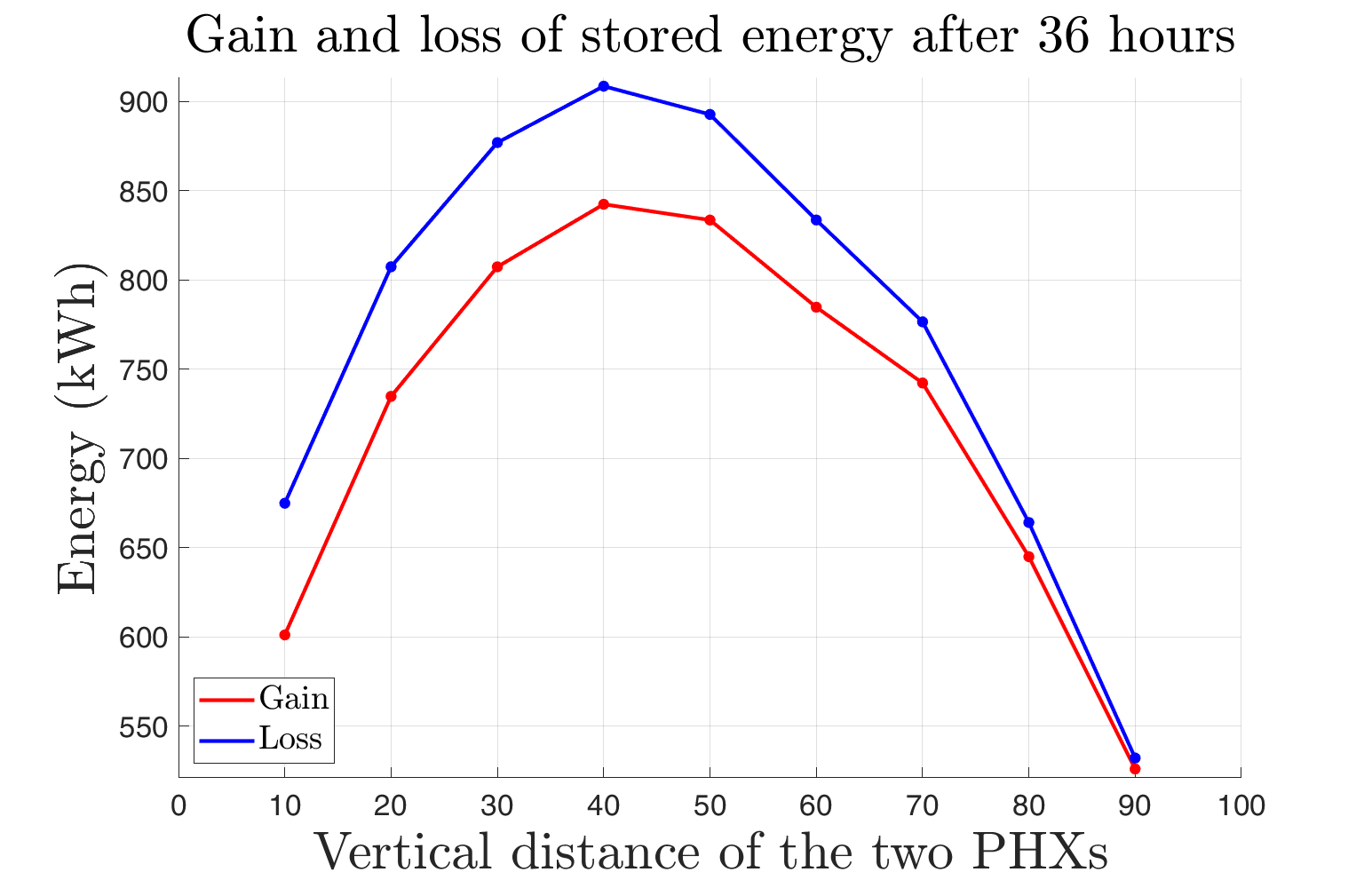}
		\mycaption{Gain and loss of stored energy   
			for a storage with two horizontal $\phxs$  of different distance $d$.  to $90~cm$.
			\newline
			Top left: Gain $\Gmf$ of stored energy during charging. ~
			Top right: Loss $-\Gmf$ of stored energy during discharging.\newline
			Bottom:   Gain  $\Gmf$ and loss  $-\Gmf$ of stored energy after $36$ hours  of charging and discharging, respectively, depending on distance $d$.}
		\label{fig:fig2g}
	\end{figure}
	
	\subsubsection{Charging and Discharging Without Waiting Periods}
	Fig.~\ref{fig:fig2c} shows for three different distances  $d$ of the two  $\phxs$  the spatial distribution of the temperature in the storage after  $36~h$ of charging (left) and discharging (right). In the top panels  the $\phxs$  are very close ($d=10~cm$). The panels in the middle show the results for two $\phxs$  at a distance $d=40~cm$ and in the bottom panels one $\phx$ is located close to the top and the other close to the  bottom boundary  ($d=90~cm$). 
	As in the experiment with only one $\phx$ it can be seen that  warming and cooling in the left part  of the storage is  slightly stronger than in the right part. It  mainly takes places in a vicinity of the $\phx$ whereas after $36~h$  temperatures in more distant regions are only slightly changed. 
	Thus, the spatial temperature distributions differ considerably for the three arrangements of two $\phxsk$. For a small distance ($d=10~cm$), we observe a strong saturation at a level close to the inlet temperature in the small region between the $\phxs$  while the region at the top is almost at the initial temperature and the region at the bottom is only slightly warmed (cooled) by the underground. For the $\phxs$  at distance $d=90~cm$, we observe an extreme saturation in the small layer between the upper $\phx$ and the top boundary while the lower $\phx$ is also warming (cooling) the underground. 
	
	Next we will have a look at aggregated characteristics. In Fig.~\ref{fig:fig2b} the average temperatures in the storage $\Qmf$ are plotted against time for distances of the $\phxs$  $d=10,20,\ldots,90~cm$. Fig.~\ref{fig:fig2g} presents the gain $\Gmf$ and loss $-\Gmf$ of thermal energy in the storage at the end of the charging and discharging  period, respectively. The figures reveal that  apart from the first  $4$ hours there is a strong impact of the $\phx$ distance. The most efficient mode of operation is obtained for  the $\phxs$  distance of $d=40~cm$. Here,  the gain (loss) of thermal energy during charging (discharging) is at maximum. These quantities strongly decay for smaller and larger distances because of the saturation effect
	which becomes stronger if $\phxs$  are arranged closer to each other or closer to the top and bottom boundary  of the storage.

	\subsubsection{Charging and Discharging With Waiting Periods}
	\begin{figure}[h]
			\centering
			\hspace*{-0.05\linewidth}
			\includegraphics[width=0.49\textwidth,height=0.35\textwidth]{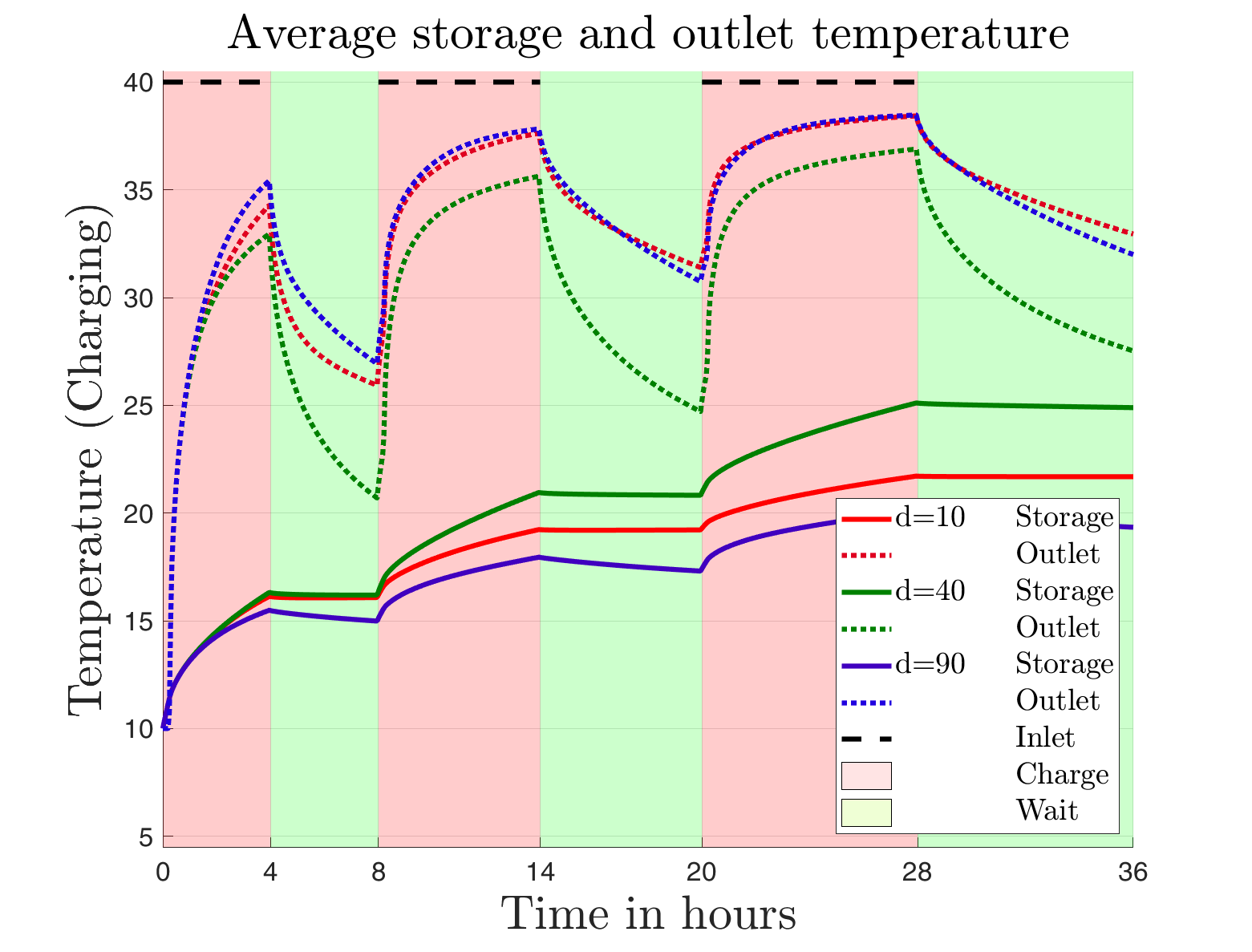}
			\includegraphics[width=0.49\textwidth,height=0.35\textwidth]{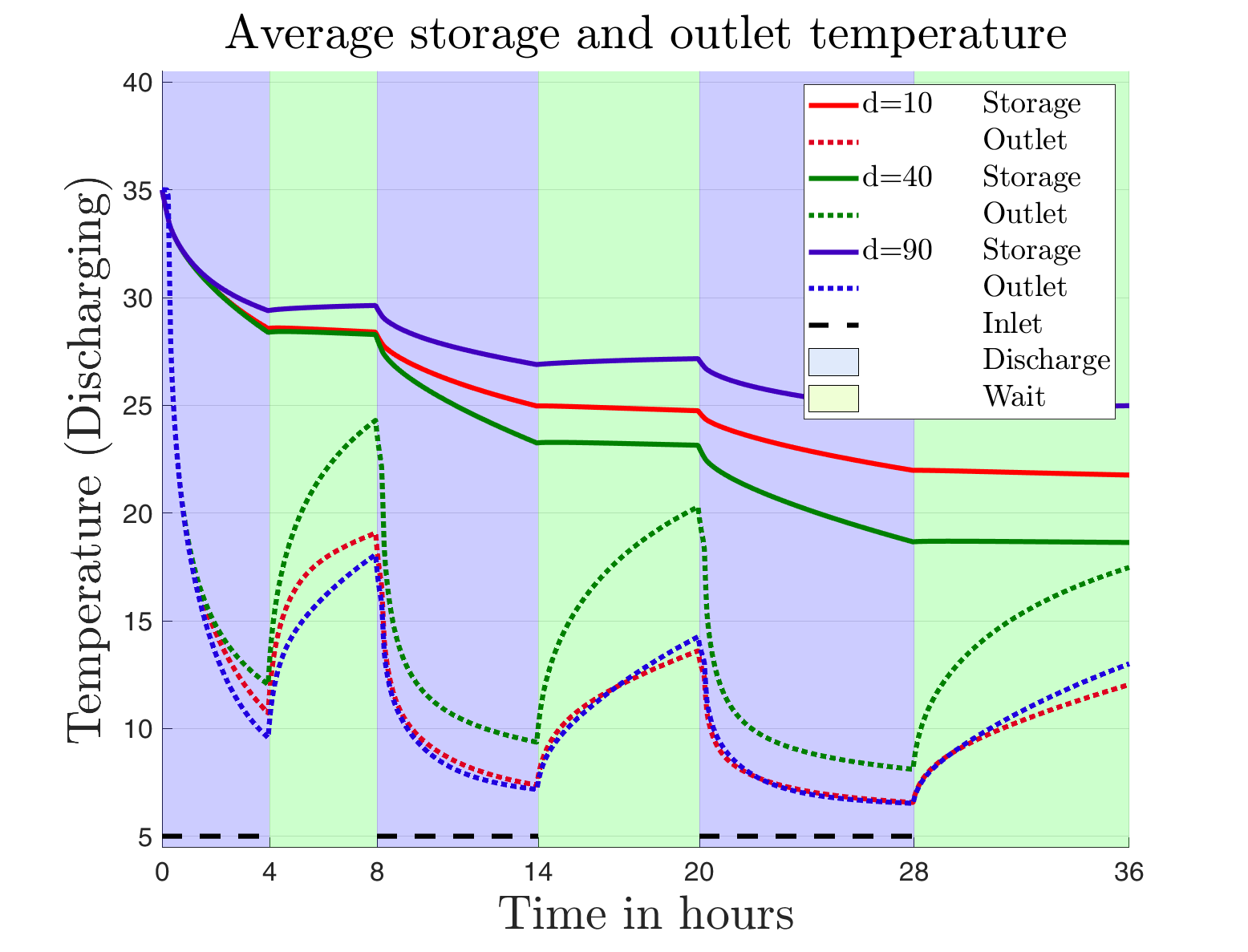}	
			\\[1ex]	
			\hspace*{-0.05\linewidth}
			\includegraphics[width=0.49\textwidth,height=0.35\textwidth]{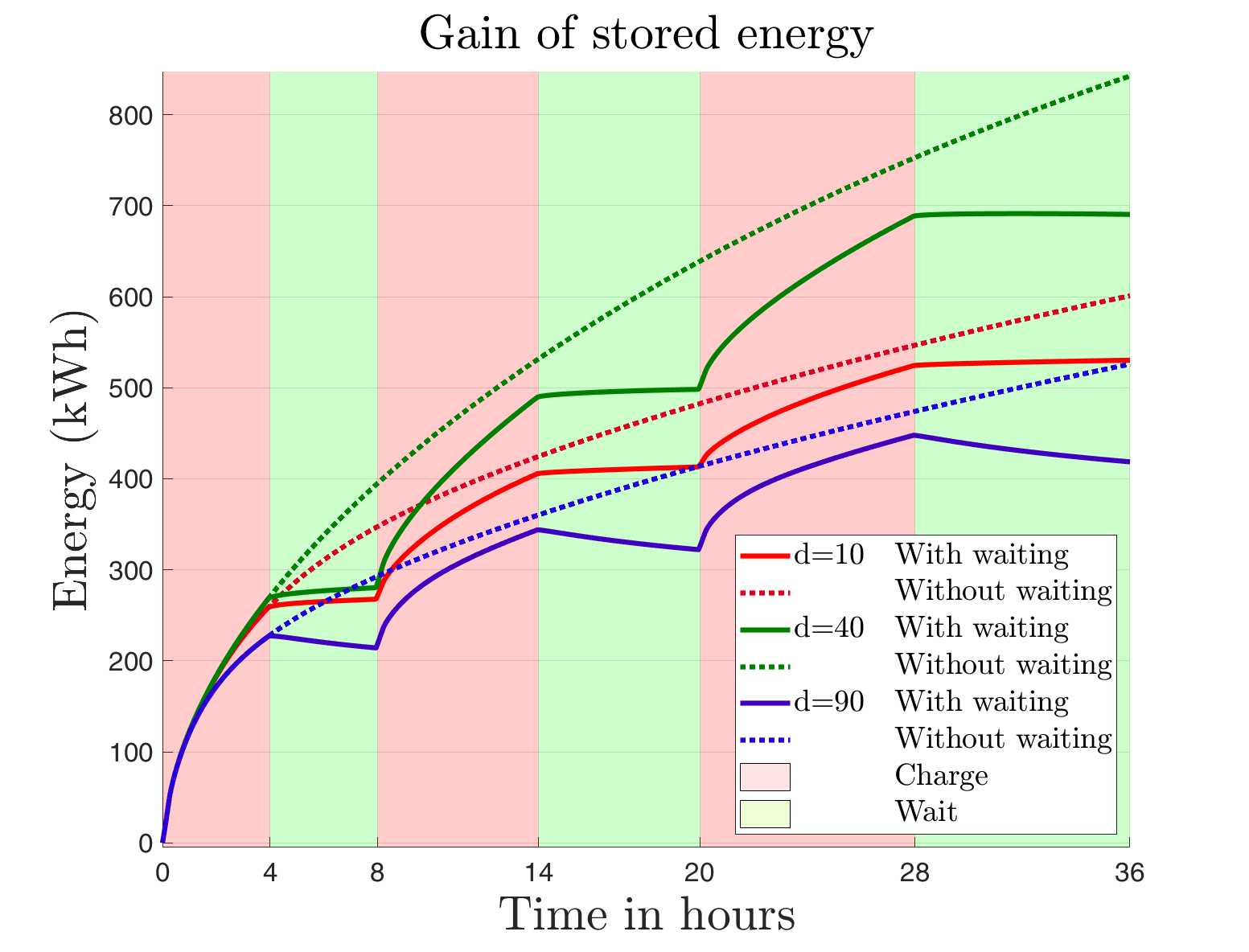}
			\includegraphics[width=0.49\textwidth,height=0.35\textwidth]{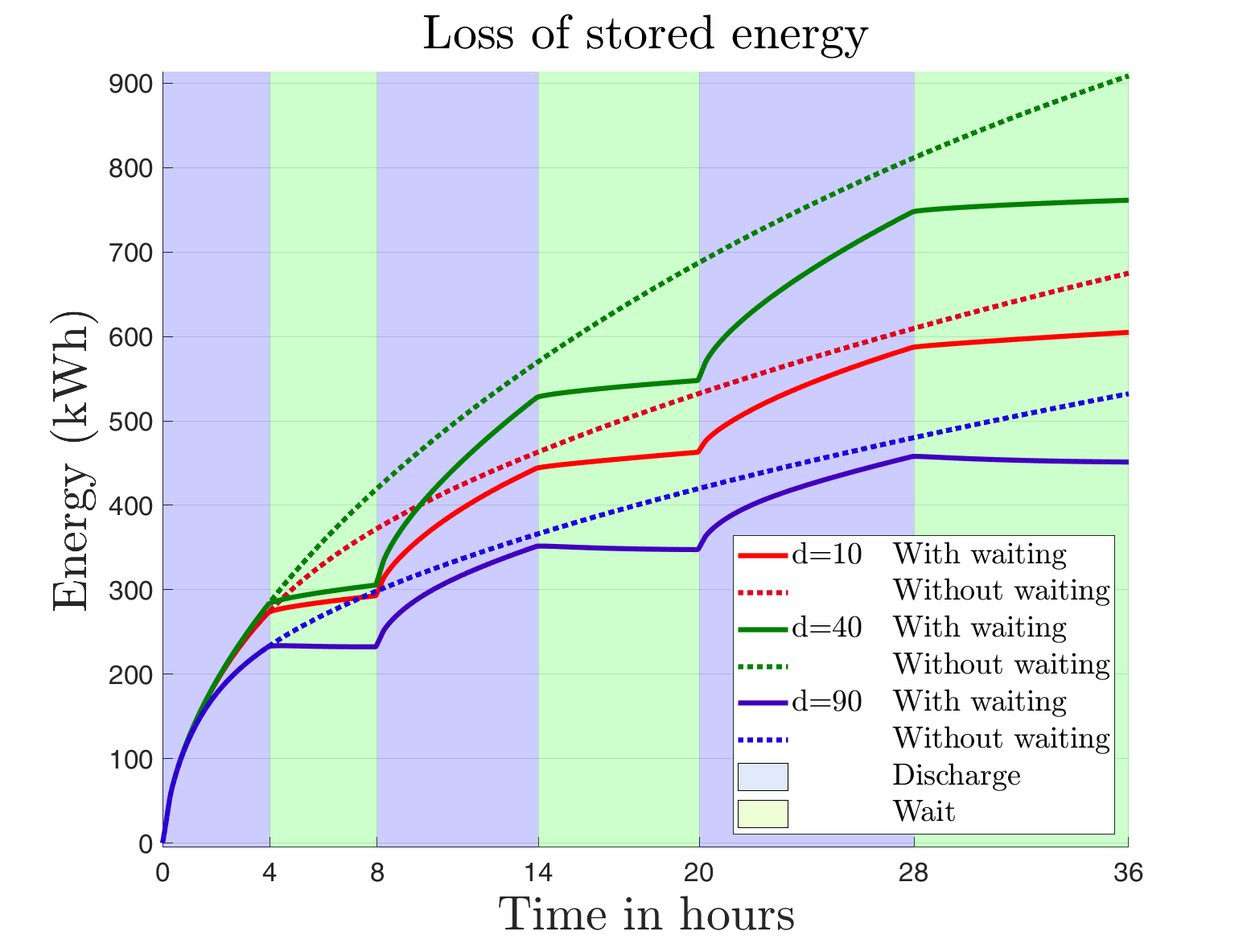}
		\mycaption{Charging and discharging during $36$ hours with several waiting periods for a storage with two horizontal $\phxs$  at distance $d=10~cm$, $d=40$, and $d=90~cm$.\newline
			Top: Aggregated characteristics $\Qmf$ and $\Qout$. ~Bottom: Gain $\Gmf$ / loss $-\Gmf$ of stored energy.\newline
			Left:  Charging. ~Right:  Discharging. 
		}	
		\label{fig:fig2e}
	\end{figure}
	The above experiments have shown how saturation effects can be mitigated by choosing an appropriate vertical distance of the two $\phxsk$. This option is only available in the design of the storage architecture and not during the operation of an already existing storage. Therefore,  we now want to examine another option, which is the  interruption of (dis)charging cycles allowing the heat injected to (extracted from) the vicinity of the $\phxs$  to propagate to the other storage regions. The idea is that after a sufficiently long waiting period the saturation in the vicinity of the $\phx$ is considerably reduced such that (dis)charging can resumed with higher efficiency. Although, the introduction of such waiting period will increase the time needed to inject (extract) a given amount of thermal energy it  reduces the saturation effect  and helps to save operational costs for  electricity used for running the pumps.
	
	In our experiments we divide the time interval $[0,T]$   into three subintervals of length $8,12,16$ hours. In each subinterval  (dis)charging is followed by a waiting period of the same length as it can be seen in Fig.~\ref{fig:fig2e} where charging, waiting and discharging periods are represented by red, green and blue background color. The top panels show the average temperatures  in the storage $\Qmf$ and at the outlet $\Qout$, respectively,  during charging and discharging. We compare  a storage with two $\phxs$  of distance $d=40~cm$ and a storage with more close-by  $\phxs$  $d=10~cm$ and two $\phxs$  at distance $d=90~cm$. Recall that  in the previous subsection we have seen that $d=40~cm$ allows for much more efficient operation than for $d=10, ~90~cm$. As expected, during the waiting periods the average temperatures at the outlet and in the $\phx$ decay after charging and rise after discharging. This is  due to the diffusion of heat in the storage, in particular the  heat flux induced by the different temperatures inside and outside the $\phxk$. During waiting the average temperature in the storage $\Qmf$ is almost constant since injection or extraction of heat is stopped. However, the heat transfer to and from  the underground at the bottom boundary continues also during waiting but the waiting periods are too short to produce  a visible change of $\Qmf$. In the two lower panels of Fig.~\ref{fig:fig2e} we compare the storage operation with and without waiting periods. We plot  the gain $\Gmf$ (loss $-\Gmf$) of thermal energy in the storage during charging (discharging) over time. 
	Note that for operation with waiting (dis)charging takes place only   50\% of the time. However,  for the ``optimal'' $\phx$ distance $d=40~cm$ the resulting gain (loss) reaches more than 80\% of the values for uninterrupted operation. For the less efficient cases of $\phxs$  at distance  $d=10~cm$ and  $\phxs$  at distance $d=90~cm$  that cause strong saturation effects  the differences are smaller and the gaps are quickly reduced to almost zero after resuming (dis)charging.

	\subsection{Storage With  Three Horizontal Straight \phxs}
	\label{sub: num3}
	\begin{figure}[h]
		
		\centering
		\hspace*{-0.01\linewidth}
		\includegraphics[width=0.49\textwidth,height=0.2\textwidth]{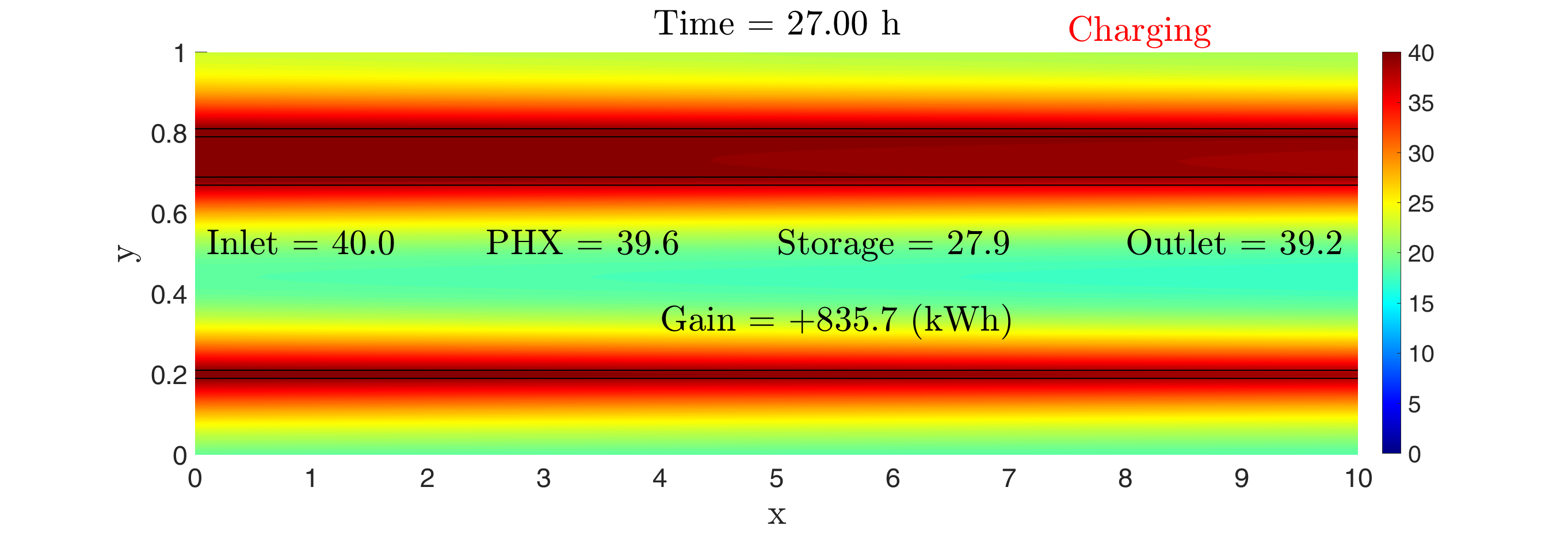}\hspace*{0.05\textwidth}
		\hspace*{-0.04\linewidth}
		\includegraphics[width=0.49\textwidth,height=0.2\textwidth]{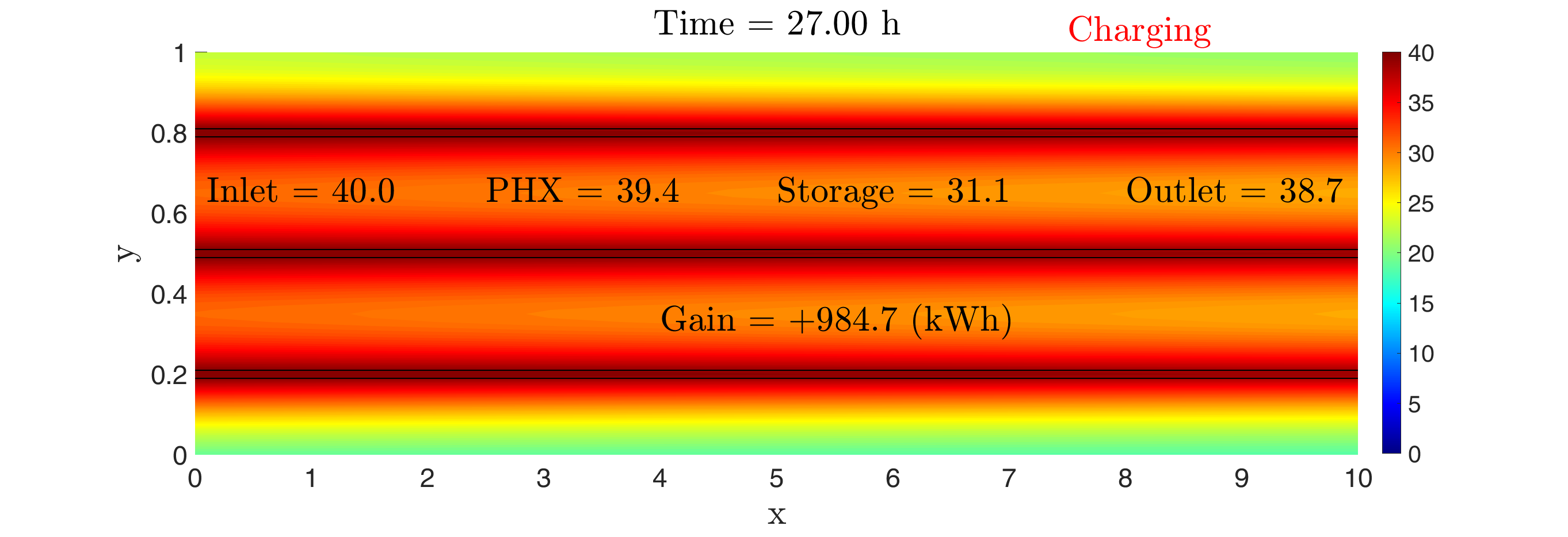} 
		
		\hspace*{-0.01\linewidth}
		\includegraphics[width=0.49\textwidth,height=0.2\textwidth]{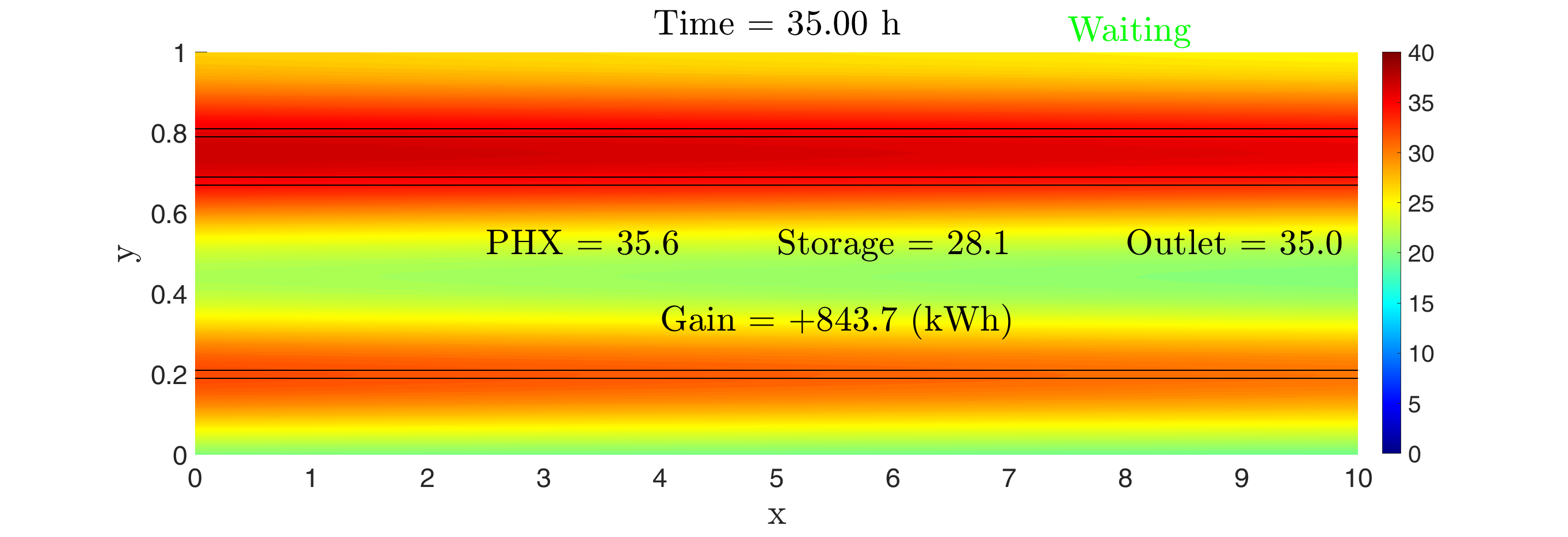}\hspace*{0.05\textwidth}
		\hspace*{-0.04\linewidth}
		\includegraphics[width=0.49\textwidth,height=0.2\textwidth]{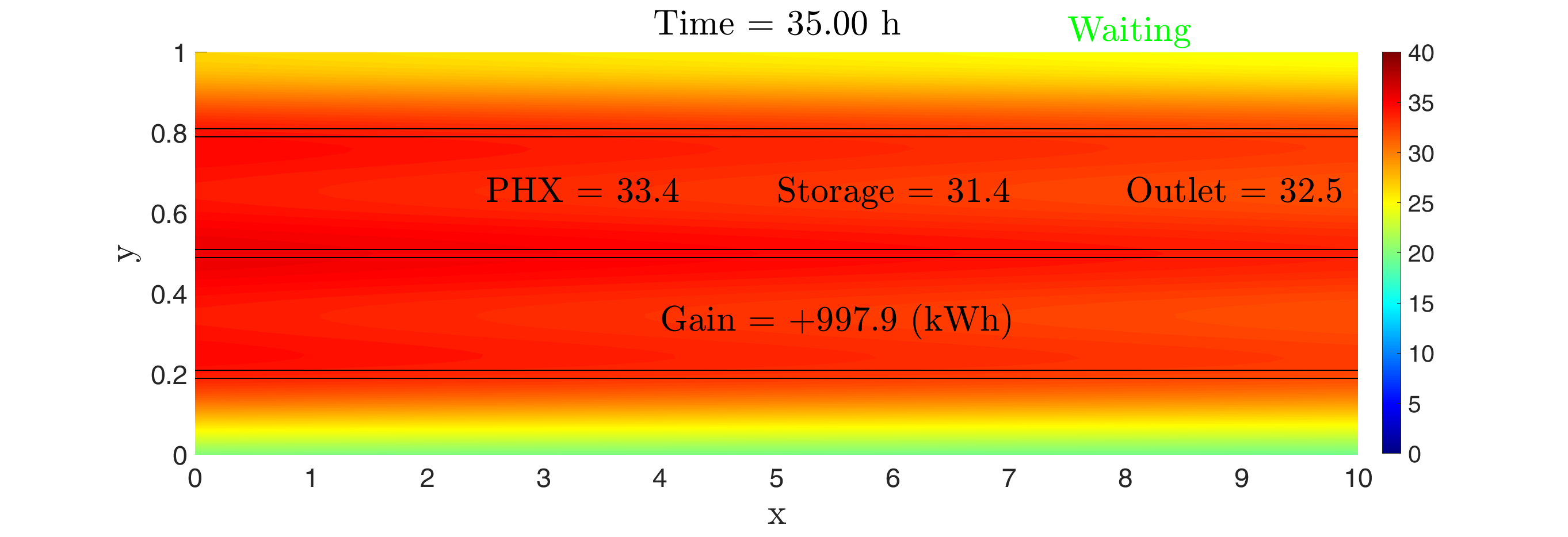}

		\hspace*{-0.01\linewidth}
		\includegraphics[width=0.49\textwidth,height=0.2\textwidth]{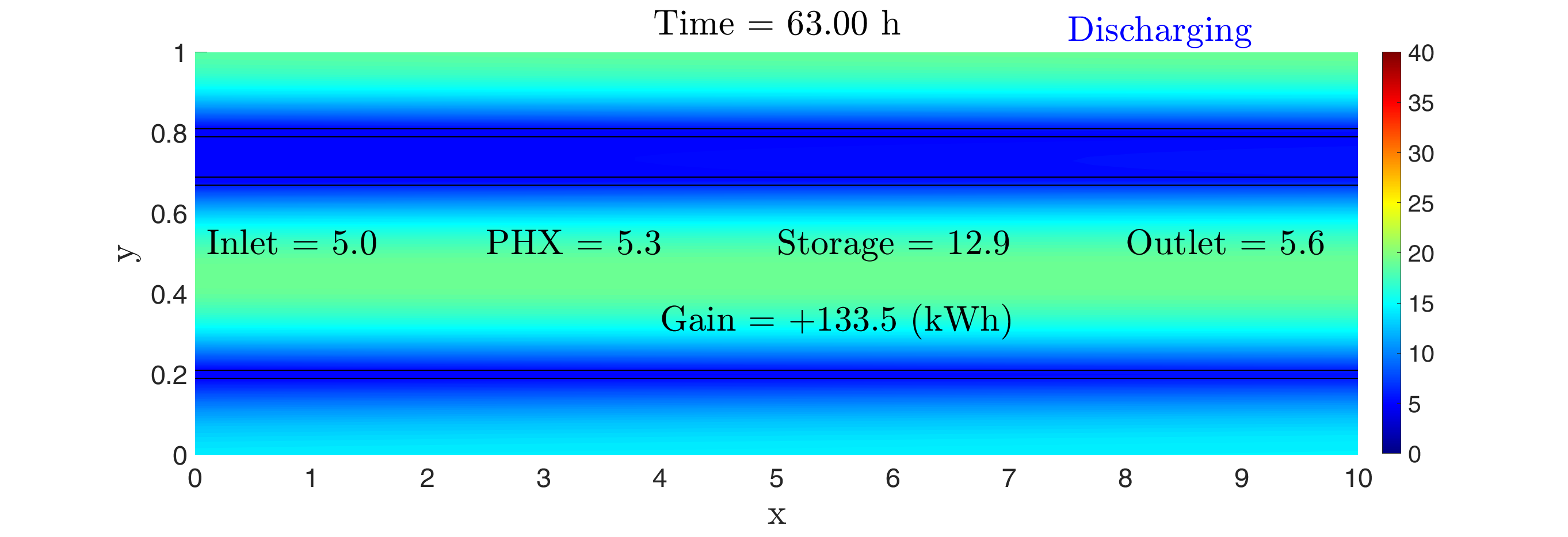}\hspace*{0.05\textwidth}
		\hspace*{-0.04\linewidth}
		\includegraphics[width=0.49\textwidth,height=0.2\textwidth]{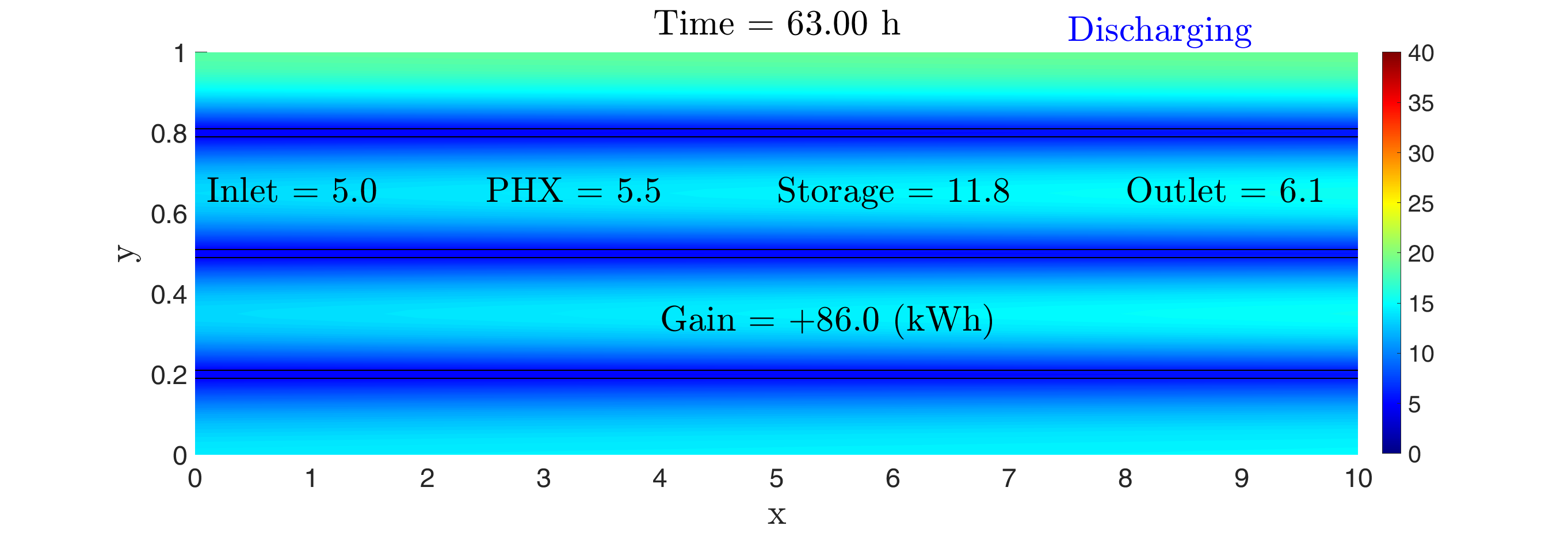} 
	\mycaption{Spatial distribution of the temperature in the storage with three horizontal \phxs
		during charging (top), waiting  (middle) and discharging (bottom) period. \newline
		Left: Non-symmetric \phxsk,~
		Right: Symmetric \phxsk.}
	\label{fig:fig5c}
\end{figure}
\begin{figure}[h]
	\centering
	\hspace*{-0.03\textwidth}	
	\includegraphics[width=0.49\textwidth,height=0.35\textwidth]{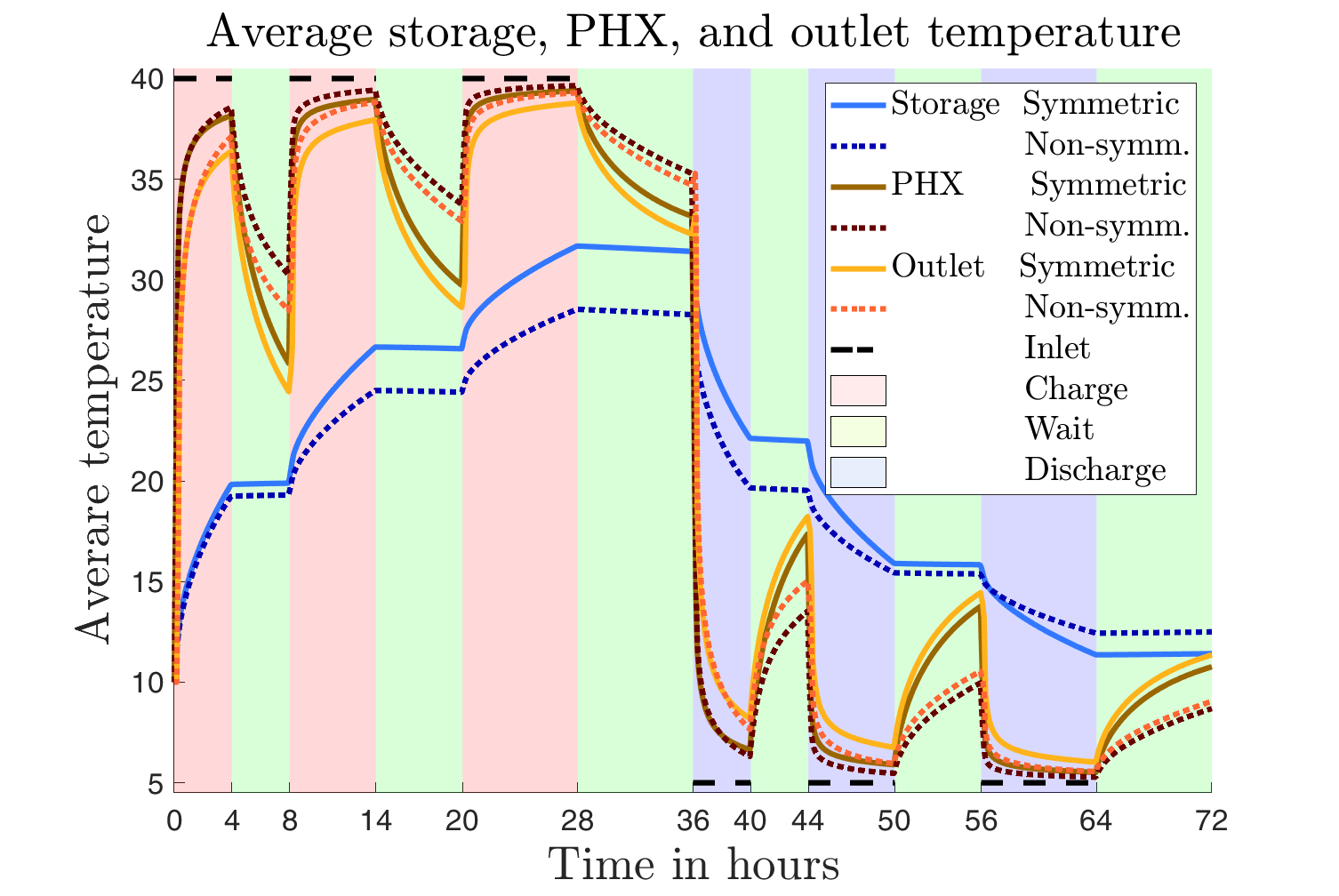}
	\hspace*{-0.01\textwidth}
	\includegraphics[width=0.49\textwidth,height=0.35\textwidth]{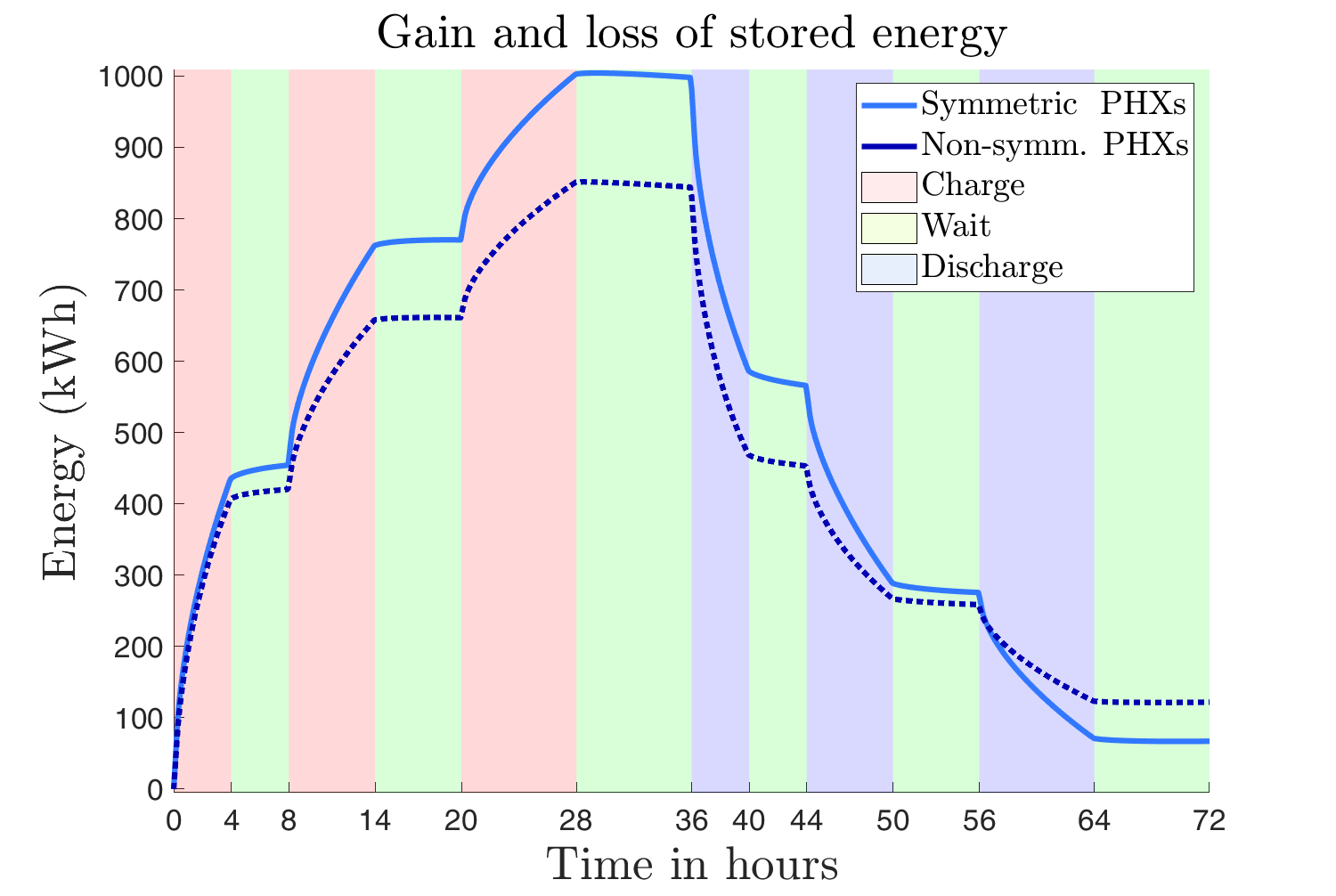}
	\mycaption{Storage with three horizontal $\phxs$	during $72$ hours with  charging,  waiting and discharging  periods. \quad 
		Left: Aggregated characteristics  $\Qmf,\Qf,\Qout$. ~~
		Right: Gain of stored energy $\Gmf$.}
	\label{fig:fig5d}
\end{figure}

In this example we add a third $\phx$ to the storage architecture and study two different $\phx$ arrangements. We proceed with the experimental design including the same waiting periods considered  in the previous subsection but now we ``glue''  together the two periods of charging and discharging each of length $36~h$. The result is a total period of length $T=72~h$ starting with a storage at temperature $Q(0,x,y)=10~\Celsius$. Within the the first 36 hours  the storage is charged  by the  moving fluid arriving at the $\phx$ inlet with temperature $\QinC(t)=40~\Celsius$.  In the second 36 hours it is  discharged  using the inlet  temperature $\QinD(t)=5~\Celsius$. The charging, waiting and discharging periods can be seen in Fig.~\ref{fig:fig5d}. Contrary to the above experiments, discharging now starts not with a temperature $35~\Celsius$ but with a non-uniformly temperature distribution which is obtained after $36\,h$ of charging (and waiting). In this more realistic setting, temperatures typically are higher in the vicinity of the $\phxs$  and lower in other regions.

Fig.~\ref{fig:fig5c} shows snapshots of the spatial temperature  distribution  during the last charging period (at $t=27 h$), during the subsequent waiting period  (at $t=35h$) and during of the last discharging period  (at $t=63h$), respectively. We compare two storage architectures with three $\phxsk$. In the first, the $\phxs$  are located symmetrically w.r.t.~the vertical mid level. For the second, the central $\phx$ was moved upwards such that we get a non-symmetric arrangement with two quite close-by $\phxs$  in the upper region. The snapshots show a strong saturation between the two upper $\phxs$  of the non-symmetric $\phx$ arrangement while for symmetric $\phxs$  the temperature distribution is much more uniform, in particular  during the waiting period as it can be seen in the middle panel for time $t=35~h$.

In Fig.~\ref{fig:fig5d} we present aggregated characteristics which are plotted over time and observe similar patterns as in the experiment with a two $\phx$ storage considered in the previous subsection.
During the waiting periods after charging the average outlet and $\phx$ temperatures decay at a faster rate  for symmetric $\phxs$   than  for non-symmetric $\phxsk$.  Vice versa they increase faster in waiting periods after discharging. This is a consequence of the stronger saturation for non-symmetric $\phxs$  which prohibits a faster cooling (warming) of the $\phx$ during waiting. For symmetric $\phxs$  the average storage temperature during charging increases faster and during discharging decreases faster than for non-symmetric $\phxsk$. This explains the similar patterns for the gain of stored energy which are plotted in the right panel. It shows that the storage with symmetric $\phxs$  (dis)charges faster than the storage with non-symmetric $\phxsk$. 

\section{Analogous Linear Time-Invariant  System}
\label{sec:analog:system}
This section is motivated by our  paper \cite{Takam2020Reduction} in which  we aim to approximate the dynamics of certain  aggregated characteristics for the infinite dimensional spatial distribution of the temperature $Q=Q(t,x,y)$  describing the storage's input-output behavior by a low-dimensional system of ODEs. Recall that the dynamics of the  spatial distribution of $Q$ is governed by the heat equation \eqref{heat_eq2}. We applied semi-discretization to that PDE  and obtained the  finite-dimensional approximation \eqref{Matrix_form1} which reads as $\frac{d Y(t)}{dt}= \mat{A}(t)Y(t)+\mat{B}(t)g(t)$ and constitutes a high-dimensional system of ODEs  for the vector function $Y$ containing the temperatures in the grid points.
In \cite{Takam2020Reduction}  that system of ODEs is the starting point for the application of model reduction techniques to find a suitable low-dimensional system of ODEs from which the aggregated characteristics can be obtained with a reasonable degree of accuracy.

Eq.~\eqref{Matrix_form1} represents a  system	of $n$  linear non-autonomous ODEs. Since some of the coefficients in the  matrices $\mat{A,B}$ resulting from the discretization of convection terms in the heat equation \eqref{heat_eq2} depend on the velocity $v_0(t)$, it follows that  $\mat{A,B}$  are time-dependent.
Thus, \eqref{Matrix_form1} does not constitute a linear time-invariant (LTI) system. The latter is a crucial assumption for most of model reduction methods such as the Lyapunov balanced truncation technique that is  considered in our paper \cite{Takam2020Reduction}. We circumvent this problem by  replacing the model for the geothermal storage by a so-called \emph{analogous model} which is LTI. 

The key idea for the construction of such an analogue is based on the observation that under the assumption of this paper our ``original model'' is already piecewise LTI. This is due to our assumption that the  fluid velocity is constant $\vconst$ during (dis)charging when the pump is on, and zero during waiting when the pump is off. This leads to the following  approximation of the original by an analogous model which is performed in two steps.

\paragraph{Approximation Step 1}
For the analogous model we assume that contrary to the original model the fluid is also moving with constant velocity $\vconst$ during pump-off periods. 
During these waiting periods   in the original model the fluid is at rest and only subject to the diffusive propagation of heat.   In order to mimic that behavior of the resting fluid by a moving fluid we assume that the temperature $\Qin$ at the \phx's inlet is equal to the average temperature of the fluid in the $\phx$ $\Qf$. From a physical point of  view we will preserve the average temperature of the fluid but a potential  temperature gradient along the $\phx$ is  not preserved  and replaced by an almost flat temperature distribution. It can be expected that the error induced by this ``mixing'' of the fluid temperature in the $\phx$ is small after sufficiently long (dis)charging periods leading to saturation with an almost constant temperature along the $\phxk$.

In the mathematical description by an initial boundary value problem for the heat equation \eqref{heat_eq2}, the above approximation leads to a modified boundary condition at the inlet. During waiting  the homogeneous Neumann boundary condition in \eqref{input} is replaced  by a non-local coupling condition such that the inlet boundary condition reads as 
\begin{align}
	Q=
	\begin{cases}
		\begin{array}{ll}	
			\Qin(t), &\text{ ~pump on,} \\
			\Qf(t), &\text{ ~pump off,} 
		\end{array}
	\end{cases} 
	\qquad  (x,y)\in 	\partial \Din.
	\label{input_analog}
\end{align}
The above condition is termed 'non-local' since the inlet temperature is not only specified by a  condition to the local temperature distribution at the inlet boundary $\partial \Din$ but it depends on the whole  spatial temperature distribution  in the fluid domain $\Df$. Semi-discretization of the above boundary condition using approximation  \eqref{eq:av:fluid} of  the average fluid temperature $\Qf=\OutputF \,Y$ formally  leads to a modification of the   input term $g(t)$ of  the system of ODEs \eqref{Matrix_form1} given in \eqref{eq:input}. That input term now reads as
\begin{align}
	g(t)=\begin{cases}
		~~(\Qin(t),~\Qg(t))^{\top}, & \quad \text{pump on},\\
		(\OutputF \,Y(t),\Qg(t))^{\top}, & \quad \text{pump off}.
	\end{cases}
	\label{eq:input_analog}
\end{align}
Further, the non-zero entries  $B_{l1}$ of the input matrix $B$ given in \eqref{eq:input_matrix} are modified. They are now no longer time-dependent but given by the constant
$		B_{l1}=\frac{\af }{h^2_x}+\frac{\vconst}{h_x} $ which was already used during pump-on periods.

\paragraph{Approximation Step 2}
From \eqref{eq:input_analog} it can be seen that the input term $g$  during pumping depends on the state vector $Y$ via $\OutputF \,Y$ and can no longer considered as exogenous.  Formally, the term  $\OutputF \,Y$ has to be included in $\mat{A}Y$ which would  lead to an additional contribution  to the system matrix  $\mat{A}$ given by $\mat{B}_{\bullet1}\OutputF$ where $\mat{B}_{\bullet1}$ denotes the first column of $\mat{B}$. Thus, the system matrix again would be time-dependent and the system not LTI.  In order to obtain an LTI system  we therefore perform a second approximation step and treat  $\Qf$ as an exogenously given quantity  (such as $\QinC,\QinD, \Qg$). This leads to a tractable approach for model reduction by the  Lyapunov balanced truncation technique applied in \cite{Takam2020Reduction}.  The latter generates low-dimensional systems depending only on the system  matrix $\mat{A}$ and the input matrix $\mat{B}$ but not on the input term $g$. Further, from an algorithmic or implementation point of view this is not a problem since given the solution $Y$ of \eqref{Matrix_form1} at time $t$, the average fluid temperature $\Qf(t)$  can be computed as a linear combination of the entries of $Y(t)$.

\begin{figure}[h]
		\centering
		\hspace*{-0.04\textwidth}
		\includegraphics[width=0.52\textwidth,height=0.45\textwidth]{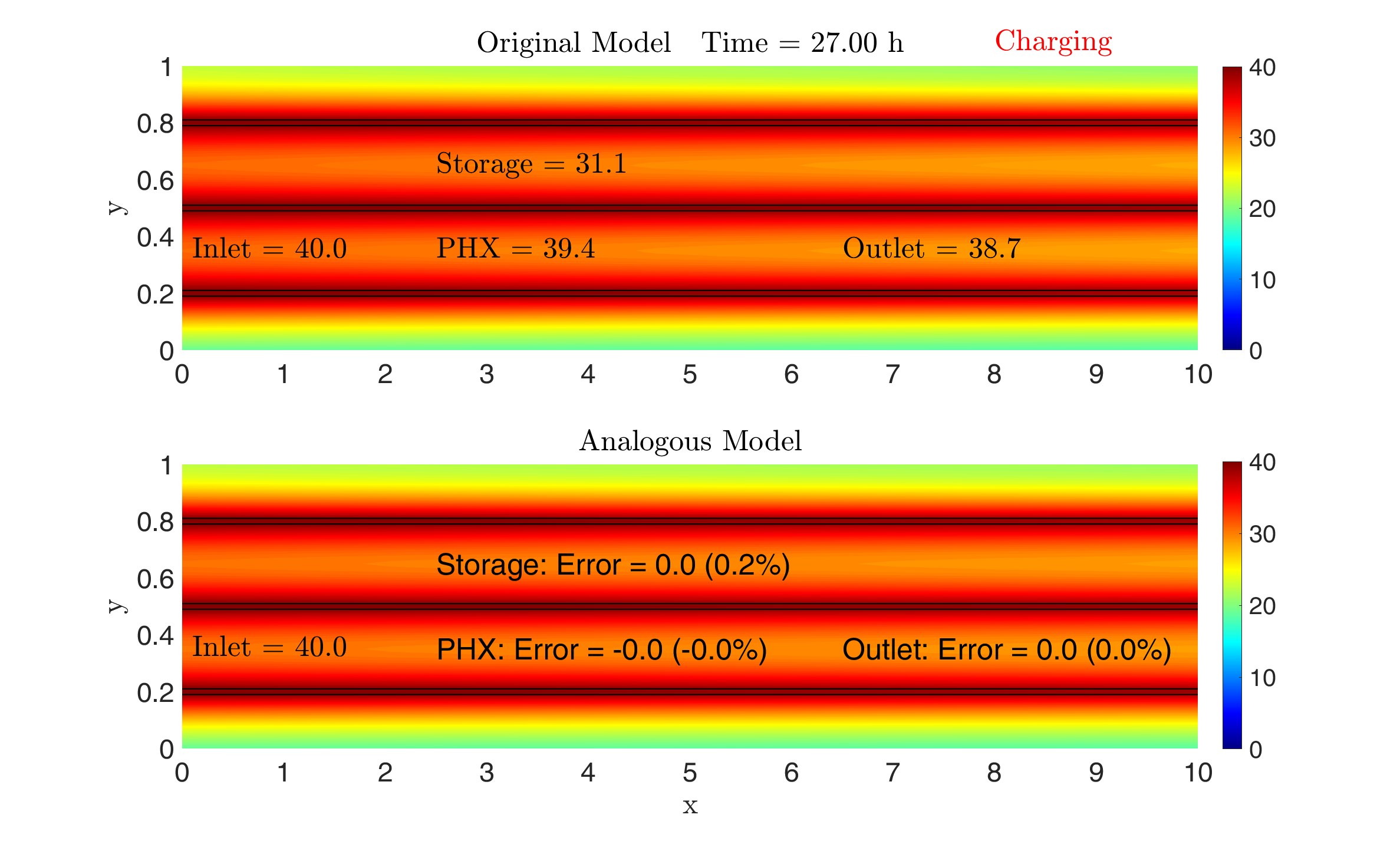}
		\hspace*{-0.05\textwidth}
		\includegraphics[width=0.52\textwidth,height=0.45\textwidth]{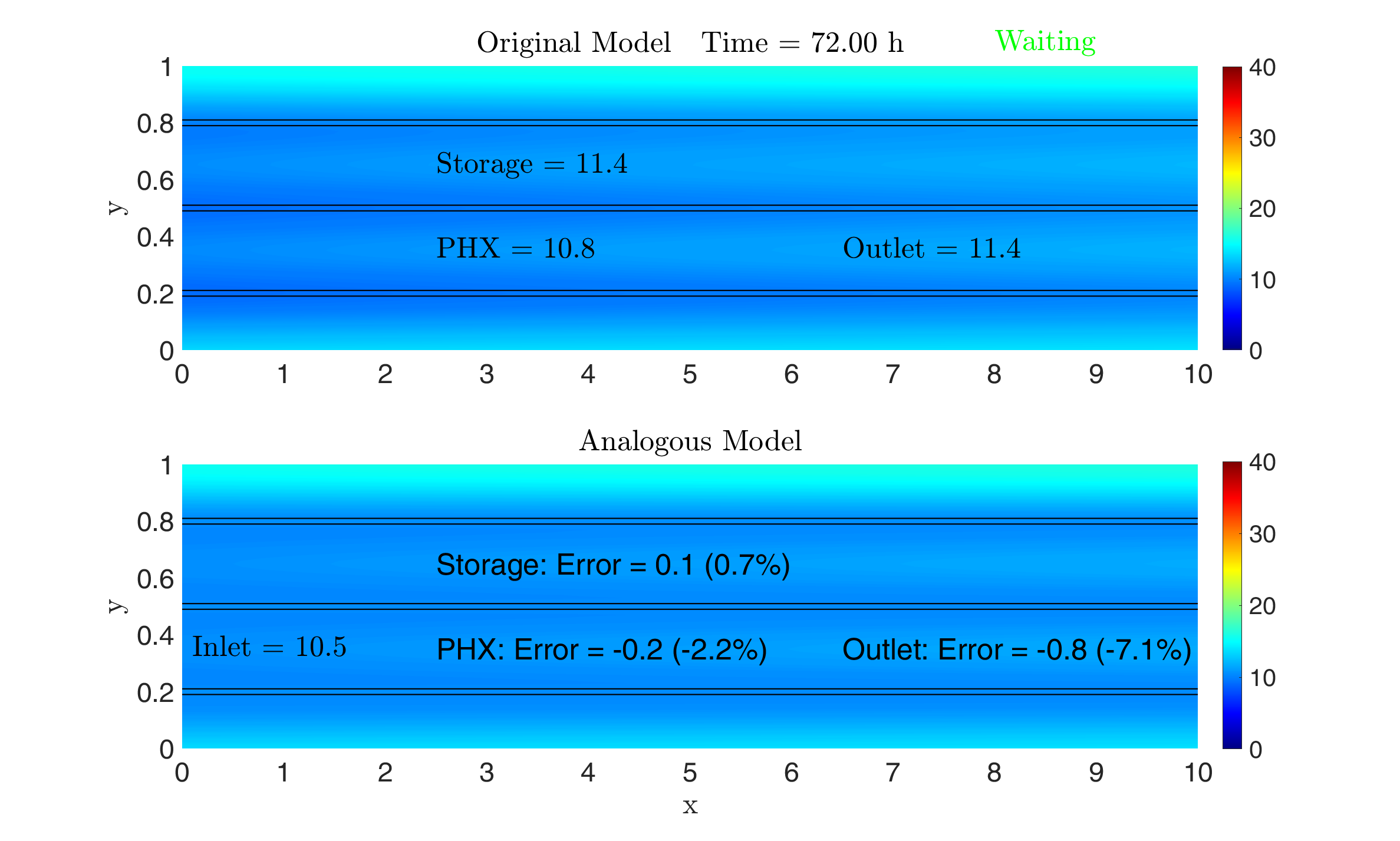}
		\mycaption{Spatial distribution of the temperature in the storage with three horizontal symmetric $\phxs$  during charging (left) and waiting (right). 
			Top: Original model.~Bottom:  Analogous model.
		}
		\label{fig:fig5e}
	\end{figure}

\begin{figure}[h]
			\centering
			\hspace*{-0.05\textwidth}
			\includegraphics[width=0.54\textwidth,height=0.4\textwidth]{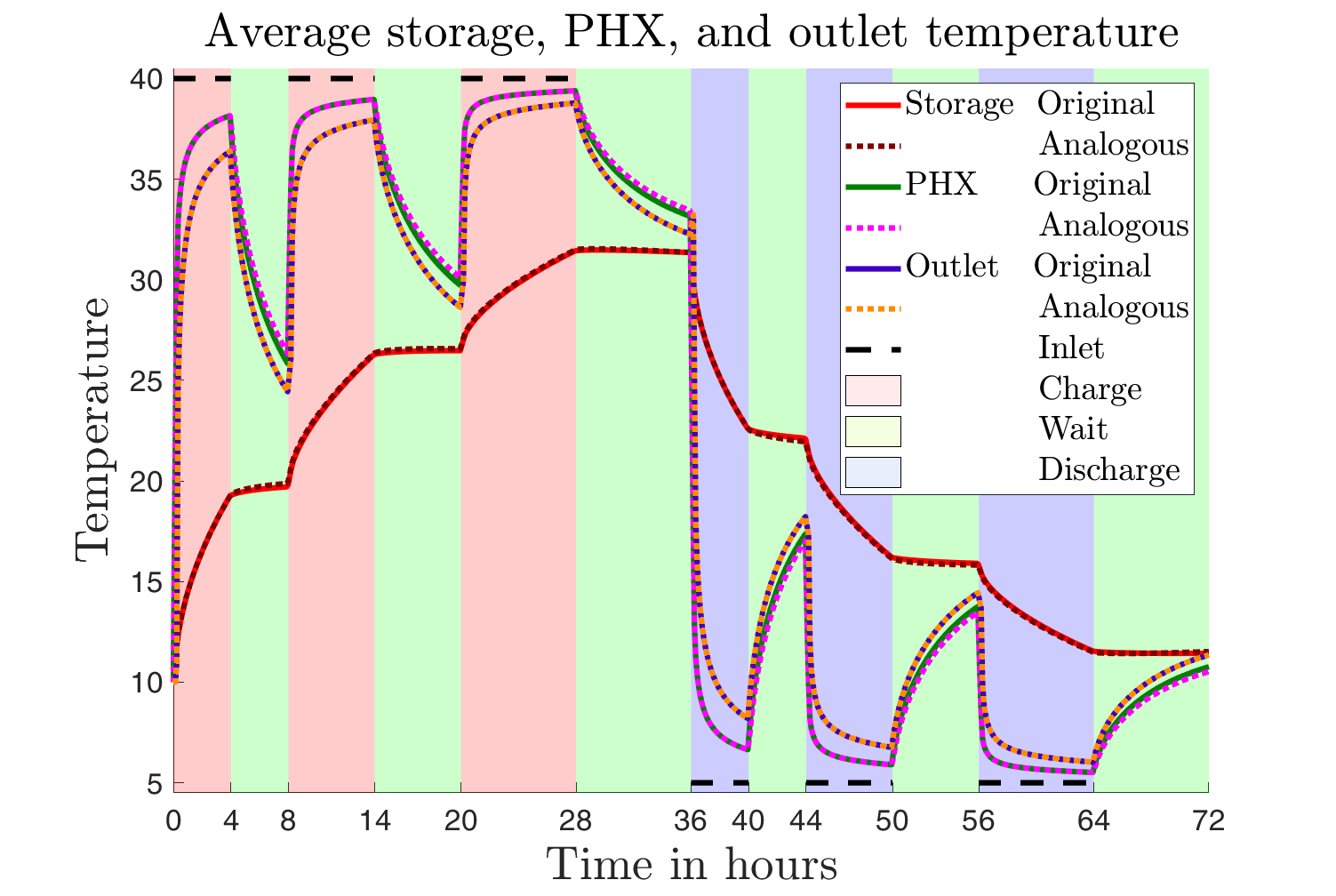}
			\hspace*{-0.05\textwidth}
			\includegraphics[width=0.54\textwidth,height=0.4\textwidth]{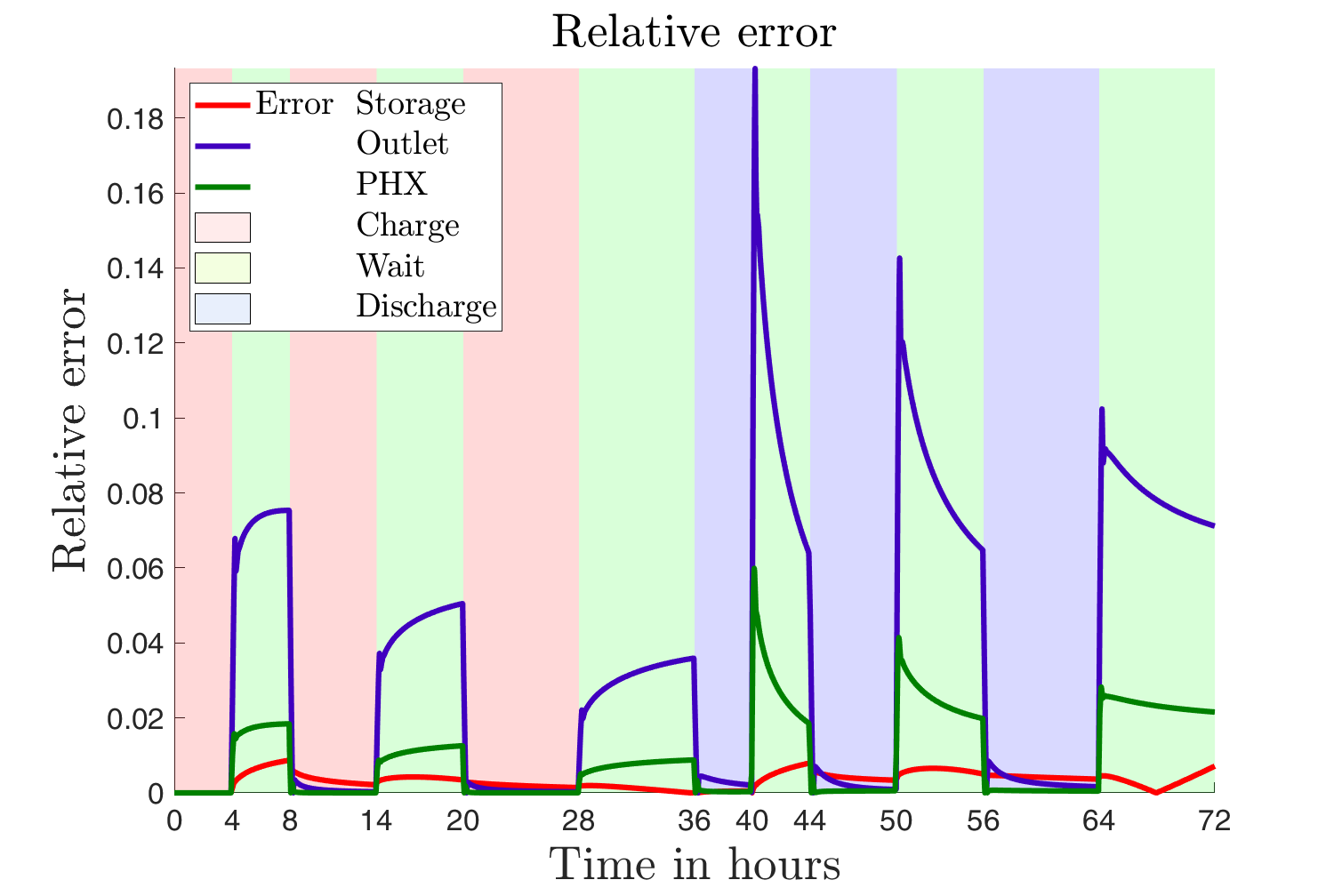}
			\mycaption{Original and analogous model of a storage with  three horizontal non-symmetric $\phxs$  during $72~h$ of charging, waiting and discharging. \newline 
				Left:~Comparison of aggregated characteristics $\Qmf,\Qf,\Qout$.\newline Right: Relative error of approximation by analogous model.}
		\label{fig:fig5f}
	\end{figure}

\paragraph{Numerical Results}
In Figs.~\ref{fig:fig5e} and \ref{fig:fig5f} we present some numerical results where we compare the spatio-temporal  temperature distribution and its aggregated characteristics of the original and the associated analogous model. These results are based on the experimental design in the Subsec.~\ref{sub: num3} for a storage architecture with three symmetric $\phxs$  and waiting periods.
Fig.~\ref{fig:fig5e} compares snapshots of the spatial temperature distribution in the storage for the original and analogous model. One snapshot is taken  during charging  and the other at the end of the last waiting period after preceding discharging periods.  At first glance there are no visible differences. A look at the aggregated characteristics  in Fig.~\ref{fig:fig5f} shows  negligible approximation errors for the average temperature in the storage $\Qmf$ and  the fluid $\Qf$. However,  the approximation of the average outlet temperature $\Qout$ suffers slightly from the replacement of a resting fluid by a moving fluid during the waiting period. The resulting ``mixing of the temperature profile''  inside the $\phx$ adjusts the outlet  to the average  in the $\phxk$. This can be seen in the right panel where the relative error for the outlet temperature dominates the errors for the two  other average temperatures in the storage and the fluid.  The experiment indicates that apart  from some noticeable approximation errors in the $\phx$ during waiting periods, in particular at the outlet, the other deviations are negligible. Finally, it can be nicely seen that during the (dis)charging periods the errors decrease and  vanish almost completely, i.e., in the long run there is no accumulation of errors.

	
	\begin{remark}\label{rem:poor_outlet_appr}
		The  poor precision of the outlet temperature approximation by the analogous model during waiting periods is of no relevance for the management and  operation of  the geothermal storage within a residential heating system. Here,  the outlet temperature is required only during charging and discharging but not during the waiting periods. The interesting quantity for which a good approximation precision is required is the average temperature in the storage and this is provided by the  analogous model.
	\end{remark}

	\section{Conclusion}
	We have investigated the numerical simulation of the short-term behavior of the spatial temperature distribution in a geothermal energy storage. The underlying initial boundary value problem for the heat equation  with a convection term  has been discretised using finite difference schemes. In a large number of numerical experiments we have shown how these simulations can support the design and operation of a geothermal storage. Examples are the dependence of the charging and discharging efficiency on the topology and arrangement  of heat exchanger $\phxs$  and on the length of charging, discharging and waiting periods.   
	
	Based on the findings of this paper we study in  \cite{Takam2020Reduction}  model reduction techniques  to derive low-dimensional approximations of aggregated characteristics of the temperature distribution  describing the input-output behavior of the storage. The latter is crucial if the geothermal storage is embedded into a  residential heating system  and the cost-optimal management of such systems is studied mathematically in terms of optimal control problems.
	
	\newpage \clearpage
	\begin{appendix}	
		\section{List of Notations}
		\label{append_a}
		\begin{longtable}{p{0.3\textwidth}p{0.68\textwidth}l}
			$Q=Q(t,x,y)$ & temperature in the geothermal storage &\\	
			$T$ & finite time horizon&\\
			$l_x$,~$l_y$,~$l_z$ &width, height and depth of the storage &\\
			$\Domainspace =(0, l_x) \times (0,l_y)$ &domain of the geothermal storage &\\
			$\Dm , ~\Df $ &   domain of  medium (soil) and  \phx fluid  &\\
			$\DInterface =\DInterface _L \cup \DInterface _U$ & interface between the $\phxs$  and the medium &\\
			$\partial \Domainspace$ &boundary of the domain&\\
			$\partial \Din$,~$\partial \Dout$ & inlet and outlet boundaries of the \phx&\\
			$\partial  \Dleft, \partial \Dright , \partial \Dtop $,~$\partial  \Dbottom $  & left, right, top and bottom boundaries of the domain&\\		
			$\mathcal{N}^*_*$  &   subsets of index pairs for grid  points  &\\
			$\mathcal{K}, \overline{\mathcal{K}}$  &   mappings $(i,j)\mapsto l$ 
			of index pairs to single indices&\\
			$v=v_0(t)(v^x, v^y)^{\top}$ & time-dependent velocity vector,&\\
			$\vconst$ &  constant velocity during pumping&\\
			$\cpf $,~$\cpm $ & specific heat capacity of the fluid  and medium &\\
			$\rhof $,~$\rhom $ & mass density of the fluid and medium&\\
			$\kappaf$,~$\kappam$ &thermal conductivity of the fluid and medium&\\
			$\af $,~$\am $& thermal diffusivity of the fluid and medium&\\
			$\heattransfer$ & heat transfer coefficient between storage and  underground &\\
			$Q_0$ &initial temperature distribution of the geothermal storage &\\
			$\Qg(t)$ & underground temperature &\\
			$\Qin(t), \QinC(t), \QinD(t)$ & inlet temperature of the $\phxk$,  during charging and discharging, &\\
			$\Qm, \Qf,\Qmf$ & average temperature in the storage medium, fluid and whole storage&\\
			$\Qout, \Qbottom$  & average temperature at the outlet   and bottom boundary&\\
			$G^*$  &  gain of thermal energy in a certain subdomain&\\
			$N_x,~N_y$,~$N_{\tau}$ & number of grid points in $x,y$ and $\tau$-direction&\\
			$h_x, h_y$,~$\tau$ & step size in $x$ and $y$-direction and the time step&\\
			$\normalvec $& outward normal to the boundary $\partial \Domainspace$&\\
			$n$ &   dimension of vector $Y$ &\\	
			$n_P$ & number of \phxs&\\	
			$\mathds{I}_n$& $n \times n$ identity matrix&\\
			$\mat{A}$ & $n \times n$ dimensional  system matrix &\\
			$\mat{B}$ & $n \times m$ dimensional input matrix    &\\
			$\mat{D}^{\pm}, ~\mat{A}_{L}, ~\mat{A}_{M},~\mat{A}_{R}$ & block matrices of matrix $\mat{A}$&\\
			$Y$&vector of temperatures at grid points&\\
			$g$& input variable of the system&\\
			$\nabla$, ~~$\Delta=\nabla \cdot \nabla$ & 	 gradient, Laplace operator&\\
			\phx &  pipe heat exchanger&\\
			LTI &  linear time invariant&
		\end{longtable}

		\section{Numerical  Computation of Aggregated Characteristics}
		\label{append:Aggregate:Num}
		\subsection{Derivation of Quadrature Formula  \eqref{quad_2D_0} }
		\label{append:quadformula}
		Rewriting  the double integral as two iterated single integrals and applying trapezoidal rule to the  outer integral we obtain (suppressing the time variable $t$)
		\begin{align*}
			J=\iint_{\mathcal{B}} Q(x,y)dxdy&= \int_{x_\iu}^{x_\io} \bigg( \int_{y_\ju}^{y_\jo}Q(x,y)\, dy\bigg)dx\\
			&\approx\int_{x_\iu}^{x_\io} h_y\bigg(\frac{1}{2}Q(x,y_\ju) +\sum_{j=\ju+1}^{\jo-1} Q(x,y_j) + \frac{1}{2}Q(x,y_\jo) \bigg) dx.
		\end{align*}
		Approximating the inner integrals  again by  trapezoidal rule we get
		\begin{align*}
			\int_{x_\iu}^{x_\io} Q(x,y_j)\, dx
			\approx h_x\bigg(\frac{1}{2}Q(x_\iu,y_j) +\sum_{i=\iu+1}^{\io-1} Q(x_i,y_j) + \frac{1}{2}Q(x_\io,y_j) \bigg), \quad j=\ju,\ldots,\jo.
		\end{align*}
		Substituting into the above expression for $J$ yields
		\begin{align*}	
			&J\approx h_xh_y \bigg(\frac{1}{4} \big[Q(x_\iu,y_\ju)+Q(x_\io,y_\ju)+Q(x_\iu,y_\jo)+Q(x_\io,y_\jo)\big]\\
			&+ \frac{1}{2} \bigg[\sum_{i=\iu+1}^{\io-1} \big[Q(x_i,y_\ju) +Q(x_i,y_\jo)\big] + \sum_{j=\ju+1}^{\jo-1} \big[Q(x_\iu,y_j) +Q(x_\io,y_j)\big]\bigg]
			+ \sum_{i=\iu+1}^{\io-1}\sum_{j=\ju+1}^{\jo-1} Q(x_i,y_j)	\bigg)\\
			&= h_xh_y \bigg(\frac{1}{4} \big[Q_{\iu\ju}(t)+Q_{\io\ju}(t)+ Q_{\iu\jo}(t)+Q_{\io\jo}(t)\big]\\
			&\phantom{=h_xh_y \bigg(} 	+ \frac{1}{2} \bigg[\sum_{i=\iu+1}^{\io-1} \big[Q_{i\ju} +Q_{i\jo}\big] + \sum_{j=\ju+1}^{\jo-1} \big[Q_{\iu j} +Q_{\io j}\big]\bigg]
			+ \sum_{i=\iu+1}^{\io-1}\sum_{j=\ju+1}^{\jo-1} Q_{ij}	\bigg).
		\end{align*}
		Since the area of the rectangle ${\mathcal{B}}$  is given by $(\io-\iu)(\jo-\ju)h_xh_y$ the average temperature $\Qav^{\mathcal{B}}$  can be approximated by
		\begin{align}
			\nonumber 
			\Qav^{\mathcal{B}} & =\frac{1}{|\mathcal{B}|}\iint_{\mathcal{B}} Q(t,x,y)dxdy \approx \sum_{(i,j)\in \mathcal{N}_\mathcal{B}} \mu_{ij}\,Q_{ij} 
		\end{align}
		with the coefficients $\mu_{ij}$ given in \eqref{mu_coeff}.
		
		\subsection{Numerical Approximation of $\Qout$ and $\Qbottom$}
		\label{append:out:bottom}
		Now we consider the average temperatures $\Qout$ and $\Qbottom$ where the temperature $Q(t,x,y)$ is averaged over one-dimensional curves  on the boundary $\partial \Domainspace$.
		Assume that $\mathcal{C}\subset \mathcal{\partial D}$ is a generic curve on one of the four outer boundaries. For the ease of exposition we restrict  $\mathcal{C} $  to be a line between the grid points  $(x_\iu,y_0)$ and $(x_\io,y_0)$ on the bottom boundary,  where $0\le \iu$, $ \iu+2\le\io\le N_x$. We denote by $\Qav^{\mathcal{C}}=\Qav^{\mathcal{C}}(t)=\frac{1}{|\mathcal{C}|}\int_{\mathcal{C}} Q(t,x,y)\,ds$ the average temperature in ${\mathcal{C}}$. Applying trapezoidal rule to the line integral  we obtain (suppressing the time variable $t$)
		\begin{align*}
			&\int_{\mathcal{C}} Q(x,y)\,ds= \int_{x_\iu}^{x_\io} Q(x,y_0)\, dx
			\approx h_x\bigg(\frac{1}{2}Q(x_\iu,y_0) +\sum_{i=\iu+1}^{\io-1} Q(x_i,y_0) + \frac{1}{2}Q(x_\io,y_0) \bigg). 
		\end{align*}
		Since the length of the curve ${\mathcal{C}}$  is given by $(\io-\iu)h_x$ the average temperature $\Qav^{\mathcal{C}}$  can be approximated by
		\begin{align}
			\label{quad_1D_0}
			\Qav^{\mathcal{C}} & =\frac{1}{|\mathcal{C}|}\int_{\mathcal{C}} Q(t,x,y)\,ds \approx \sum_{(i,j)\in \mathcal{N}_\mathcal{C}} \mu_{ij}\,Q_{ij}, 
		\end{align}
		where $\mathcal{N}_\mathcal{C}=\{ (i,j): i=\iu,\ldots, \io, j=0\}$  and the coefficients $\mu_{ij}$ of the above quadrature formula are given by
		\begin{align*}
			\mu_{ij} &= \frac{1}{(\io-\iu)} \left\{
			\begin{array}{cl@{\hspace*{1em}}l}
				1, &  \text{for }~~\iu<i<\io, ~~j=0, & \text{(inner points)}\\[1ex]
				\frac{1}{2}, & \text{for }~~i\,=\iu,\io,		&  \text{(end points).}
			\end{array} \right.
		\end{align*}
		Using the same notation and approach  as above we can rewrite approximation \eqref{quad_1D_0} as 
		\begin{align}
			\label{quad_1D_1}	\Qav^{\mathcal{C}} & \approx \sum_{(i,j)\in \mathcal{N}_\mathcal{C}^0} \mu_{ij}\,Q_{ij} +\sum_{(i,j)\in \overline{\mathcal{N}_\mathcal{C}}^0} \mu_{ij}\,Q_{ij} 
			= DY+\overline{D} \,\overline{Y},
		\end{align}
		where the matrices $D$ and $\overline{D}$ are defined as in \eqref{matrix_D} with $\mathcal{N}_\mathcal{B}^0$ and $\overline{\mathcal{N}_\mathcal{B}^0}$  replaced by  $\mathcal{N}_\mathcal{C}^0$ and $\overline{\mathcal{N}_\mathcal{C}}^0$, respectively.
		Note that in our finite difference scheme the grid values of boundary points are not contained in $Y$. Thus, we have $\mathcal{N}_\mathcal{C}^0=\varnothing$ and $D= 0_{1\times n}$.
		Finally, substituting $\overline Y=\overline C Y$ into \eqref{quad_1D_1} yields a representation of the average temperature $\Qav^{\mathcal{C}}$ as a linear combination of entries of the vector $Y$ which reads as
		\begin{align}
			\label{quad_1D_3}
			\Qav^{\mathcal{C}} & \approx    \Cav^{\mathcal{C}} \,Y \quad \text{with}\quad \Cav^{\mathcal{C}}=D+\overline{D}\,\overline{C}.
		\end{align}
		For $\mathcal{C}=\partial  \Dbottom $, i.e.,~$\iu=0,\io=N_x$ the above representation directly gives the approximation of $\Qbottom=\Cav^{\partial \Dbottom } \,Y $. For the average temperature $\Qout$ at the outlet  of a storage with $n_P$ $\phxs$  the outlet boundary $\Dout$ splits  into   $n_P$ disjoint  curves  $\Dout_j, j=1,\ldots,n_P$. Then we can apply \eqref{quad_1D_3} to derive the approximation 
		\begin{align}
			\nonumber
			\Qout=\frac{1}{|\partial \Dout|} \sum_{j=1}^{n_P} |\partial \Dout_j| \Qav^{\partial\Dout_j}
			\approx \OutputOut \,Y \quad\text{where}\quad  \OutputOut= \frac{1}{|\partial \Dout|} \sum_{j=1}^{n_P} |\partial \Dout_j| \,\Cav^{\Dout_j}.
		\end{align}

	\end{appendix}	
	
	\begin{acknowledgements}
		The authors thank  Thomas Apel (Universität der Bundeswehr München), Martin Bähr, Michael Breuss, Carsten Hartmann,  Gerd Wachsmuth (BTU Cottbus--Senftenberg), Andreas Witzig (ZHAW Winterhur), Karsten Hartig (Energie-Concept Chemnitz), Dietmar Deunert, Regina Christ (eZeit Ingenieure Berlin) for valuable discussions	that improved this paper.\\
		P.H.~Takam gratefully acknowledges the  support by the German Academic Exchange Service \linebreak[4] (DAAD) within the project ``PeStO – Perspectives in Stochastic Optimization and Applications''. \\
		R.~Wunderlich gratefully acknowledges the  support by the Federal Ministry of Education and Research (BMBF) within the project ``05M2022 - MONES: Mathematische Methoden für die Optimierung von	Nahwärmenetzen und Erdwärmespeichern''.\\
		The work of O. Menoukeu Pamen was supported with funding provided by the Alexander von Humboldt Foundation, under the programme financed by the German Federal Ministry of Education and Research entitled German Research Chair No 01DG15010.
	\end{acknowledgements}
	

\end{document}